%% file: vertex-4.tex
\documentclass[10pt]{amsart}
\usepackage{amsmath, amscd, amssymb}
\usepackage{graphicx, psfrag}
\usepackage[frame,cmtip,arrow,matrix,line,graph,curve]{xy}
\usepackage{graphpap, color}
\usepackage[mathscr]{eucal}

\numberwithin{equation}{section}

\newcommand{\A }{\mathbb{A}}
\newcommand{\CC}{\mathbb{C}}
\newcommand{\EE}{\mathbb{E}}
\newcommand{\LL}{\mathbb{L}}
\newcommand{\PP}{\mathbb{P}}
\newcommand{\QQ}{\mathbb{Q}}
\newcommand{\RR}{\mathbb{R}}
\newcommand{\ZZ}{\mathbb{Z}}

\newcommand{\bL}{\mathbf{L}}
\newcommand{\bT}{\mathbf{T}}

\newcommand{\bm}{\mathbf{m}}
\newcommand{\bn}{\mathbf{n}}
\newcommand{\bp}{\mathbf{p}}
\newcommand{\bq}{\mathbf{q}}
\newcommand{\bs}{\mathbf{s}}
\newcommand{\bt}{\mathbf{t}}
\newcommand{\bw}{\mathbf{w}}
\newcommand{\bx}{\mathbf{x}}

\newcommand{\cal}{\mathcal}

\def\cA{{\cal A}}
\def\cB{{\cal B}}
\def\cC{{\cal C}}
\def\cD{{\cal D}}
\def\cE{{\cal E}}
\def\cF{{\cal F}}
\def\cH{{\cal H}}

\def\cL{{\cal L}}
\def\cM{{\cal M}}
\def\cN{{\cal N}}
\def\cO{{\cal O}}
\def\cP{{\cal P}}
\def\cQ{{\cal Q}}
\def\cR{{\cal R}}
\def\cT{{\cal T}}
\def\cU{{\cal U}}
\def\cV{{\cal V}}
\def\cW{{\cal W}}
\def\cX{{\cal X}}
\def\cY{{\cal Y}}
\def\cZ{{\cal Z}}

\def\fD{\mathfrak{D}}

\def\fY{\mathfrak{Y}}

\def\ff{\mathfrak{f}}
\def\fm{\mathfrak{m}}
\def\fp{\mathfrak{p}}
\def\ft{\mathfrak{t}}

\def\fv{\mathfrak{v}}

\def\fl{\mathfrak{l}}

\def\frev{\mathfrak{rev}}

\newcommand{\tC}{\tilde{C}}

\newcommand{\tF}{\tilde{F}}

\newcommand{\tz}{\tilde{z}}
\newcommand{\tu}{\tilde{u}}

\def\tilcW{{\tilde\cW}}
\def\tilcM{\tilde\cM}

\newcommand{\hD}{\hat{D}}
\newcommand{\hL}{\hat{L}}

\newcommand{\hY}{\hat{Y}}

\newcommand{\hu}{\hat{u}}

\newcommand{\eGa}{\check{\Ga}}

\newcommand{\Mbar}{{\cM}}

\newcommand{\GX}{G^\bu_\chi(\Gamma,\vec{d},\vmu)}

\newcommand{\tMi}{\Mbar^{\bu\sim}_{\chi^i,\nu^i,\mu^i} }
\newcommand{\tMP}{\Mbar^{\bu\sim}_{\chi,\nu,\mu}}

\newcommand{\GYfm}{G^\bu_{\chi,\vmu}(\Gamma)}
\newcommand{\Gdmu}{G^\bu_{\chi,\vd,\vmu}(\Gamma)}
\newcommand{\Mdmu}{\cM^\bu_{\chi,\vd,\vmu}(\Yr)}

\def\Mfml{\cM_{\chi,\vd,\vmu}^\bu(\Yr,\hL)}

\def\Mmu{\cM^\bu_{\chi,\vmu}(\Gamma)}
\def\mapright#1{\,\smash{\mathop{\lra}\limits^{#1}}\,}
\def\twomapright#1{\,\smash{\mathop{-\!\!\!\lra}\limits^{#1}}\,}

\def\Mw{\cM_{\chi,\vd,\vmu}^\bu(W\urel,L)}

\def\mwtdef{\cM^\soe\ldef}
\def\MY{\cM_{\chi,\vd,\vmu}^\bu(\hat{\cY}) }

\newcommand{\three}[1]{ ({#1}_1,{#1}_2,{#1}_3) }

\newcommand{\pair}{(\vx,\vnu)}
\newcommand{\vmu}{{\vec{\mu}}}
\newcommand{\vnu}{{\vec{\nu}}}
\newcommand{\vx}{\vec{\chi}}

\newcommand{\vd}{\vec{d}}
\newcommand{\vsi}{\vec{\sigma}}
\newcommand{\vn}{\vec{n}}

\newcommand{\up}[1]{ {{#1}^1,{#1}^2,{#1}^3} }

\def\dual{^{\vee}}

\def\sta{^\ast}

\def\virt{^{\mathrm{vir}}}
\def\upmo{^{-1}}
\def\sta{^{\ast}}
\def\ori{^{\mathrm{o}}}
\def\urel{^{\mathrm{rel}}}
\def\pri{^{\prime}}
\def\virtt{^{\mathrm{vir},\te}}

\def\sta{^*}

\newcommand{\lo}[1]{ {{#1}_1,{#1}_2,{#1}_3} }

\newcommand{\xn}[1]{{#1}_{\vx,\vnu}}

\def\lra{\longrightarrow}

\def\lsta{_{\ast}}

\def\Oplus{\mathop{\oplus}}
\def\las{_{a\ast}}
\def\ldef{_{\mathrm{def}}}

\newcommand{\Si}{\Sigma}
\newcommand{\Ga}{\Gamma}

\newcommand{\ep}{\epsilon}
\newcommand{\lam}{\lambda}
\newcommand{\si}{\sigma}

\newcommand{\xnm}{ {-\chi^i +\ell(\nu^i)+ \ell(\mu^i)}}

\newcommand{\pa}{\partial}
\newcommand{\bu}{\bullet}

\newcommand{\TY}{u^*\bl\Omega_{Y_\bm}(\log \hD_\bm)\br^\vee }
\newcommand{\zero}{H_{\mathrm{et}}^0(\mathbf{R}_l^{v\bu})}
\newcommand{\one}{H_{\mathrm{et}}^1(\mathbf{R}_l^{v\bu})}

\newcommand{\ee}{{\bar{e}}}

\def\begeq{\begin{equation}}
\def\endeq{\end{equation}}
\def\and{\quad{\rm and}\quad}
\def\bl{\bigl(}
\def\br{\bigr)}

\def\mh{\!:\!}
\def\sub{\subset}
\def\lab#1{\label{#1}[{#1}]\  }
\def\Ao{{\mathbb A}^{\!1}}

\def\Zt{\ZZ^{\oplus 2}}
\def\hatYrel{\hat Y_\Gamma\urel}

\def\Po{{\mathbb P^1}}
\def\and{\quad\text{and}\quad}
\def\mapright#1{\,\smash{\mathop{\lra}\limits^{#1}}\,}
\def\lmapright#1{\,\smash{\mathop{-\!\!\!\lra}\limits^{#1}}\,}

\def\pkt{\phi_{k,t}}
\def\wm{W[\bm]}
\def\dm{D[\bm]}

\def\soe{{T_\eta}}

\def\te{T}
\def\tz{T}

\def\bbl{\Bigl(}
\def\bbr{\Bigr)}

 \DeclareMathOperator{\stab}{stab}
\DeclareMathOperator{\Ext}{Ext} 
\DeclareMathOperator{\ext}{\cE \it{xt}}
\DeclareMathOperator{\Hom}{Hom} 
\DeclareMathOperator{\Aut}{Aut} \DeclareMathOperator{\image}{Im}
 \DeclareMathOperator{\id}{id}

\DeclareMathOperator{\rank}{rank} \DeclareMathOperator{\spec}{Spec}

\def\Mdy{\cM^\bu(\hat Y)}
\def\Mdw{\cM^\bu(W)}
\def\Mdyy{\cM^\bu(Y)}
\def\mwtdef{\Xi(\eta)}
\def\tmt{(\tilde M_\infty)^T}

\newcommand{\Yr}{\hat{Y}\urel}

\newtheorem{prop}{Proposition}[section]
\newtheorem{theo}[prop]{Theorem}
\newtheorem{lemm}[prop]{Lemma}
\newtheorem{coro}[prop]{Corollary}
\newtheorem{rema}[prop]{Remark}

\newtheorem{defi}[prop]{Definition}
\newtheorem{conj}[prop]{Conjecture}

\def\Gk{G_{(k_i)}}
\def\lalp{_\alpha}
\def\lzalp{_{0,\alpha}}
\def\bone{{\mathbf 1}}

\let\lab=\label

\begin{document}

\title{A mathematical theory of the topological vertex}
\author{Jun Li}
\address{Department of Mathematics, Stanford University,
Stanford, CA 94305, USA} \email{jli@math.stanford.edu}
\author{Chiu-Chu Melissa Liu}
\address{Department of Mathematics,
Northwestern University, Evanston, IL 60208, USA}
\email{ccliu@math.northwestern.edu}
\author{Kefeng Liu}
\address{Center of Mathematical Sciences, Zhejiang University, Hangzhou,
China;
Department of Mathematics,
University of California at Los Angeles, Los Angeles, CA 90095-1555, USA }
\email{liu@cms.zju.edu.cn, liu@math.ucla.edu }
\author{Jian Zhou}
\address{Department of Mathematical Sciences,
Tsinghua University, Beijing, 100084, China}
\email{jzhou@math.tsinghua.edu.cn}

\begin{abstract}
We have developed a mathematical theory of the topological
vertex---a theory that was originally proposed by
M. Aganagic, A. Klemm, M. Mari\~{n}o, and C. Vafa
on effectively computing Gromov-Witten invariants of smooth toric Calabi-Yau
threefolds derived from duality between open string theory of
smooth Calabi-Yau threefolds and Chern-Simons theory on three manifolds.
\end{abstract}
\maketitle

\input sec1.tex

\input sec2.tex

\input sec3.tex

\input sec4.tex

\input sec5.tex

\input sec6.tex

\input sec7.tex

\input sec8.tex

\end{document}

%% file: sec1.tex
\section{Introduction}\label{sec:introduction}
In \cite{AKMV}, M. Aganagic, A. Klemm, M. Mari\~{n}o and C. Vafa
proposed a theory on computing Gromov-Witten invariants in all
genera of any smooth toric Calabi-Yau threefold derived from duality
between open string theory of smooth Calabi-Yau threefolds and
Chern-Simons theory on 3-manifolds. In summary their theory
says that
\begin{enumerate}

\item[O1.]  There exist certain open Gromov-Witten invariants
       counting holomorphic maps from bordered Riemann surfaces
       to $\CC^3$ with boundary mapped to three specific Lagrangian
       submanifolds. The {\em  topological vertex}
       $$C_{\vmu}(\lam;\bn)$$
       is a generating function of such invariants; it depends
       on a triple of partitions $\vmu=(\up{\mu})$ and a triple of integers
       $\bn=\three{n}$.

\item[O2.] The Gromov-Witten invariants of any smooth toric Calabi-Yau
      threefold can be expressed in terms of
      $C_{\vmu}(\lam;\bn)$ by explicit gluing algorithms.

\item[O3.] By the duality between Chern-Simons theory and Gromov-Witten
      theory, the topological vertex is given by
\begin{equation}\label{eqn:O3}
C_{\vmu}(\lam;\bn)=q^{\frac{1}{2}(\sum_{i=1}^3\kappa_{\nu^i}n_i)}
\cW_\vmu(q),
\end{equation}
where $q=e^{\sqrt{-1}\lam}$ and $\cW_\vmu(q)$ is a combinatorial
expression related to the Chern-Simons link invariants of a
particular link. (cf. Section \ref{sec:W}.)
\end{enumerate}

As was demonstrated in the work of many, for instance
\cite{Pen, Kon1, Kon2}, this algorithm is extremely efficient in deriving the
structure result of the Gromov-Witten invariants of toric
Calabi-Yau threefolds.

The purpose of this paper is to provide a mathematical theory
for this algorithm. To achieve this, we need to provide a
mathematical definition of the open Gromov-Witten invariants
referred to in O1; secondly, we need to establish the gluing
algorithms O2; finally, the duality O3 must be established
mathematically.

In this paper, we shall complete the first two steps as
outlined. Our theory is based on relative Gromov-Witten theory
\cite{Li-Rua, Ion-Par1, Ion-Par2, Li1, Li2}. Our results can be
summarized as follows.
\begin{enumerate}
\item[R1.] We introduce the notion of formal toric Calabi-Yau (FTCY) graphs,
      which is a refinement and generalization of the graph associated to a
      toric Calabi-Yau threefold. An FTCY graph $\Gamma$ determines a
      relative FTCY threefold $Y\urel=(\hY,\hD)$; it can be degenerated to
      indecomposable ones.

\item[R2.] We define {\em formal relative Gromov-Witten invariants}
      for relative
      FTCY threefolds ({\bf Theorem \ref{thm:R2}}). These invariants include
      as special cases Gromov-Witten invariants of smooth toric Calabi-Yau threefolds.

\item[R3.] We show that the formal relative Gromov-Witten invariants
      in R2 satisfy
      the degeneration formula of relative Gromov-Witten invariants
      of projective varieties ({\bf Theorem \ref{thm:R3}}).

\item[R4.] Any smooth relative FTCY threefold can
      be degenerated to indecomposable ones, whose isomorphism classes
      are determined by a triple of integers $\bn=\three{n}$.
      By degeneration formula,
      the formal relative Gromov-Witten invariants in R2 can be expressed in
      terms of the generating function $\tilde{C}_{\vmu}(\lam;\bn)$ of
      that of an indecomposable FTCY threefold
      ({\bf Proposition \ref{thm:R4}}).
      This degeneration formula coincides with the gluing algorithms described in O2.

\item[R5.]  We evaluate $\tC_{\vmu}(\lam;\bn)$
({\bf Proposition \ref{thm:tCframe}}, {\bf Theorem \ref{thm:R5}}):
\begin{equation}\label{eqn:R5}
\tC_{\vmu}(\lam;\bn)=
q^{\frac{1}{2}(\sum_{i=1}^3\kappa_{\nu^i}n_i )} \tilcW_\vmu(q),
\end{equation}
where $\tilcW_\vmu(q)$ is a combinatorial expression  defined by \eqref{eqn:tW}
in Section \ref{sec:W}.
\end{enumerate}

\smallskip

In R4, we shall define $\tilde{C}_\vmu(\lam;\bn)$ as local
relative Gromov-Witten invariants of a formal Calabi-Yau $(\hat
Z,\hat D)$ that is the infinitesimal neighborhood of a configuration
$C_1\cup C_2\cup C_3$ of three $\PP^1$'s meeting at a point $p_0$ in
a relative Calabi-Yau threefold $(Z,D)$;
the stable maps have ramification partition $\mu^i$ around the
relative divisor $D$. Since $\hat Z$ is formal, we shall define the
local invariants $\tilde{C}_{\vmu}(\lam;\bn)$ via localization
formula. The most technical part of the paper is to show that such
local invariants exist as {\em topological} invariants; namely they
are rational numbers independent of equivariant parameters ({\bf
Theorem \ref{thm:tri-inv}}, invariance of the topological vertex).

Our results R1-R5, together with a conjectural identity
$\tilcW_\vmu(q)=\cW_\vmu(q)$ ({\bf Conjecture \ref{conj:WW}}), will
provide a complete mathematical theory of the topological vertex
theory. The conjecture holds when one of the partitions, say
$\mu^3$, is empty ({\bf Corollary \ref{thm:R5-two}}); it also holds
for all low degree cases we have checked.

By virtual localization, $\tC_{\vmu}(\lam;\bn)$ can be expressed in
terms of Hodge integrals ({\bf Proposition \ref{thm:GtC}}),
which combined with (\ref{eqn:R5}) provides us a formula of a
generating function $G^\bu_\vmu(\lam;\bw)$ of the three-partition
Hodge integrals ({\bf Theorem \ref{thm:Gthree}}):
\begin{equation}\label{eqn:3partition}
G^\bu_\vmu(\lam;\bw)=\sum_{|\nu^i|=|\mu^i|}\prod_{i=1}^3\frac{\chi_{\nu^i}(\mu^i)}{z_{\mu^i}}
q^{\frac{1}{2}(\sum_{i=1}^3 \kappa_{\nu^i}\frac{w_{i+1}}{w_i})}
\tilde{\cW}_\vnu(q).
\end{equation}
(See Section \ref{sec:pre} for precise definitions involved in
\eqref{eqn:3partition}). This generalizes a formula of two-partition
Hodge integrals ({\bf Theorem \ref{thm:Gtwo}}) proved in
\cite{LLZ2}. Our evaluation of three-partition Hodge integrals,
thus $\tC_{\vmu}(\lam;\bn)$, relies on the cut-and-join equations
of three partition Hodge integrals, which turn out to be
equivalent to the invariance of the topological vertex
mentioned above. Indeed, the cut-and-join equations
are so strong that they reduce the evaluation of general
three-partition Hodge integrals first to the case where $\mu^3=\emptyset$,
and eventually to the case where $\mu^1=(d)$ and $\mu^2=\mu^3=\emptyset$.

An important class of toric Calabi-Yau threefolds consists of local
toric surfaces in a Calabi-Yau threefold; they are the total space
of the canonical line bundle of a projective toric surface (e.g.
$\cO_{\PP^2}(-3)$). In this case, only
$\tC_{\mu^1,\mu^2,\emptyset}(\lam;\bn)$ (or two-partition Hodge
integrals) are involved. The algorithm in this case was described in
\cite{Aga-Mar-Vaf}; explicit formula was given in \cite{Iqb} and
derived in \cite{Zho3} by localization, using the formula of
two-partition Hodge integrals.

It is worth mentioning that,
assuming the existence of $C_\vmu(\lam;\bn)$ and the validity of
open string virtual localization,  Diaconescu and  Florea related
$C_\vmu(\lam;\lo{n})$ (at certain fractional $n_i$) to
three-partition Hodge integrals, and derived the gluing algorithms
in O2 by localization \cite{Dia-Flo}.

Maulik, Nekrasov, Okounkov, and Pandharipande conjectured a
correspondence between the Gromov-Witten
and Donaldson-Thomas theories for any non-singular
projective threefold \cite{MNOP1, MNOP2}.
This correspondence can also be formulated for
certain non-compact threefolds in the presence
of a torus action; the correspondence for toric Calabi-Yau threefolds is
equivalent to the validity of the topological vertex \cite{MNOP1, Oko-Res-Vaf}.
For non-Calabi-Yau toric threefolds the building block is the equivariant vertex
(see \cite{MNOP1, MNOP2, Pan-Tho1, Pan-Tho2}) which depends on equivariant 
parameters\footnote{During the revision of this paper, Maulik, Oblomkov, Okounkov,
and Pandharipande announced a proof of GW/DT correspondence for all toric threefolds \cite{MOOP}.}.

The rest of this paper is organized as follows. In Section
\ref{sec:pre}, we recall some definitions and previous results, and
introduce some generating functions. R1 is carried out in Section
\ref{sec:FTCY}. R2 is carried out in Section \ref{sec:invariants};
the case when the relative FTCY threefold is indecomposable gives
the mathematical definition of topological vertex, and the proof of
its invariance (Theorem \ref{thm:tri-inv}) is given in Section
\ref{sec:tri-inv}. In Section \ref{sec:hodge}, we express the
topological vertex in terms of three-partition Hodge integrals and
double Hurwitz numbers. In Section \ref{sec:glue}, we establish R3
and R4. In Section \ref{sec:comb}, we derive the combinatorial
expression in R5. Some examples of the identity
$\cW_\vmu(q)=\tilcW_\vmu(q)$ are listed in Section
\ref{sec:lowdegree}.

\bigskip

\paragraph{{\bf Acknowledgments.}}
We wish to thank A. Aganagic, A. Klemm, M. Mari\~{n}o, and C. Vafa
for their explanation of \cite{AKMV} in public lectures and private
conversations. We also wish to thank S. Katz and E. Diaconescu for
helpful conversations. We thank A. Klemm for writing a maple program
which computes both $\cW_\vmu(q)$ and $\tilcW_\vmu(q)$. The first
author is supported by NSF grant DMS-0200477 and DMS-0244550.
The third author is supported by the NSF and the Guggenheim
foundation.

%% file: sec2.tex
\section{Definitions and Previous Results} \label{sec:pre}
In this section, we recall some definitions and
previous results, and introduce some generating functions. We begin
with the partitions and representations of symmetric groups.

\subsection{Partitions and Representations of Symmetric Groups}\label{sec:W}
Recall that a \emph{partition} $\mu$ of a nonnegative integer $d$ is
a sequence of positive integers
$$
\mu=(\mu_1\geq \mu_2 \geq \cdots \geq \mu_h>0)\quad\text{such that}\
\ d=\mu_1+\ldots+\mu_h.
$$
We write $\mu\vdash d$ or $|\mu|=d$, and call $\ell(\mu)=h$ the
\emph{length} of the partition. For convenience, we denote by
$\emptyset$ the empty partition; thus $|\emptyset|=
\ell(\emptyset)=0$. The order of $\Aut(\mu)$, the group of
permutations of $\mu_1,\mu_2,\cdots$ that leave $\mu$ fixed, is
$$ 
|\Aut(\mu)|=\prod_{j}m_j(\mu)!,\quad\text{where}\quad
m_j(\mu)=\#\{i:\mu_i=j\}.
$$
The \emph{transpose} of $\mu$ is a partition $\mu^t$ defined by $
m_i(\mu^t)=\mu_i-\mu_{i+1}$. Note that
$$
|\mu^t|=|\mu|,\ \, (\mu^t)^t=\mu,\ \ \ell(\mu^t)=\mu_1.
$$

A partition $\mu$ also corresponds to a conjugacy class in $S_d$,
where $S_d$ is the permutation group of $d=|\mu|$ elements in the
obvious way. With this understanding, the cardinality $z_\mu$ of the
centralizer of any element in this conjugacy class is of the form,
$$
z_\mu=a_\mu|\Aut(\mu)|,\quad \text{where}\ a_\mu=\mu_1\cdots
\mu_{\ell(\mu)}.
$$

We let $\cP$ denote the set of partitions; we let
$$
\cP_+=\cP-\{\emptyset\},\ \
\cP^2_+=\cP^2-\{(\emptyset,\emptyset)\},\ \
\cP^3_+=\cP^3-\{(\emptyset,\emptyset,\emptyset)\}.
$$
Given a triple of partitions $\vmu=(\up{\mu})\in\cP^3$, we define
$$
\ell(\vmu)=\sum_{i=1}^3\ell(\mu^i),\ \
\Aut(\vmu)=\prod_{i=1}^3\Aut(\mu^i).
$$

For any partition $\nu$, we let $\chi_\nu$ denote the irreducible
character of $S_{|\nu|}$ indexed by $\nu$, and let $\chi_\nu(\mu)$
be the value of $\chi_\nu$ on the conjugacy class determined by the
partition $\mu$. Recall that the Schur functions $s_\mu$ are related to
the Newton functions $p_i(x)=x_1^i+x_2^i+\cdots$ by
$$
s_\mu(x)=\sum_{|\nu|=|\mu|}\frac{\chi_\mu(\nu)}{z_\nu}p_\nu(x),
\quad x=(x_1,x_2,\ldots).
$$

The Littlewood-Richardson coefficients $c_{\mu\nu}^\eta$, which are
nonnegative integers, and skew Schur functions are determined
according to the rules
$$
s_\mu s_\nu =\sum_{\eta} c^\eta_{\mu\nu}s_\eta\and
s_{\eta/\mu}=\sum_\nu c^\eta_{\mu\nu}s_\nu.
$$

In order to define the combinatorial expressions $\cW_\vmu(q)$ and
$\tilde{\cW}_\vmu(q)$ in O3 and R5 in the introduction (Section
\ref{sec:introduction}), we need to introduce more notation. We define
$[m] = q^{m/2} - q^{-m/2}$, and define
\begin{equation}\label{eqn:kappa}
\kappa_\mu=\sum_{i=1}^{\ell(\mu)}\mu_i(\mu_i-2i+1),
\end{equation}
which for transpose partitions satisfies $\kappa_{\mu^t}=-\kappa_\mu$.

We define
\begin{equation}\label{eqn:Wone}
\cW_{\mu}(q) = q^{\kappa_{\mu}/4}\prod_{1 \leq i < j \leq \ell(\mu)}
\frac{[\mu_i - \mu_j + j - i]}{[j-i]} \prod_{i=1}^{\ell(\mu)}
\prod_{v=1}^{\mu_i} \frac{1}{[v-i+\ell(\mu)]}
\end{equation}
and define
\begin{equation}\label{eqn:Wtwo}
\cW_{\mu, \nu}(q) = q^{|\nu|/2} \cW_{\mu}(q) \cdot
s_{\nu}(\cE_{\mu}(q,t)),
\end{equation}
where
$$
\cE_{\mu}(q,t) = \prod_{j=1}^{l(\mu)}
\frac{1+q^{\mu_j-j}t}{1+q^{-j}t} \cdot \biggl(1 +
\sum_{n=1}^{\infty} \frac{t^n}{\prod_{i=1}^n (q^i-1)}\biggr).
$$
We also denote
$$c_{\rho^1(\rho^3)^t}^{\mu^1(\mu^3)^t} = \sum_{\eta}
c_{\eta\rho^1}^{\mu^1}c_{\eta(\rho^3)^t}^{(\mu^3)^t}.
$$

\begin{defi}  For $\vmu=(\up{\mu})$,
we define
\begin{equation}\label{eqn:Wthree}
\cW_{\vmu}(q) =
q^{\kappa_{\mu^2}/2+\kappa_{\mu^3}/2}\sum_{\rho^1,
\rho^3}c_{\rho^1(\rho^3)^t}^{\mu^1(\mu^3)^t}
\frac{\cW_{(\mu^2)^t\rho^1}(q)\cW_{\mu^2(\rho^3)^t}(q)}{\cW_{\mu^2}(q)}\, .
\end{equation}
\end{defi}

\begin{defi}
For a partition $\mu=(\mu_1\geq\mu_2\geq \cdots)$ we let $2\mu$ be the
partition $(2\mu_1\geq 2\mu_2\geq \cdots)$. For $\vec{\rho}=(\up{\rho})$, we
define
\begin{equation}\label{eqn:tW}
\begin{aligned}
\tilde{\cW}_{\vec{\rho}}(q) =&
q^{-(\kappa_{\rho^1}-2\kappa_{\rho^2}-\frac{1}{2}\kappa_{\rho^3})/2}
 \sum_{\nu^+, \eta^1,\nu^1,\eta^3,\nu^3} c^{\nu^+}_{(\nu^1)^t\rho^2}
c_{(\eta^1)^t\nu^1}^{\rho^1}c_{\eta^3(\nu^3)^t}^{\rho^3}\\
&\cdot q^{(-2\kappa_{\nu^+}  - \frac{\kappa_{\nu^3}}{2})/2}
\cW_{\nu^+, \nu^3}(q) \sum_\mu
\frac{1}{z_\mu}\chi_{\eta^1}(\mu)\chi_{\eta^3}(2\mu)\, .
\end{aligned}
\end{equation}
\end{defi}

We have the following identities (see \cite{Zho4}):
\begin{equation}\label{threetwo}
\cW_{\mu,\nu,\emptyset}(q)= \cW_{\emptyset,\mu,\nu}(q)=
\cW_{\nu,\emptyset,\mu}(q)=q^{\kappa_\nu/2}\cW_{\mu,(\nu)^t}(q).
\end{equation}
\begin{equation}
\cW_{\mu,\nu}(q)=\cW_{\nu,\mu}(q),\quad
\cW_{\mu,\emptyset}(q)=\cW_\mu(q).
\end{equation}

\subsection{Double Hurwitz Numbers}\label{sec:H}
We now come to the generating function of the double Hurwitz
numbers. Let $\mu^+,\mu^-$ be partitions of $d$; let
$H^\bu_{\chi,\mu^+,\mu^-}$ be the weighted counts of Hurwitz covers
of the sphere of the type $(\mu^+,\mu^-)$ by possibly disconnected
Riemann surfaces of Euler characteristic $\chi$. We form the
generating function
$$
\Phi^\bu_{\mu^+,\mu^-}(\lam)= \sum_\chi
\lam^{-\chi+\ell(\mu^+)+\ell(\mu^-)} \frac{H^\bu_{\chi,\mu^+,\mu^-}
}{(-\chi+\ell(\mu^+)+\ell(\mu^-))!}.
$$
By Burnside formula,
\begin{equation}\label{eqn:burnside}
\Phi^\bu_{\mu^+,\mu^-}(\lam)= \sum_{|\nu|=d}e^{\kappa_\nu \lam/2}
\frac{\chi_\nu(\mu^+) }{z_{\mu^+} } \frac{\chi_\nu(\mu^-)
}{z_{\mu^-} }.
\end{equation}
Using the orthogonality of characters
\begin{equation}\label{eqn:ortho}
\sum_\rho
\frac{\chi_\mu(\rho)\chi_\nu(\rho)}{z_\rho}=\delta_{\mu\nu},
\end{equation}
it is straightforward to check that (\ref{eqn:burnside}) implies the
following two identities:
\begin{equation}\label{eqn:Hsum}
\Phi^\bu_{\mu^+,\mu^-}(\lam_1+\lam_2)
=\sum_{|\nu|=d}\Phi^\bu_{\mu^+,\nu}(\lam_1)z_\nu
\Phi^\bu_{\nu,\mu^-}(\lam_2)
\end{equation}
and
\begin{equation}\label{eqn:Hinitial}
\Phi^\bu_{\mu^+,\mu^-}(0)=\frac{\delta_{\mu^+\mu^-}}{z_{\mu^+} }.
\end{equation}
Equation \eqref{eqn:Hsum} is a sum formula for double Hurwitz
numbers; Equation \eqref{eqn:Hinitial} gives the initial values for
the double Hurwitz numbers.

Following \cite{Gou, Gou-Jac, Gou-Jac-Vai}, we introduce a differential equation
that has a generating function of
$\Phi^\bu_{\mu^+,\mu^-}$ as its unique solution
which satisfies the initial condition \eqref{eqn:Hinitial}.
This equation is similar to \cite[Lemma 2.2]{Gou-Jac}
and \cite[Lemma 3.1]{Gou-Jac-Vai}.

We let $p^\pm=(p^\pm_1,p^\pm_2,\ldots)$ be formal variables, and for a
partition $\mu$ we let $p^\pm_\mu=p^\pm_{\mu_1}\cdots
p^\pm_{\mu_{\ell(\mu)} }$. We then define a generating function
$$
\Phi^\bu(\lam;p^+,p^-)=1+\sum_{d=1}^\infty
\sum_{|\mu^\pm|=d}\Phi^\bu_{\mu^+,\mu^-}(\lam)p^+_{\mu^+}
p^-_{\mu^-},
$$
and differential operators
$$
C^\pm=\sum_{j,k}(j+k)p^\pm_j p^\pm_k\frac{\pa}{\pa p^\pm_{j+k}},\quad
J^\pm=\sum_{j,k}jkp^\pm_{j+k}\frac{\pa^2}{\pa p^\pm_j\pa p^\pm_k}.
$$
They form a {\em cut-and-join equation} for double Hurwitz numbers:
\begin{equation}\label{eqn:Hcj}
\frac{\pa \Phi^\bu}{\pa\lam} =\frac{1}{2}(C^+ +
J^+)\Phi^\bu=\frac{1}{2}(C^- + J^-)\Phi^\bu.
\end{equation}
The generating function $\Phi^\bu(\lam;p^+,p^-)$ is the unique solution to this
system satisfying the initial value
$$
\Phi^\bu(0;p^+,p^-)=1+\sum_{\mu\in\cP^+}\frac{p^+_\mu
p^-_\mu}{z_\mu}.
$$

\subsection{Three-Partition Hodge Integrals}
\label{sec:G}

We shall introduce three-partition Hodge integrals in this
subsection.

For the three-partition Hodge integrals we need to work with
Deligne-Mumford moduli stack $\Mbar_{g,n}$ of stable $n$-pointed
nodal curves of genus $g$. Over this moduli stack, we let
$\pi:\Mbar_{g,n+1}\to \Mbar_{g,n}$ be the universal curve, let
$s_i:\Mbar_{g,n}\to\Mbar_{g,n+1}$ be its section that corresponds to
the $i$-th marked points of the family, and let $\omega_\pi$ be the
relative dualizing sheaf. The commonly known $\lambda$ classes and
the $\psi$-class are defined using these morphisms. One forms the
Hodge bundle $\EE=\pi_*\omega_\pi$; its $j$-th Chern class
$\lambda_j=c_j(\EE)$ is the $\lambda$-class. One then form the pull
back line bundle $\LL_i=s_i^*\omega_\pi$; its first Chern class
$\psi_i=c_1(\LL_i)$ is the $\psi$-class. A Hodge integral is then an
integral of the form
$$\int_{\Mbar_{g, n}} \psi_1^{j_1}
\cdots \psi_n^{j_n}\lambda_1^{k_1} \cdots \lambda_g^{k_g}.
$$

We now introduce three-partition Hodge integrals. Let $\lo{w}$
be formal variables. In this Subsection, and in
Section \ref{sec:hodge}, \ref{sec:glue}, \ref{sec:comb}, we use the following
convention:
\begin{equation}\label{eqn:w-convention}
\bw=(\lo{w}),\quad w_1+w_2+w_3=0,\quad w_4=w_1.
\end{equation} 
For $\vmu=(\up{\mu})\in \cP^3_+$, we let
$$
d^1_\vmu=0,\quad 
d^2_\vmu=\ell(\mu^1),\quad
d^3_\vmu=\ell(\mu^1)+\ell(\mu^2).
$$
We define the {\em three-partition Hodge integral}
to be
\begin{equation}\label{eqn:Gg}
G_{g,\vmu}(\bw) = \frac{(-\sqrt{-1})^{\ell(\vmu)} }{|\Aut(\vmu)|}
 \prod_{i=1}^3\prod_{j=1}^{\ell(\mu^i)}
\frac{\prod_{a=1}^{\mu^i_j-1}(\mu^i_j w_{i+1} + a w_i) }
     {(\mu^i_j-1)!w_i^{\mu^i_j-1} } 
\int_{\Mbar_{g,\ell(\vmu)}
} \prod_{i=1}^3\frac{\Lambda_g^\vee(w_i)w_i^{\ell(\vmu)-1} }
{\prod_{j=1}^{\ell(\mu^i)}(w_i(w_i-\mu^i_j\psi_{d^i_\vmu+j}))}
\end{equation}
where $\Lambda_g^\vee(u)=u^g-\lam_1 u^{g-1}+\cdots + (-1)^g\lam_g$.

It is clear from the definition that 
$G_{g,\vmu}(\bw)$ obeys the cyclic symmetry:
\begin{equation}\label{eqn:w_cyclic}
 G_{g,\mu^1,\mu^2,\mu^3}(\lo{w})
=G_{g,\mu^2,\mu^3,\mu^1}(w_2,w_3,w_1).
\end{equation}

Note $\sqrt{-1}^{\ell(\vmu)}G_{g,\vmu}(\bw)$ is a rational function in 
$w_1, w_2, w_3$ with $\QQ$-coefficients, and is homogeneous of degree $0$,
so it suffices to work with $\bw=(1,\tau,-\tau-1)$. For such weights, we shall
write
$$
G_{g,\vmu}(\tau)=G_{g,\vmu}(1,\tau,-\tau-1).
$$
Then \eqref{eqn:w_cyclic} becomes 
$$
G_{g,\up{\mu}}(\tau) =
G_{g,\mu^2,\mu^3,\mu^2}(-1-\frac{1}{\tau}) =
G_{g,\mu^3,\mu^1,\mu^2}(\frac{-1}{\tau+1}). 
$$

We now form three-partition generating functions. We let $\lam$ and
$p^i=(p^i_1,p^i_2,\ldots)$ be formal variables; given a partition
$\mu$ we define $p^i_\mu=p^i_1\cdots p^i_{\ell(\mu)}$; (note
$p^i_\emptyset=1$); for $p^1$, $p^2$ and $p^3$, we abbreviate
$$
\bp=(\up{p})\and \bp_\vmu=p^1_{\mu^1}p^2_{\mu^2}p^3_{\mu^3}.
$$
The generating functions are
$$
G_\vmu(\lam;\bw)=\sum_{g=0}^\infty
\lam^{2g-2+\ell(\vmu)}G_{g,\vmu}(\bw)
\quad \textup{and}\quad
G(\lam;\bp;\bw)=\sum_{\vmu\in\cP^3_+} G_\vmu(\lam;\bw)\bp_\vmu.
$$
The generating fucntions for not necessarily connected domain curves
are
\begin{equation}\label{eqn:Gmu}
G^\bu(\lam;\bp;\bw)=\exp(G(\lam;\bp;\bw)) =1+\sum_{\vmu\in\cP_+^3 }
G^\bu_\vmu(\lam;\bw)\bp_\vmu
\end{equation}
and
\begin{equation}\label{eqn:Gchimu}
G^\bu_\vmu(\lam;\bw) = \sum_{\chi\in2\ZZ,\chi\leq 2\ell(\vmu)}
\lam^{-\chi+\ell(\vmu)}G^\bu_{\chi,\vmu}(\bw).
\end{equation}
Finally, we define $ G_\vmu(\lam;\tau)$, $G(\lam;\bp;\tau)$,
$G^\bu(\lam;\bp;\tau)$ and $G^\bu_\vmu(\lam;\tau)$ similarly.
We will relate $G^\bu_\vmu(\lam;\tau)$ to $\tilde{\cW}_{\vmu}(q)$
in Theorem \ref{thm:Gthree}.

%% file: sec3.tex
\section{Relative Formal Toric Calabi-Yau Threefolds}\label{sec:FTCY}

In this section, we will introduce formal toric Calabi-Yau (FTCY)
graphs, and construct their associated relative FTCY threefolds.

\subsection{Toric Calabi-Yau Threefolds}\label{sec:toric}

Given a smooth toric Calabi-Yau threefold $Y$, we let $Y^1$ (resp. $Y^0$) be
the union of all 1-dimensional (resp. $0$-dimensional) $(\CC^*)^3$
orbit closures in $Y$. We assume that
\begin{equation}\label{eqn:Zone}
\textit{ $Y^1$ is connected and $Y^0$ is nonempty.}
\end{equation}
Under the above condition, we will find a distinguished
subtorus $T\subset (\CC^*)^3$ and use the $T$-action to construct a
planar trivalent graph $\Ga_Y$. FTCY graphs that will be defined
in Section \ref{sec:graph} are generalization of the planar
trivalent graphs  associated to smooth toric Calabi-Yau
threefolds just mentioned.

We first describe the distinguished subtorus $T$.
We pick a fixed point $p\in Y^0$ and look at
the $(\CC^*)^3$ action on the tangent space $T_p Y$ and its top
wedge $\Lambda^3 T_p Y$. Clearly, the later defines a homomorphism
$\alpha_p: (\CC^*)^3 \to \CC^*$, which by the Calabi-Yau condition
and the connectedness of $Y^1$ is independent of choice $p$. We
define
$$
T\stackrel{\mathrm{def}}{=}
\mathrm{Ker}\alpha_p \cong (\CC^*)^2.
$$

We next desribe the planar trivalent graph $\Ga_Y$.
We let $\Lambda_T$ be the group of irreducible characters of $T$, i.e.,
$$ 
\Lambda_T \stackrel{\mathrm{def}}{=}\Hom(T,\CC^*)\cong 
\ZZ^{\oplus 2}.
$$
We let $T_\RR\cong U(1)^2$  be the maximal compact subgroup of $T$; let
${\ft}_\RR$ and $\ft_\RR^\vee$ be its Lie algebra and its dual; let
$\mu: Y\to {\ft}_\RR^\vee$ be the moment map of the
$T_\RR$-action on $Y$. Because of the canonical isomorphism
${\ft}_\RR\dual\cong \Lambda_T\otimes_\ZZ \RR$, 
the image of $Y^1$ under $\mu$ forms a {\it planar trivalent graph}
$\Ga$ in $\Zt\otimes_\ZZ \RR$. Some examples are shown in Figure 1.
\begin{figure}[h]\label{figure1}
\begin{center}
\psfrag{C3}{\small $\CC^3$}
\psfrag{O(-1)+O(-1)}{\small$\cO_\Po(-1)\oplus\cO_\Po(-1)$}
\psfrag{O(-3)}{\small $\cO_{\PP^2}(-3)$}
\includegraphics[scale=0.6]{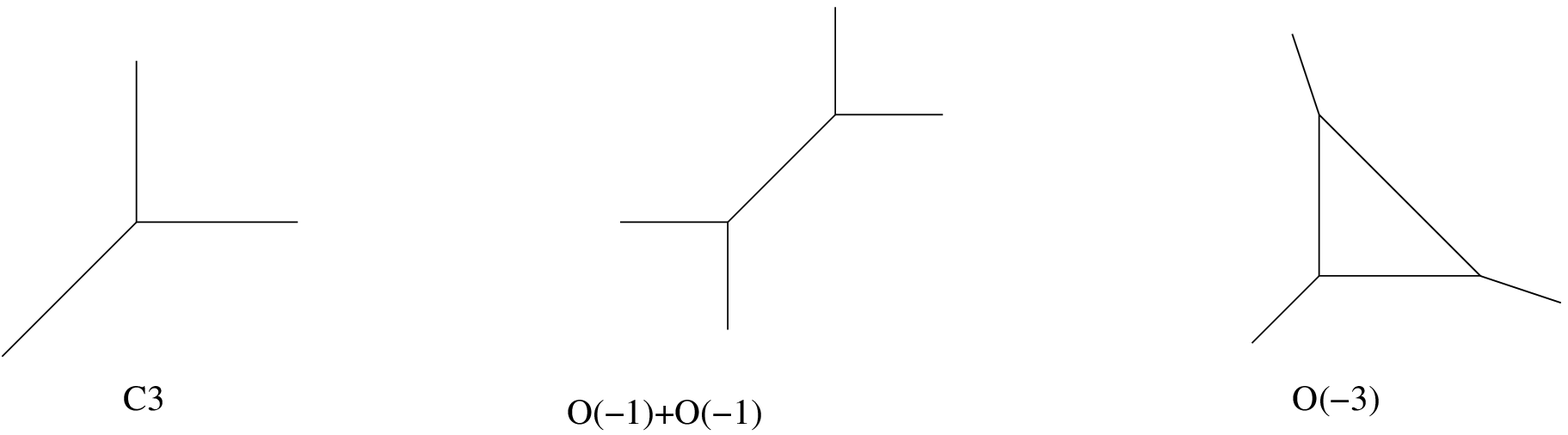}
\end{center}
\caption{ }
\end{figure}

The graph $\Ga_Y$ encodes the information of $Y$ in that its edges
and vertices correspond to irreducible components of $Y^1$ and fixed
points $Y^0$;  the slope of an edge determines the $T$-action on the
corresponding component of $Y^1$. 

Let $\hY$ be the formal completion of $Y$ along $Y^1$. Then
$\hY$ is a smooth formal Calabi-Yau scheme and inherits the
$T$-action on $Y$. The formal Calabi-Yau 
scheme $\hY$ together with the $T$-action can be reconstructed from 
the graph $\Ga_Y$ (cf. (a) in Section 
\ref{sec:localC} below).
The construction of a relative FTCY threefold from a FTCY
graph (given in Section \ref{sec:threefold}) can be viewed
as generalization of this reconstructing procedure.

\subsection{Relative Toric Calabi-Yau Threefolds}\label{sec:localC}
A smooth relative toric Calabi-Yau threefold is a pair $(Y,D)$ where
$Y$ is a smooth toric threefold, $D$ is a possibly disconnected,
smooth $(\CC^*)^3$ invariant divisor of $Y$, such that the relative Calabi-
Yau condition holds: 
$$ 
\Lambda^3 \Omega_Y(\log D)\cong \cO_Y.
$$
A toric Calabi-Yau threefold can be viewed as a relative Calabi-Yau
threefold where the divisor $D$ is empty. 

We now describe in details
three examples of relative toric Calabi-Yau threefolds and their associated graphs,  
as they are the building blocks of the definitions and constructions
in the rest of Section \ref{sec:FTCY}:
\begin{enumerate}
\item[(a)] $Y$ is the total space of $\cO_\Po(-1+n)\oplus \cO_\Po(-1-n)$.
\item[(b)] $Y$ is the total space of $\cO_\Po(n)\oplus \cO_\Po(-1+n)$;
$D$ is its fiber over $q_1=[1,0]\in\Po$.
\item[(c)] $Y$ is the total space of $\cO_\Po(n)\oplus \cO_\Po(-n)$;
$D$ is the union of its fibers over $q_0=[0,1]$ and 
$q_1=[1,0]$ in $\Po$.
\end{enumerate} 
\begin{figure}[h]\label{figure2}
\begin{center}
\psfrag{w1}{\footnotesize $w_1$}
\psfrag{-w1}{\footnotesize $-w_1$}
\psfrag{w2}{\footnotesize $w_2$}
\psfrag{-w2}{\footnotesize $-w_2$}
\psfrag{w3}{\footnotesize $-w_1-w_2$}
\psfrag{w4}{\footnotesize $w_1+(1-n)w_1$}
\psfrag{w5}{\footnotesize $-w_2+ nw_1$}
\psfrag{f}{\footnotesize $w_2-nw_1$}
\psfrag{-f}{\footnotesize $-w_2+nw_1$}
\psfrag{v0}{\footnotesize $v_0$}
\psfrag{v1}{\footnotesize $v_1$}
\psfrag{Epp}{\footnotesize (a) $\cO_\Po(-1+n)\oplus \cO_\Po(-1-n)$}
\psfrag{Epf}{\footnotesize (b) $\cO_\Po(n)\oplus \cO_\Po(-1-n)$}
\psfrag{Eff}{\footnotesize (c) $\cO_\Po(n)\oplus \cO_\Po(-n)$}
\includegraphics[scale=0.55]{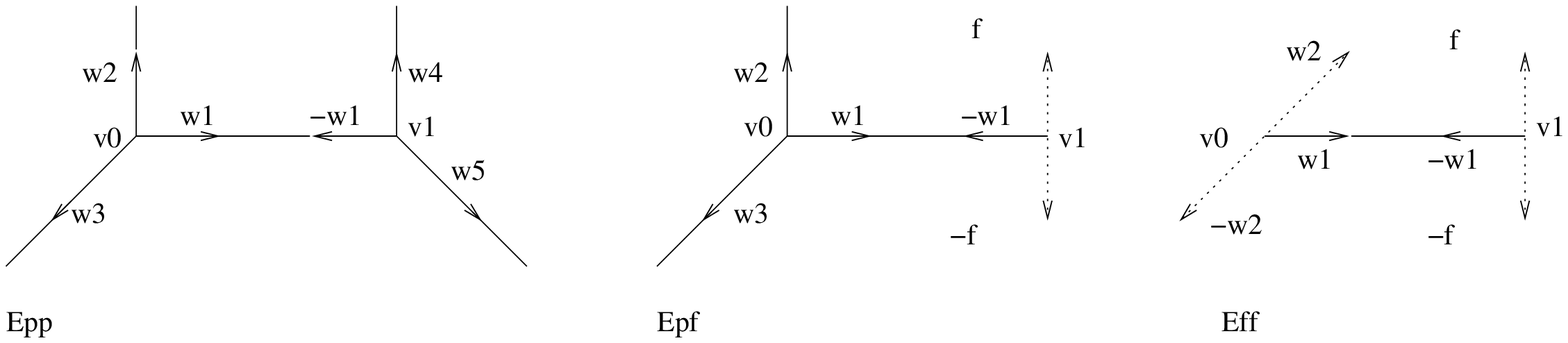}
\end{center}
\caption{ }
\end{figure} 
In Figure 2, the edge connecting the two vertices $v_0$ and $v_1$ corresponds to the 
zero section $\Po$, which
is a 1-dimensional $(\CC^*)^3$ orbit closure in $Y$; the 
vertices $v_0, v_1$ correspond to the $(\CC^*)^3$ fixed points
$q_0,q_1\in \Po$, respectively. 

In Case (a), $Y$ is a toric Calabi-Yau
threefold, so we may specify a subtorus $T$ as in Section \ref{sec:toric}.
The weights of the $T$-action on the fibers of
$T\Po, \cO_\Po(-1+n), \cO_\Po(-1-n)$ at the $T$-fixed point
$q_0\in\Po$ are given by $w_1,w_2,w_3\in \Lambda_T\cong \Zt$, respectively.
We must have $w_1+w_2+w_3=0$ because $T$ acts on $\Lambda^3 T_{q_0} Y$ trivially.
The weights of the $T$-action on the fibers of  
$T\Po, \cO_\Po(-1+n), \cO_\Po(-1-n)$ at the other fixed point
$q_1\in \Po$ are given by $-w_1,w_2+(1-n)w_1, w_3+(1+n)w_1= -w_2+nw_1$, respectively.
From the graph in Figure 2(a), one can read off the degrees
of the two summand of $N_{\Po/ Y}$ and the $T$-action on $Y$; 
$Y$ together with the $T$-action can be reconstructed from the graph.

Similarly, from the graph in Figure 2(b) (resp. (c)), one can reconstruct
the pair $(Y,D)$ in (b) (resp. (c)) together with the $T$-action; the weights of
the $T$-action at fixed points can be read off from the graph 
as follows:

\medskip

\begin{tabular}{c|c|c|c||c|c|c|c}
(b)      & $T\Po$ & $\cO_\Po(n)$ & $\cO_\Po(-1-n)$ & 
(c)      & $T\Po$ & $\cO_\Po(n)$ & $\cO_\Po(-n)$    \\ \hline
$q_0$      & $w_1$  & $w_2$        & $w_3=-w_1-w_2$ &
$q_0$      & $w_1$  & $w_2$        & $-w_2$  \\ \hline
$q_1$ & $-w_1$ & $w_2-nw_1$   & $\begin{array}{l}w_3+(1+n)w_1\\=-w_2+n w_2\end{array}$ &
$q_1$ & $-w_1$ & $w_2-nw_1$   & $-w_2+nw_1$
\end{tabular} 

\subsection{FTCY Graphs}\label{sec:graph}
In this subsection, we will introduce formal toric Calabi-Yau (FTCY)
graphs, which are graphs together with local embedding into
$\RR^2$ endowed with the standard orientation and integral lattice
$\Zt\sub\RR^2$. 

As will be clear later, assigning a slope to an edge depends
on the orientation of the edge. For book keeping purpose, we shall
associate to each edge two (oppositely) oriented edges; for an
oriented edge we can talk about its initial and terminal vertices. To
recover the graph, we simply identify the two physically identical
but oppositely oriented edges as one (unoriented) edge. This
leads to the following definition.

\begin{defi}[graphs]\label{CY2}
A graph $\Ga$ consists of a set of {\em oriented edges}
$E\ori(\Ga)$, a set of {\em vertices} $V(\Ga)$, an {\em orientation
reversing map} $\frev\mh E\ori(\Ga)\to E\ori(\Ga)$, an {\em initial
vertex map} $\fv_0\mh E\ori(\Ga)\to V(\Ga)$ and a {\em terminal
vertex map} $\fv_1\mh E\ori(\Ga)\to V(\Ga)$, such that $\frev$ is a
fixed point free involution; that both $\fv_0$ and $\fv_1$ are
surjective and $\fv_1=\fv_0\circ\frev$. We say $\Gamma$ is {\em
weakly trivalent} if $|\fv^{-1}(v)|\leq 3$ for $v\in V(\Gamma)$.
\end{defi}

For simplicity, we will abbreviate $\frev(e)$ to $-e$. Then the
equivalence classes $E(\Ga)=E\ori(\Ga)/\{\pm 1\}$ is the set of
edges of $\Gamma$ in the ordinary sense. 
In case $\Gamma$ is weakly trivalent, we shall denote by $V_1(\Ga)$, 
$V_2(\Ga)$ and  $V_3(\Ga)$ the sets of univalent, bivalent, and trivalent 
vertices of $\Ga$; we shall also define
$$
E^{\ff}(\Ga)=\{ e\in  E\ori(\Ga)\mid
\fv_1(e)\in V_1(\Ga)\cup V_2(\Ga)\}
$$ 
which consists of oriented edges whose terminal edges are not
trivalent. Finally, we fix a standard basis $\{u_1,u_2\}$ of $\Zt$ such that
the ordered basis $(u_1,u_2)$ determines the orientation on $\RR^2$.

\begin{defi} [FTCY graphs] \label{def:FTCYgraph}
A  {\em formal toric Calabi-Yau (FTCY) graph} is a weakly trivalent
graph $\Ga$ together with  a {\em position} map 
$$
\fp: E\ori(\Ga)\lra\Zt-\{0\}
$$
and a {\em framing} map
$$
\ff: E^{\ff}(\Ga)  \lra \Zt-\{0\}, 
$$
such that (see Figure 3)
\begin{enumerate}
\item[T1.] $\fp$ is anti-symmetric: $\fp(-e)=-\fp(e)$.
\item[T2.] $\fp$ and $\ff$ are balanced:
\begin{itemize} 
\item For a bivalent vertex $v\in V_2(\Ga)$ with
      $\fv_1\upmo(v)=\{e_1,e_2\}$,  
$\fp(e_1)+\fp(e_2)=0$ and $\ff(e_1)+ \ff(e_2)=0$.  
\item For a trivalent  $v\in V_3(\Ga)$ with $\fv_0\upmo(v)=\{e_1,e_2,e_3\}$, 
        $\fp(e_1)+\fp(e_2)+\fp(e_3)=0$.
\end{itemize} 
\item[T3.] $\fp$ and $\ff$ are primitive:
\begin{itemize}
\item  For a trivalent vertex $v\in V_3(\Ga)$
        with $\fv_0\upmo(v)=\{e_1,e_2, e_3\}$, any two vectors in
        $\{\fp(e_1),\fp(e_2),\fp(e_3)\}$ form an integral basis of
        $\Zt$
\item  For $e\in E^{\ff}(\Ga)$,
        $\fp(e)\wedge\ff(e)=u_1\wedge u_2$. 
\end{itemize} 
\end{enumerate} 
\begin{figure}[h]\label{figure3}
\begin{center}
\psfrag{V1}{\footnotesize $\fv_1^{-1}(v)=\{e\}$}
\psfrag{V2}{\footnotesize $\fv_1^{-1}(v)=\{e_1,e_2\}$}
\psfrag{V3}{\footnotesize $\fv_0^{-1}(v)=\{e_1,e_2,e_3\}$}
\psfrag{v}{\footnotesize $v$}
\psfrag{pe}{\footnotesize $\fp(e)$}
\psfrag{p1}{\footnotesize $\fp(e_1)$}
\psfrag{p2}{\footnotesize $\fp(e_2)$}
\psfrag{p3}{\footnotesize $\fp(e_3)$}
\psfrag{e1}{\footnotesize $e_1$}
\psfrag{e2}{\footnotesize $e_2$}
\psfrag{e3}{\footnotesize $e_3$}
\psfrag{e}{\footnotesize $e$}
\psfrag{f1}{\footnotesize $\ff(e_1)$}
\psfrag{f2}{\footnotesize $\ff(e_2)$}
\psfrag{fe}{\footnotesize $\ff(e)$}
\includegraphics[scale=0.7]{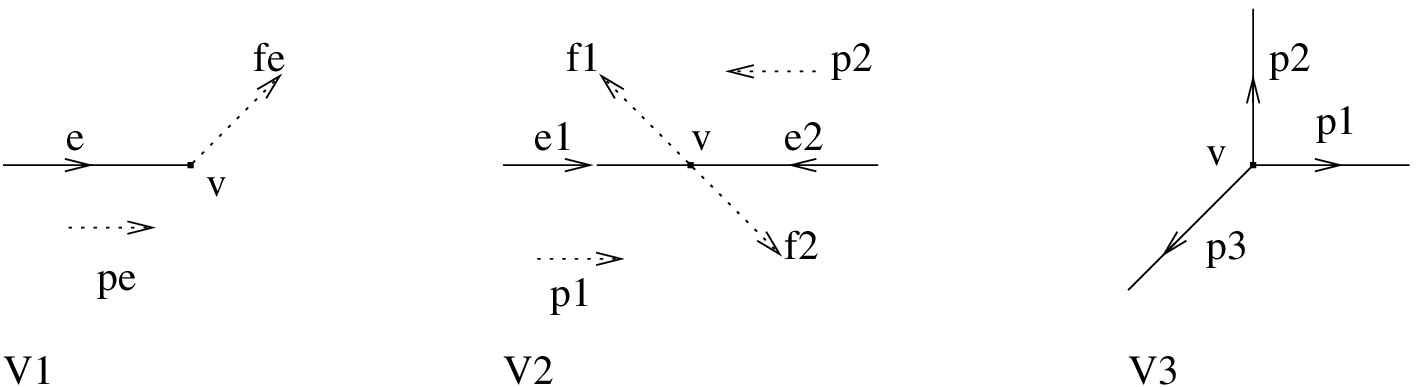}
\end{center}
\caption{ }
\end{figure}
We say $\Ga$ is a {\em regular} FTCY graph if it has no bivalent
vertex.
\end{defi}

The neighborhood of an oriented edge $e$ of an FTCY graph $\Ga$
looks like (a), (b), or (c) in Figure 2 (if we add vectors
$\{-\ff(e)\mid \fv_1(e)\in V_1(E)\}$). So to each edge we may
associate a relative Calabi-Yau threefold $(Y^e,D^e)$ where
$Y^e$ is the total space of the direct sum of two line
bundles over $\Po$. Since $e$ is oriented, we may assign
$L^e$, one of the two line bundles, to $e$ as follows:
in Figure 2(a), (b), or
(c), if $\fp(e)=w_1$ points to the right, so that
$v_0$ (resp. $v_1$) is the initial (resp. terminal) vertex of $e$, 
then the weight  of $T$-action
on $L^e_{q_0}$ (resp. $L^e_{q_1}$) is given by the upward vector
at $v_0$ (resp. $v_1$), denoted by $\fl_0(e)$ (resp. $\fl_1(e)$).
Then $Y^e$ is the total space of $L^e\oplus L^{-e}\to \Po$.
This is the geometric interpretation of the following definition:
\begin{defi} \label{def:l}
Let $\Ga$ be an FTCY graph.
We define $\fl_0, \fl_1: E\ori(\Ga)\longrightarrow \Zt$
as follows:
$$
\fl_0(e)=\begin{cases}
-\ff(-e), & \fv_0(e)\notin V_3(\Ga),\\
\fp(e_{01}), & \fv_0(e)\in V_3(\Ga).
\end{cases}
\quad\quad
\fl_1(e)=\begin{cases}
\ff(e), & \fv_1(e)\notin V_3(\Ga),\\
\fp(e_{11}), & \fv_1(e)\in V_3(\Ga).
\end{cases}
$$
Here $e_{i1}$ is the unique oriented edge such that
$\fv_0(e_{i1})=\fv_i(e)$ and $\fp(e)\wedge \fp(e_{i1})=u_1\wedge u_2$.
\end{defi} 
\begin{figure}[h]\label{figure4}
\begin{center}
\psfrag{w1}{\footnotesize $\fp(e)$}
\psfrag{-w1}{\footnotesize $\fp(-e)$}
\psfrag{w2}{\footnotesize $\fl_0(e)$}
\psfrag{w3}{\footnotesize $\fl_1(-e)$}
\psfrag{w4}{\footnotesize $\fl_1(e)$}
\psfrag{w5}{\footnotesize $\fl_0(-e)$}
\psfrag{Epp}{\footnotesize (a) $\fv_0(e),\fv_1(e)\in V_3(\Ga)$}
\psfrag{Epf}{\footnotesize (b) $\begin{array}{l}\fv_0(e)\in V_3(\Ga)\\
\fv_1(e)\in V_1(\Ga)\cup V_2(\Ga)\end{array}$}
\psfrag{Eff}{\footnotesize (c) $\fv_0(e),\fv_1(e)\in V_1(\Ga)\cup V_2(\Ga)$}
\includegraphics[scale=0.55]{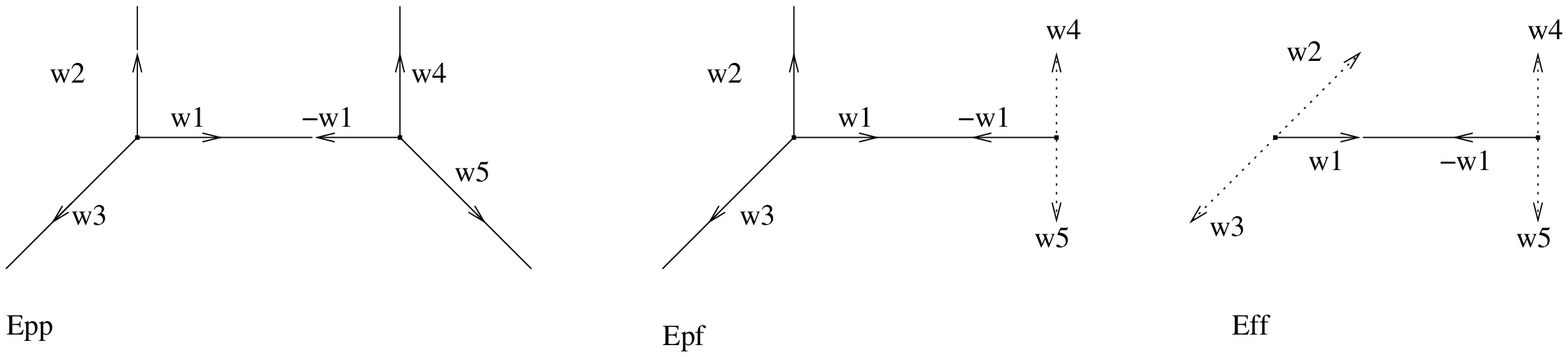}
\end{center}
\caption{ }
\end{figure} 
The degree of the line bundle $L^e$ determines an integer $n^e$:
\begin{equation}\label{eqn:degL}
\deg L^e =\begin{cases} n^e-1, & v_1\in V_3(\Ga),\\
n^e, & v_1\notin V_3(\Ga). \end{cases}
\end{equation}
This motivates the following definition.
\begin{defi}\label{def:vn}
We define $\vn:E\ori(\Ga)\to \ZZ$ by
$$
\fl_1(e)-\fl_0(e)=\begin{cases}
(1-\vn(e))\fp(e), & \fv_1(e)\in    V_3(\Ga),\\
  -\vn(e) \fp(e), & \fv_1(e)\notin V_3(\Ga).
\end{cases}
$$
We write $n^e$ for $\vn(e)$.
\end{defi}
Note that $n^{-e}=-n^{e}$.

\subsection{Operations on FTCY Graphs} \label{sec:operation}
In this subsection, we define four operations on FTCY graphs:
smoothing, degeneration, normalization, and gluing. These operations
extend natural operations on toric Calabi-Yau threefolds.

The first operation is the smoothing of a bivalent vertex
$v\in V_2(\Gamma)$. In doing this, we shall eliminate the vertex $v$
and combine the two edges attached to $v$.
\begin{defi}[smoothing]
The {\em smoothing} of $\Gamma$ along a bivalent vertex $v\in V_2(\Ga)$ is a
graph $\Gamma_v$ that has vertices $V(\Gamma)-\{v\}$, oriented edges
$E\ori(\Gamma)/\sim$ with the equivalence $\pm e_1\sim\mp e_2$ for
$\{e_1,e_2\}=\fv_1\upmo(v)$. The maps $\fv_0$, $\fv_1$, $\fp$ and
$\ff$ descends to corresponding maps on $\Ga_v$, making it a FTCY
graph.
$($See Figure 5: $\Ga_3$ is the smoothing of $\Ga_2$ along $v$.$)$
\end{defi}

The reverse of the above construction is called a degeneration.
\begin{defi}[degeneration]
Let $\Gamma$ be a FTCY graph and let $e\in E\ori(\Gamma)$ be an
edge. We pick a lattice point $\ff_0\in \Zt$ so that
$\fp(e)\wedge\ff_0=u_1\wedge u_2$. The {\em degeneration} of $\Gamma$ at
$e$ with framing $\ff_0$ is a graph $\Gamma_{e,\ff_o}$ whose edges
are $E\ori(\Gamma)\cup\{\pm e_1,\pm e_2\}-\{\pm e\}$ and whose
vertices are $V(\Gamma)\cup\{v_0\}$; its initial vertices
$\tilde\fv_0$, terminal vertices $\tilde\fv_1$, position map
$\tilde\fp$ and framing map $\tilde\ff$ are identical to those of
$\Gamma$ except
\begin{eqnarray*}
&& \tilde\fv_0(e_1)=\fv_0(e),\
\tilde\fv_1(e_1)=\tilde\fv_1(e_2)=v_0,\
\tilde\fv_0(e_2)=\fv_1(e),\\
&& \tilde\fp(e_1)=-\tilde\fp(e_2)=\fp(e),\ \tilde\ff(e_1)=\tilde\ff(e_2)=\ff_0.
\end{eqnarray*}
$($See Figure 5: $\Ga_2$ is the degeneration of $\Ga_3$ at $e$ with framing
$\ff_0$.$)$
\end{defi}

The normalization is to separate a graph along a bivalent vertex
and the gluing is its inverse.
\begin{defi}[normalization]
Let $\Ga$ be a FTCY graph and let $v\in V_2(\Ga)$ be a bivalent
vertex. The {\em normalization} of $\Ga$ at $v$ is a graph $\Ga^v$
whose edges are the same as that of $\Ga$ and whose vertices are
$V(\Ga)\cup\{v_1,v_2\}-\{v\}$; its associated maps $\tilde\fv_0$,
$\tilde\fv_1$, $\tilde\fp$ and $\tilde\ff$ are identical to that of
$\Ga$ except for $\{e_1,e_2\}=\fv_1\upmo(v)$, $
\tilde\fv_1(e_1)=v_1$ and $\tilde\fv_1(e_2)=v_2$.
$($See  Figure 5: $\Ga_1$ is the normalization of $\Ga_2$ at $v$.$)$
\end{defi}

\begin{defi}[gluing]
Let $\Ga$ be a FTCY graph and let $v_1, v_2\in V_1(\Ga)$ be two
univalent vertices of $\Ga$. Let $\ff_i=\ff(e_i)$, where
$\{e_i\}=\fv_1^{-1}(v_i)$. Suppose $\fp(e_1)=-\fp(e_2)$ and
$\ff_1=-\ff_2$. We then identify $v_1$ and $v_2$ to form a single
vertex, and keep the framing $\ff(e_i)=\ff_i$. The resulting graph
$\Ga^{v_1,v_2}$ is called the {\em gluing} of $\Ga$ at $v_1$ and
$v_2$. $($See Figure 5: $\Ga_2$ is the gluing of $\Ga_1$ at $v_1$ and $v_2$.$)$
\end{defi}

\begin{figure}[h]\label{figure5}
\begin{center}
\psfrag{v}{\footnotesize $v$} 
\psfrag{v1}{\footnotesize $v_1$}
\psfrag{v2}{\footnotesize $v_2$} 
\psfrag{e}{\footnotesize $e$}
\psfrag{e1}{\footnotesize $e_1$}
\psfrag{e2}{\footnotesize $e_2$}
\psfrag{f0}{\footnotesize $\ff_0$} 
\psfrag{-f0}{\footnotesize $-\ff_0$} 
\psfrag{f1}{\footnotesize $\ff_1$} 
\psfrag{-f1}{\footnotesize $-\ff_1$} 
\psfrag{f2}{\footnotesize $\ff_2$} 
\psfrag{-f2}{\footnotesize $-\ff_2$} 
\psfrag{Ga}{\footnotesize
$\Ga_2=(\Ga_1)^{v_1,v_2}=(\Ga_3)_{\ff_0,e}$}
\psfrag{res}{\footnotesize $\Ga_1=(\Ga_2)^v$}
\psfrag{smo}{\footnotesize $\Ga_3=(\Ga_2)_v$}
\includegraphics[scale=0.7]{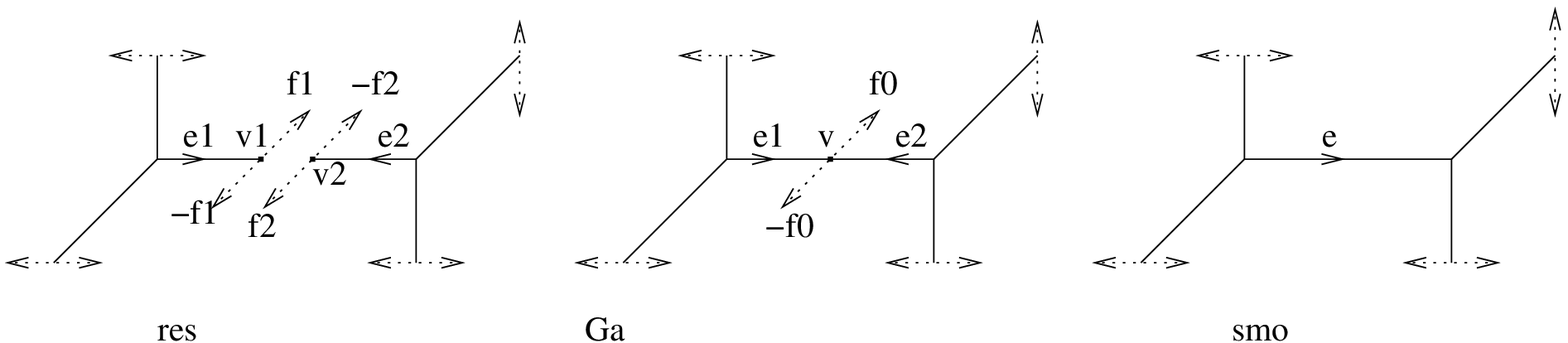}
\end{center}
\caption{ }
\end{figure}
\smallskip

It is straightforward to generalize smoothing and normalization to
subset $A$ of $V_2(\Gamma)$. Given $A\subset V_2(\Gamma)$, let
$\Gamma_A$ denote the smoothing of $\Gamma$ along $A$,  and let
$\Gamma^A$ denote the normalization of $\Gamma$ along $A$. There are
surjective maps
$$
\pi_A: E(\Gamma)\to E(\Gamma_A),\quad \pi^A: V(\Gamma^A) \to
V(\Gamma).
$$

\subsection{Relative FTCY Threefolds}
\label{sec:threefold}

In this subsection we will introduce relative formal toric
Calabi-Yau (FTCY) threefolds.

Given a FTCY graph $\Ga$, we will construct a pair $\Yr=(\hY,\hD)$,
where $\hY$ is a threefold, possibly with normal crossing
singularities, $\hD\sub \hY$ is a relative divisor, so that
$\Yr=(\hY,\hD)$ is a formal relative Calabi-Yau threefold:
\begin{equation}
\wedge^3 \Omega_{\hY}(\log \hD)\cong\cO_{\hY}.
\end{equation}
The pair $(\hY,\hD)$ admits a $T$-action so that the action on
$\Lambda^3 T_p\hY$ is trivial for any fixed point $p$. As a set, the
scheme $\hY$ is a union of $\Po$'s, each associated to an edge of
$\Ga$; two $\Po$ intersect exactly when their associated edges share
a common vertex; the normal bundle to each $\Po$ in $\hY$ and the
$T$-action on $\hY$ are dictated by the data encoded in the graph
$\Ga$. We will also specify a $T$-invariant divisor $\hL\subset
\hD$.

In the following construction, we will use the notation introduced
in Section \ref{sec:graph}.

\subsubsection{Edges}\label{sec:edge}
Let $e\in E\ori(\Gamma)$ with $v_0$ and $v_1$ its initial and
terminal vertices. We let $T$ acts on $\Po$ by
$$
t\cdot[X,Y]=[\fp(e)(t)X,Y],\quad t\in T.
$$
Here we view $\fp(e)$ as an element in $\Lambda_T=\Hom(T,\CC^*)$. We
denote the two fixed points by $q_0=[0,1]$ and $q_1=[1,0]$. Next we
let $L^e\to \Po$ be the line bundle of 
$$ 
\deg L^e=\begin{cases} 
n^e-1,& v_1\in V_3(\Ga),\\
n^e,  & v_1\notin V_3(\Ga),
\end{cases}
$$
where $n^e=\vn(e)$ is defined in Definition \ref{def:vn}.
We let $L^{-e}$ be the line bundle on $\Po$ of degree $n^{-e}-1$
(resp. $n^{-e}$) when $v_0$ is trivalent (resp. non-trivalent). We
then assign the $T$-action at $L^e_{q_0}$ be $\fl_0(e)$ and at
$L^e_{q_1}$ be $\fl_1(e)$; assign the $T$-action at $L^{-e}_{q_1}$
be $\fl_0(-e)$ and at $L^{-e}_{q_0}$ be $\fl_1(-e)$ (see Figures 4).

Next, we let $\Si(e)$ be the formal completion of the total space of
$L^e\oplus L^{-e}$ along its zero section. The $T$-actions on $L^e$
and on $L^{-e}$ induce a $T$-action on $\Si(e)$. By construction,
there is a $T$-isomorphism
\begin{equation}\label{22.5}
\Si(e)\cong \Si(-e)
\end{equation}
that sends $[0,1]\in \Si(e)$ to $[1,0]\in\Si(-e)$ and sends the
first summand $L^e$ in $N_{\Po/\Si(e)}$ to the second summand
$L^e$ in $N_{\Po/\Si(-e)}$.

It is clear that
\begin{equation}\label{22.6}
\wedge^3\Omega_{\Si(e)}\cong p\sta\cO_{\Po}(c),
\end{equation}
where $p\mh \Si(e)\to\Po$ is the projection, and 
$c =  \# (\{ v_0, v_1\} \cap V_3(\Ga))-2$.

\subsubsection{Gluing along Trivalent Vertices}\label{sec:trivalent}
Given $v\in V_3(\Gamma)$ with $\fv_0^{-1}(v)=\{e_1,e_2,e_3\}$ so
indexed so that $\fp(e_1),\fp(e_2),\fp(e_3)$ is in counter-clockwise
order.

To glue the formal scheme $\Si(e_1)$, $\Si(e_2)$ and $\Si(e_3)$, we
introduce a $T$-scheme
$$
\Si(v)=\spec \CC[\![x_1,x_2,x_3]\!],\quad\text{with $T$-action}\
t\cdot x_i=\fp(e_i)(t) x_i\quad \forall\ t\in T,
$$
and gluing morphisms
\begin{equation}\lab{glu-01}
\psi_{e_k v}: \Si(v)\lra \Si(e_k)
\end{equation}
according to the following rule. First, we let $\hat \Si(e_k)$ be
the formal completion of $\Si(e_k)$ along $q_0\in\Po\sub \Si(e_k)$.
$\hat \Si(e_k)$ is a formal $T$-scheme and is $T$-isomorphic to the
$T$-scheme
$$\spec\CC[\![y_1,y_2,y_3]\!]; \quad t\cdot
y_i=\fp(e_{i+k})(t)y_i
$$
\footnote{Here we agree that $e_{k+3}=e_k$; same for $x_{i+3}=x_i$
later.}such that $L^{e_k}_{q_0}$, $L^{-e_k}_{q_1}$, $T_{q_0}\PP^1$
are mapped to $\CC\frac{\pa}{\pa y_1}$, $\CC\frac{\pa}{\pa y_2}$,
$\CC\frac{\pa}{\pa y_3}$, respectively. The gluing morphism
$\psi_{e_kv}$ is the composite of
\begin{equation}\label{22.1}
\spec\CC[\![y_1,y_2,y_3]\!]\lra \hat\Si(e_k)\lra \Si(e_k)
\end{equation}
with the $T$-isomorphism
$$\Si(v)\equiv\spec\CC[\![x_1,x_2,x_3]\!]\lra
\spec\CC[\![y_1,y_2,y_3]\!]
$$
defined by $y_i\mapsto x_{k+i}$.

Using the morphisms $\psi_{e_kv}$, we can glue $\Si(e_1)$ and
$\Si(e_2)$ and then glue $\Si(e_3)$ onto it via the cofiber products
$$
\begin{CD}
\Si(e_1) @>>> \Si(e_1)\coprod_{\Si(v)} \Si(e_2)\\
@AAA @AAA\\
\Si(v) @>>> \Si(e_2)
\end{CD}
\qquad
\begin{CD}
\Si(e_3) @>>> \Si(e_1)\coprod_{\Si(v)} \Si(e_2)\coprod_{\Si(v)} \Si(e_3)\\
@AAA @AAA\\
\Si(v) @>>> \Si(e_1)\coprod_{\Si(v)} \Si(e_2)
\end{CD}
$$
Since the gluing map $\psi_{e_k v}$ are $T$-equivalent, the
$T$-actions on $\Si(e_k)$ descend to the glued scheme.

\subsubsection{Gluing along Bivalent Vertices}\label{sec:bivalent}
Next we glue $\Si(e_1)$ and $\Si(e_2)$ in case
$\{e_1,e_2\}=\fv_0^{-1}(v)$ for a $v\in V_2(\Gamma)$. Note that $
\fl_0(e_1)+\fl_0(e_2)=0$. Let $\Si(v)$ be the formal $T$-scheme
$$
\Si(v)=\spec\CC[\![ x_1,x_2]\!],\quad t\cdot x_i=\fl_0(e_i)(t)x_i;
$$
we let the gluing morphism $\psi_{e_kv}$ of \eqref{glu-01} be the
composite of (\ref{22.1}) with the $T$-morphism
$$
\spec\CC[\![ x_1,x_2]\!]\lra \spec\CC[\![y_1,y_2,y_3]\!]
$$
defined via $y_3\mapsto 0$, $y_1$ and $y_2$ map to $x_1$ and $x_2$
respectively in case $k=1$ and to $x_2$ and $x_1$ respectively in
case $k=2$. We can glue $\Si(e_1)$ and $\Si(e_2)$ along $\Si(v)$ via
the cofiber product as before.

\subsubsection{Univalent vertices}\label{sec:univalent}
Lastly, we consider the case $e_1=\fv_0\upmo(v)$ for a $v\in
V_1(\Ga)$. Note that $\fl_0(e)+\fl_1(-e)=0$. Let $\Si(v)$ be the formal $T$-scheme
$$
\Si(v)=\spec\CC[\![ x_1,x_2]\!],\quad t\cdot
x_1=\fl_0(e_1)(t)x_1,\quad t\cdot x_2=\fl_1(-e_1)(t)x_2;
$$
and define $\psi_{e_1v}$ in \eqref{glu-01}. We let $\hD^v$ be the
image divisor $\psi_{e_1v}(\Si(v))\sub \Si(e_1)$ and consider it as
part of the relative divisor of the formal Calabi-Yau scheme $\Yr$
we are constructing.

Let $L(v)\sub \Si(v)$ be the divisor defined by $x_2=0$, and let
$\hL^v=\psi_{e_1v}(L(v))\sub \hD^v$.

\subsubsection{Final step}
Now it is standard to glue all $\Si(e)$ to form a scheme $\hat Y$.
We first form the disjoint union
$$
\coprod_{e\in E\ori(\Ga)} \Si(e);
$$
because of (\ref{22.5}), the orientation reversing map
$E\ori(\Gamma) \to E\ori(\Gamma)$ defines a fixed point free
involution
$$
\tau: \coprod_{e\in E\ori(\Ga)} \Si(e)\lra \coprod_{e\in E\ori(\Ga)}
\Si(e);
$$
we define $\tilde Y$ be its quotient by $\tau$. Next, for each
trivalent vertex $v$ with $\fv_0\upmo(v)=\{e_1,e_2,e_3\}$, we glue
$\Si(e_1), \Si(e_2), \Si(e_3)$ along $\Si(v)$; for each bivalent
vertex $v$ with $\fv_0\upmo(v)=\{e_1,e_2\}$, we glue $\Si(e_1)$ and
$\Si(e_2)$ along  $\Si(v)$. We denote by $\hY$ the resulting scheme
after completing all the gluing associated to all trivalent and
bivalent vertices. The $T$-action on $\Si(e)$'s descends to a
$T$-action on $\hY$. Finally, for each univalent vertex $v$ with
$e=\fv_0\upmo(v)$, we let $\hD^v\sub \Si(e)$ be the divisor defined
in Section \ref{sec:univalent}. The (disjoint) union of all such $\hD^v$ form a
divisor $\hD$ that is the relative divisor of $\hY$. Since $\hD$ is
invariant under $T$, the pair $\Yr=(\hY,\hD)$ is a $T$-equivariant
formal scheme. Because of (\ref{22.6}), we have
$$
\wedge^3 \Omega_{\hY}(\log\hD)\cong \cO_{\hY};
$$
hence $\Yr=(\hY,\hD)$ is a formal toric Calabi-Yau scheme.

Following the construction, the scheme $\hY$ is smooth away from the
images $\psi_{ev}(\Si(v))$ associated to bivalent vertices $v$, and
has normal crossing singularities there. Therefore $\hY$ is smooth
when $\Gamma$ has no bivalent vertices. The relative divisor $\hD$
is the union of smooth divisor $\hD^v$ indexed by $v\in V_1(\Ga)$.
Within each divisor $\hD^v$ there is a divisor $\hL^v\sub\hD^v$
defined as in Section \ref{sec:univalent}.

For later convenience, we introduce some notation. Let  $\ee$ denote
the equivalence class $\{e,-e\}$ in $E(\Ga)$, and let $C^\ee$ denote
the projective line in $\hY$ coming from the $\Po$ in $\Si(e)$. For
$v\in V_1(\Ga)$, let $z^v$ denote the point in $\hD^v$ coming from
the closed point $q_0$ in $\Si(e)$, where $\fv_0(e)=v$.

%% file: sec4.tex
\section{Definition of Formal Relative Gromov-Witten Invariants}
\label{sec:invariants}

In this section, we will define relative Gromov-Witten invariants of
relative FTCY threefolds; the case when the relative FTCY threefold
is indecomposable gives the mathematical definition of topological
vertex.

\subsection{Moduli Spaces of Relative Stable morphisms}
\label{sec:moduli}

Let $\Ga$ be a FTCY graph and let $\Yr=(\hY,\hD)$ be its associated
scheme. Clearly, the degrees and the ramification patterns of
relative stable morphisms to $\Yr$ are characterized by {\em
effective classes} of $\Ga$:
\begin{defi}[effective class]\label{def:effective}
Let $\Ga$ be a FTCY graph. An {\em effective class} of $\Ga$
is a pair of functions $\vd\mh E(\Gamma)\to \ZZ_{\geq 0}$ and
$\vmu\mh V_1(\Gamma) \to \cP$ that satisfy
\begin{enumerate}
\item $|\vmu(v)|=\vd(\ee)$ if $v\in V_1(\Gamma)$ and $\fv_1(e)=v$;
\item $\vd(\ee_1)=\vd(\ee_2)$ if $v\in V_2(\Gamma)$ and
      $\fv_0\upmo(v)=\{e_1,e_2\}$.
\end{enumerate}
We write $\mu^v$ for $\vmu(v)$, $d^\ee$ for $\vd(\ee)$.
\end{defi}

To show that an effective class does characterize a relative stable
morphism, a quick review of its definition is in order. An ordinary
relative morphism $u$ to $(\hY,\hD)$ consists of
\begin{itemize}
\item a possibly disconnected nodal curve $X$
\item distinct smooth points
$\{q_j^v \mid v\in V_0(\Gamma), 1\leq j\leq \ell(\mu^v)\}$ in $X$
such that each connected component of X contains at least one of
these points,
\item a morphism $u\mh X\to \hY$
\end{itemize}
so that
\begin{itemize}
\item $u\upmo(\hD^v)=\sum_{j=1}^{\ell(\mu^v)} \mu_j^v q_j^v$ for
some positive integers $\mu_j^v$;
\item $u$ is {\em pre-deformable} along the singular loci
$$
\coprod_{v\in V_2(\Gamma)}\Si(v)
$$
of $\Yr$, i.e, if $v\in V_2(\Gamma)$ and $\fv_0^{-1}(v)=\{e_1,e_2\}$
, then $u^{-1}(\Si(v))$ consists of nodes of $X$, and for each $y\in
u^{-1}(\Si(v))$, $u|_{u^{-1}(\Si(e_1))}$ and $u|_{u^{-1}(\Si(e_2))}$
have the same contact order to $\Si(v)$ at $y$;
\item $u$ coupled with the marked points $q_i^v$ is a stable morphism
      in the ordinary sense.
\end{itemize}

Unless otherwise specified, all the stable morphisms in this paper
are with not necessarily connected domain.

Since
$$
H_2(\hY;\ZZ)=\bigoplus_{\ee\in E(\Gamma)} \ZZ[C^\ee],
$$
the morphism $u$ defines a map $\vd\mh E(\Gamma)\to \ZZ$ via
\begin{equation}\label{2.5}
u\lsta([X])=\sum_{\ee\in E(\Gamma)} \vd(\ee)[C^\ee].
\end{equation}
The integers $\mu_j^v$ form a partition
$$
\mu^v=(\mu^v_1,\cdots,\mu^v_{\ell(\mu^v)})
$$
and the map $\vmu\mh V_1(\Gamma)\to\cP$ is
$$
\vmu(v)=\mu^v.
$$
With this definition, the requirement (1) in Definition
\ref{def:effective} follows from (\ref{2.5}) and (2) holds since $u$
is pre-deformable.

To define relative stable morphisms to $\Yr$, we need to work with
the expanded schemes of $\Yr$ introduced in \cite{Li1}. In the case
studied, they are the associated formal schemes of the expanded
graphs of $\Ga$.

\begin{defi}
Let $\Ga$ be a FTCY graph. A {\em flat chain of length $n$} in $\Ga$
is a subgraph $\eGa\sub\Ga$ that has $n$ edges $\pm e_1,\cdots,\pm
e_n$, $n+1$ univalent or bivalent vertices $v_0,\cdots,v_n$ with
identical framings $\ff$ (up to sign) so that
$$
\fv_0(e_1)=v_0;\quad
\fv_1(e_i)=\fv_0(e_{i+1})=v_i\ i=1,\cdots,n-1;\quad
\fv_1(e_n)=v_n,
$$
and that all $\fp(e_i)$ are identical.
\end{defi}
\begin{figure}[h]\label{figure10}
\begin{center}
\psfrag{v0}{\footnotesize $v_0$} 
\psfrag{v1}{\footnotesize $v_1$}
\psfrag{v2}{\footnotesize $v_2$} 
\psfrag{vn}{\footnotesize $v_n$}
\psfrag{e1}{\footnotesize $e_1$} 
\psfrag{e2}{\footnotesize $e_2$}
\psfrag{en}{\footnotesize $e_n$} 
\psfrag{f1}{\footnotesize $\ff(e_1)$} 
\psfrag{f2}{\footnotesize $\ff(e_2)$}
\psfrag{f3}{\footnotesize $\ff(e_n)$}
\includegraphics[scale=0.7]{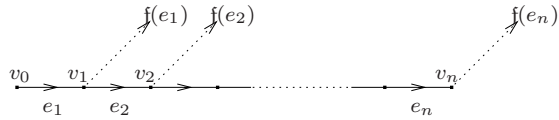}
\end{center}
\caption{A flat chain of length $n$}
\end{figure}

\begin{defi}
A {\em contraction} of a FTCY graph $\Ga$ along a flat chain
$\eGa\sub\Ga$ is the graph after eliminating all edges and bivalent
vertices of $\eGa$ from $\Ga$, identifying the univalent vertices of
$\eGa$ while keeping their framings unchanged.
\end{defi}

Given a FTCY graph $\Ga$ and a function
$$
\bm:V_1(\Ga)\cup V_2(\Ga)\longrightarrow \ZZ_{\geq 0},
$$
the expanded graph $\Ga_\bm$ is obtained by replacing each $v\in
V_1(\Ga)\cup V_2(\Ga)$ by a flat chain $\eGa^v_{m^v}$ of length
$m^v=\bm(v)$ with framings $\pm\ff(e)$, where $\fv_1(e)=v$. In
particular $\Ga_{\mathbf{0}}=\Ga$, where $\mathbf{0}(v)=0$ for all
$v\in V_1(\Ga)\cup V_2(\Ga)$. The original graph $\Ga$ can be
recovered by contracting $\Ga_\bm$ along the flat chains
$$
\{ \eGa^v_{m^v} \mid v\in V_1(\Ga)\cup V_2(\Ga) \}.
$$

We now study their associated Calabi-Yau scheme. We denote by
$(\hY,\hD)$ the associated Calabi-Yau scheme of $\Ga$ and by
$(\hY_\bm,\hD_\bm)$ that of $\Ga_\bm$. We recover the original
scheme $\hY$ by shrinking the irreducible components of $\hY_\bm$
associated to the flat chains that are contracted. This way we
define a projection
$$
\pi_\bm: \hY_\bm\lra \hY.
$$
We define a relative automorphism of $\hY_\bm$ to be an automorphism
of $\hY_\bm$ that is also a $\hY$-morphism; an automorphism of a
relative morphism $u\mh X\to (\hY_\bm,\hD_\bm)$ is a pair of a
relative automorphism $\si$ of $\hY_\bm$ and an automorphism $h$ of
$X$ so that
$$
u\circ h=\sigma\circ u.
$$

\begin{defi}
A {\em relative morphism} to $\hatYrel$ is an ordinary relative
morphism to $(\hY_\bm,\hD_\bm)$ for some $\bm$; it is {\em stable}
if its automorphism group is finite.
\end{defi}

Note that an effective class $(\vd,\vmu)$ of an FTCY graph $\Ga$ can
also be viewed as an effective class of any degeneration of $\Ga$,
and in particular, an effective class of $\Ga_\bm$. We fix a FTCY
graph $\Ga$, an effective class $(\vd,\vmu)$, and an even integer
$\chi$. We then form the moduli space $\Mdmu$ of all stable relative
morphisms $u$ to $\Yr$ that satisfy
\begin{itemize}
\item $\chi(\cO_X)=\chi/2$, where $X$ is the domain
      curve of $u$;
\item the associated effective class of $u$ is $(\vd,\vmu)$.
\end{itemize}

Since $\hY$ is a formal Calabi-Yau threefold with possibly normal
crossing singularity and smooth singular loci, the moduli space
$\Mdmu$ is a formal Deligne-Mumford stack with a perfect obstruction
theory \cite{Li1, Li2}.

\begin{lemm}
The virtual dimension of $\Mdmu$ is $\sum_{v\in
V_1(\Gamma)}\ell(\mu^v)$.
\end{lemm}

\begin{proof}
The proof is straightforward and will be omitted.
\end{proof}

\subsection{Equivariant Degeneration}\label{sec:T-degenerate}
Let $T$ act on $\Po\times\Ao$ by
$$
t \cdot ([X_0,X_1], s)= ([\fp(t)X_0,X_1], s),
$$
where $\fp\in\Lambda_T=\Hom(T,\CC^*)$. Let $\fY$ be the blowup of
$\Po\times \Ao$ at $([0,1],0)$. The $T$-action on $\Po\times \Ao$
can be lifted to $\fY$ such that the projection $\fY\to\Po\times
\Ao$ is $T$-equivariant; composition with the projection $\Po\times
\Ao\to \Ao$ gives a $T$-equivariant family of curves
$$
\fY\lra \Ao
$$
such that $\fY_s\cong\Po$ for $s\neq 0$ and $\fY_0\cong
\Po\sqcup\Po$.

The above construction can be generalized as follows. Let $\Ga$ be a
FTCY graph and let
$$
V_2(\Ga)=\{v_1,\ldots,v_n\}
$$

Then we have a $T$-equivariant family
\begin{equation}\label{eqn:Ftarget}
(\hat{\cY},\hat{\cD}) \to \A^n
\end{equation}
such that
$$
(\hat{\cY},\hat{\cD})_{\mathbf{0}}=\hat{\cY}_{(0,\ldots,0)}\cong \Yr
\quad \textup{and} \quad
(\hat{\cY},\hat{\cD})_\bs=(\hat{\cY},\hat{\cD}_{(s_1,\ldots,s_n)})
\cong \hY_{\Gamma_{ \{ v_i\mid s_i\neq 0\} } }\urel.
$$
Recall that given a subset $A\subset V_2(\Ga)$,
$\Ga_A$ is the smoothing of $\Ga$ along $A$ 
(Section \ref{sec:operation}). The $T$-action on  $\A^n$ is
trivial, and the $T$-action on each fiber is consistent with the one
described in Section \ref{sec:FTCY}.

By the construction in \cite{Li2}, there is a $T$-equivariant family
\begin{equation}\label{eqn:Fmoduli}
\MY\to \A^n
\end{equation}
such that $\MY_\bs=\cM^\bu_{\chi,\vd,\vmu}(\hat{\cY}_\bs)$. In
particular,
$$
\MY_{\mathbf{0}}=\Mdmu.
$$
The total space $\MY$ is a formal Deligne-Mumford stack with a
perfect obstruction theory $[\bT^1\to \bT^2]$ of virtual dimension
$$
\sum_{v\in V_1(\Ga)}\ell(\mu^v) + |V_2(\Ga)|.
$$
For each $v\in V_2(\Ga)$ there is line bundle $\bL^v$ over $\MY$
with a section $s^v: \MY\lra \bL^v$ such that
$$
\Mdmu=\MY_{\mathbf{0}}
$$
is the zero locus
$$
\{s^v=0 \mid v\in V_2(\Ga)\}\subset \MY.
$$
The pair $(\bL^v, s^v)$ corresponds to $(\bL_0, \mathbf{r}_0)$ in
\cite[Section 3]{Li2}.

\subsection{Perfect Obstruction Theory}\label{sec:tan-obs}

Let $\Ga$ be a FTCY graph, and let $(\vd,\vmu)$ be an effective
class of $\Ga$. We briefly describe the perfect obstruction theory
on $\Mdmu$ constructed in \cite{Li2}.

Let $\MY\to \A^{|V_2(\Ga)|}$, $[\bT^1\to\bT^2]$, and $\{ \bL^v\mid
v\in V_2(\Ga)\}$ be defined as in Section \ref{sec:T-degenerate}.
Let $[\tilde{\cT}^1\to \tilde{\cT}^2]$ be the perfect obstruction
theory on $\Mdmu$.  Let $u: (X,\bq) \longrightarrow
(\hY_\bm,\hD_\bm)$ represent a point in $\Mdmu\subset \MY$, where
$$
\bq=\{ q^v_j\mid v\in V_1(\Ga), 1\leq j\leq  \ell(\mu^v) \}.
$$
We have the following exact sequence of vector spaces at $u$:
\begin{equation}\label{eqn:TTL}
0\lra \tilde{\cT}_u^1\lra \bT^1_u \lra \bigoplus_{v\in
V_2(\Ga)}\bL_u^v \lra \tilde{\cT}_u^2 \lra \bT^2_u \lra 0.
\end{equation}
We will describe $\bT^1_u$, $\bT^2_u$, and $\bL_u^v$ explicitly.
When $\Ga$ is a regular FTCY graph, i.e., $V_2(\Ga)=\emptyset$, the
line bundles $\bL^v$ do not arise, and $\MY=\Mdmu$.

We first introduce some notation. Given $\bm: V_1(\Ga)\cup V_2(\Ga)
\to \ZZ_{\geq 0}$, let $\eGa^v_{m^v}$ be the flat chain of length
$m^v=\bm(v)$ associated to $v\in V_1(\Ga)\cup V_2(\Ga)$, and let
$$
V(\eGa^v_\bm)=\{ \bar{v}^v_0,\ldots, \bar{v}^v_{m^v} \},
$$
where $\bar{v}^v_{m^v}\in V_1(\Ga_\bm)$ if $v\in V_1(\Ga)$.

Let $v\in V_1(\Ga)$ and $0\leq l\leq m^v-1$, or let $v\in V_2(\Ga)$
and  $0\leq l\leq m^v$. We define a line bundle $L^v_l$ on the
divisor $\hD^v_l=\Si(\bar{v}_l)$ in $Y_\bm$ by
$$
L^v_l=N_{\hD^v_l/\Si(e_v)}\otimes N_{\hD^v_l/\Si(e_v')}
$$
where $\fv_0^{-1}(\bar{v}^v_l)=\{e_v,e_v'\}$. Note that $L^v_l$ is a
trivial line bundle on $\hD^v_l$.

With the above notation, we have
\begin{equation}\label{eqn:Lv}
\bL^v_u=\bigotimes_{l=0}^{m^v} H^0(\hD^v_l, L^v_l).
\end{equation}
The tangent space $\bT^1_u$ and the obstruction space $\bT^2_u$ to
$\MY$ at the moduli point
$$
[u:(X,\bq)\lra (\hY_\bm,\hD_\bm)]
$$
are given by the following two exact sequences:
\begin{equation}\label{eqn:exactI}
0 \lra \Ext^0(\Omega_X(R_\bq), \cO_X) \lra H^0(\mathbf{D}^\bu)\lra
\bT^1_u  \lra
\end{equation}
\begin{equation*}\label{eqn:exactI}
 \qquad\qquad \lra \Ext^1(\Omega_X(R_\bq), \cO_X) \lra H^1(\mathbf{D}^\bu)\lra
  \bT^2_u \lra  0
\end{equation*}
and
\begin{equation*}
0 \lra H^0\bigl(\TY\bigr)\lra
H^0(\mathbf{D}^\bu)\lra\qquad\qquad\qquad\qquad\qquad\qquad\qquad
\end{equation*}
\begin{equation}\label{eqn:exactII}
\lra\Oplus_{\scriptstyle v\in V_1(\Ga)\atop\scriptstyle 0\leq l\leq
m^v-1} \!\!\! \zero \oplus \Oplus_{\scriptstyle v\in
V_2(\Ga)\atop\scriptstyle 0\leq l\leq m^v}\!\!\zero \lra
H^1\bigl(\TY\bigr) \lra
\end{equation}
\begin{equation*}
\qquad\qquad\qquad\qquad\qquad \lra H^1(\mathbf{D}^\bu) \lra\!\!
\Oplus_{\scriptstyle v\in V_1(\Ga)\atop\scriptstyle 0\leq l\leq
m^v-1}\!\! \one \oplus \Oplus_{\scriptstyle v \in V_2(\Ga)\atop
\scriptstyle 0\leq l\leq m^v}\!\!\one\lra0.
\end{equation*}
where
$$
R_\bq=\sum_{v\in V_1(\Ga)}\sum_{j=1}^{\ell(\mu^v)}q_j^v,
$$\begin{equation}
\zero \cong  \bigoplus_{q\in u^{-1}(\hD^v_l) } T_q
(u^{-1}(\Si(e_v)))\otimes T^*_q (u^{-1}(\Si(e_v'))) \cong
\CC^{\oplus n^v_l}
\end{equation}
for
$\fv_0^{-1}(\bar{v}^v_l)=\{e_v,e_v'\}$,
\begin{equation}\label{eqn:one}
\one \cong \left. H^0(\hD^v_l,L^v_l)^{\oplus n^v_l}\right/
H^0(\hD^v_l, L^v_l),
\end{equation}
and $n^v_l$ is the number of nodes over $\hD^v_l$. In
(\ref{eqn:one}),
$$
H^0(\hD^v_l,L^v_l)\lra H^0(\hD^v_l,L^e_l)^{\oplus n^v_l}
$$
is the diagonal embedding.

We refer the reader to  \cite{Li2} for the definitions of
$H^i(\mathrm{D}^\bu)$  and the maps between terms in
(\ref{eqn:exactI}), (\ref{eqn:exactII}).

\subsection{Formal Relative Gromov-Witten Invariants}\label{sec:formalGW}
Usually, the relative Gromov-Witten invariants are defined as
integrations of the pull back classes from the target and the
relative divisor. In the case studied, the analogue is to integrate
a total degree $2\sum_{v\in V_1(\Gamma)}\ell(\mu^v)$ class from the
relative divisor $\hD$. The class we choose is the product of the
``Poincar\'{e} dual'' of the divisor $\hL^v\sub \hD^v$, one for each
marked point $q^v_i$. Equivalently, we consider the moduli space
$$
\Mfml=\left\{(u,X,\{q_j^v\})\in\Mdmu\mid u(q_j^v)\in \hL^v\right\}.
$$
Its virtual dimension is zero. More precisely, let $[\cT^1\to
\cT^2]$ be the perfect obstruction theory on $\Mfml$, and let
$[\tilde{\cT}^1\to \tilde{\cT}^2]$ be the perfection obstruction
theory on $\Mdmu$. Given a moduli point
$$
[u:(X,\bx)\to (\hY_\bm,\hD_\bm) ]\in \Mfml\subset \Mdmu,
$$
we have
\begin{equation}\label{eqn:TTN}
\cT^1_u -\cT^2_u = \tilde{\cT}^1_u -\tilde{\cT}^2_u- \bigoplus_{v\in
V_1(\Ga)}\bigoplus_{j=1}^{\ell(\mu^v)} (N_{\hL^v/\hD^v})_{u(q^v_j)}
\end{equation}
as virtual vector spaces.

In the rest of this subsection (Section \ref{sec:formalGW}), we fix
$\chi,\Ga,\vd,\vmu$, and write $\cM$ instead of $\Mfml$. We now
define the formal relative Gromov-Witten invariants of $\Yr$ by
applying the virtual localization to the moduli scheme $\cM$. We use
the equivariant intersection theory developed in \cite{Edi-Gra} and
the localization in \cite{Edi-Gra2, Gra-Pan}.

Since $\Yr$ is toric, the moduli space $\cM$ and its obstruction
theory are $T$-equivariant. We consider the fixed loci $\cM^T$ of
the $T$-action on $\cM$. Its coarse moduli space is projective. The
virtual localization is an integration of the quotient equivariant
Euler classes. When $[u]$ varies in a connected component of
$\cM^T$, the vector spaces $\cT^1_u$ and $\cT^2_u$ form two vector
bundles. We denoted them by $\cT^1$ and $\cT^2$. Since the
obstruction theory are $T$-equivariant, both $\cT^i$ are
$T$-equivariant. We let $\cT^{i,f}$ and $\cT^{i,m}$ be the fixed and
the moving parts of $\cT^i$. Since the fixed part $\cT^{i,f}$ induces a
perfect obstruction theory of $\cM^T$, it defines a virtual cycle
$$
[\cM^T]\virt\in A_*(\cM^T),
$$
where $A_*(\cM^T)$ is the Chow group with rational coefficients.

The perfect obstruction theory  $[\cT^{1,f}\to \cT^{2,f}]$ together
with the trivial $T$-action defines a $T$-equivariant virtual cycle
$$
[\cM^T]\virtt \in A_*^T(\cM^T).
$$
Since $T$ acts on $\cM^T$ trivially, we have
\begin{equation}\label{eqn:ATA}
A_*^T(\cM^T)\cong A_*(\cM)\otimes \Lambda_T
\end{equation}
where $\Lambda_T=\Hom(T,\CC^*)\cong A_*^T(\mathrm{pt})\cong
\QQ[u_1,u_2]$. Under the isomorphism (\ref{eqn:ATA}), we have
$$
[\cM^T]\virtt = [\cM]\virt \otimes 1.
$$

The moving part $\cT^{i,m}$ is the virtual normal bundles of
$\Mfml^T$. Let
$$
e^T(\cT^{i,m})\in A^*_T(\cM^T)
$$
be the $T$-equivariant Euler class of $\cT^{i,m}$, where
$A^*_T(\cM^T)$ is the  $T$-equivariant operational Chow group (see
\cite[Section 2.6]{Edi-Gra}). For $i=1,2$, $e^T(\cT^{i,m})$ lies in
the subring
$$
A^*(\cM^T)\otimes \QQ[u_1,u_2] \subset A^*_T(\cM^T)
$$
and is invertible in
$$
A^*(\cM^T)\otimes \QQ[u_1,u_2]_\fm \subset A^*_T(\cM^T)\otimes
\QQ[u_1,u_2]_\fm.
$$
where $\QQ[u_1,u_2]_{\fm}$ is $\QQ[u_1,u_2]$ localized at the 
ideal $\fm=(u_1,u_2)$.

For later convenience, we introduce some notation. Let $X$ be a
Deligne-Mumford stack with a $T$-action, and let $X^T$ be the
$T$-fixed points. Recall that
\begin{equation}\label{eqn:localize}
A_*^T(X)\otimes \QQ[u_1,u_2]_\fm \cong A_*^T(X^T)\otimes
\QQ[u_1,u_2]_\fm \cong A_*(X^T)\otimes \QQ[u_1,u_2]_\fm.
\end{equation}
The degree of a zero cycle defines a map  $\deg:A_0(X^T)\to \QQ$. We
define
$$
\deg_\fm: A_d(X^T)\otimes \QQ[u_1,u_2]_\fm  \to \QQ[u_1,u_2]_\fm
$$
by
$$
a\otimes b \mapsto \left\{ \begin{array}{ll}
\deg (a)b & d=0,\\
0 & d\neq 0.
\end{array}\right.
$$
This gives a ring homomorphism
$$
\deg_\fm: A_*^T(X)\otimes \QQ[u_1,u_2]_\fm \cong A_*(X^T)\otimes
\QQ[u_1,u_2]_\fm \to \QQ[u_1,u_2]_\fm.
$$
Given $c\in A^*_T(X)\otimes \QQ[u_1,u_2]_\fm$ and $\alpha\in
A_*^T(X)\otimes \QQ[u_1,u_2]_\fm$, define
$$
\int_{\alpha} c= \deg_\fm(c\cap \alpha)\in \QQ[u_1,u_2]_\fm.
$$

Following the lead of the virtual localization formula
\cite{Gra-Pan}, we define
\begin{defi}[formal relative Gromov-Witten invariants]\label{df:Fuu}
\begin{equation}\label{eqn:FN}
F^{\bu\Gamma}_{\chi,\vd,\vmu}(u_1,u_2)=\frac{1}{|\Aut(\vmu)|}
\int_{[\cM^T]\virtt}\frac{e^T(\cT^{1,m})}{e^T(\cT^{2,m})}
\end{equation}
where we view $[\cM^T]$ as an element in $A_*^T(\cM^T)\otimes
\QQ[u_1,u_2]_\fm$.
\end{defi}
Note that
$$
\frac{e^T(\cT^{1,m})}{e^T(\cT^{2,m})} \cap [\cM^T]\virtt\in \bbl
A_*^T(\cM^T)\otimes \QQ[u_1,u_2]_\fm \bbr_0
$$
where $\bbl A_*^T(\cM^T)\otimes \QQ[u_1,u_2]_\fm \bbr_0$ is the
degree zero part of the graded ring  $A_*^T(\cM^T)\otimes
\QQ[u_1,u_2]_\fm$. Therefore,
$$
F^{\bu\Ga}_{\chi,\vd,\vmu}(u_1,u_2)\in
(\QQ[u_1,u_2]_\fm)^0=\QQ(u_1/u_1)
$$
where $(\QQ[u_1,u_2]_\fm)^0$ is the degree zero part of the graded
ring $\QQ[u_1,u_2]_\fm$.

\begin{rema}
For our purpose of defining $F^{\bu\Ga}_{\chi,\vd,\vmu}(u_1,u_2)$,
we may consider the equivariant Borel-Moore homology
$H_*^T(\cM)=H_*^{T_\RR}(\cM)$ instead of the equivariant Chow group
$A_*^T(\cM)$, and consider the equivariant cohomology
$H^*_T(\cM)=H^*_{T_\RR}(\cM)$ instead of the equivariant operational
Chow group $A^*_T(\cM)$, where $\cM$ is any of the moduli spaces
involved in the above discussions.
\end{rema}

If $\cM$ were a proper Deligne-Mumford stack then its obstruction
theory would define a virtual cycle
\begin{equation}\label{eqn:false}
[\cM]\virt\in A_0(\cM)
\end{equation}
and
$$
F^{\bu\Ga}_{\chi,\vd,\vmu}(u_1,u_2)=\frac{1}{|\Aut(\vmu)|}\deg[\cM]\virt\in
\QQ
$$
would be a topological invariant independent of $u_1,u_2$. However,
$\cM$ is not proper, so (\ref{eqn:false}) does not exist.
Nevertheless, we will show that
\begin{theo}\label{thm:R2}
The function $F^{\bu\Gamma}_{\chi,\vd,\vmu}(u_1,u_2)$ is independent
of $u_1, u_2$; hence is a rational number depending only on
$\Gamma$, $\chi$, $\vd$ and $\vmu$.
\end{theo}

In Section \ref{sec:hodge} and Section \ref{sec:glue}, we will
reduce the invariance of  $F^{\bu\Ga}_{\chi,\vd,\vmu}(u_1,u_2)$
(Theorem \ref{thm:R2}) to the invariance for a special topological
vertex (Theorem \ref{thm:tri-inv}).

%% file: sec5.tex
\section{Invariance of the Topological Vertex}\label{sec:tri-inv}

We begin with the notion of topological vertex and topological
vertex with standard framing.

\begin{defi}[topological vertex and standard framing]\label{def:vertex}
A topological vertex is a FTCY graph that has one trivalent vertex
and three univalent vertices (see Figure 10 in Section
\ref{sec:hodge}). We say a topological vertex has a
standard framing if its three edges $e_1$, $e_2$ and $e_3$ that
share $v_0$ as their initial vertices have their position and
framing maps satisfying (see Figure 7)
$$
\ff(e_1)=\fp(e_2),\quad \ff(e_2)=\fp(e_3) \and \ff(e_3)=\fp(e_1).
$$
\end{defi}

\begin{figure}[h]\label{figure7}
\begin{center}
\psfrag{v0}{\footnotesize $v_0$} 
\psfrag{v1}{\footnotesize $v_1$}
\psfrag{v2}{\footnotesize $v_2$} 
\psfrag{v3}{\footnotesize $v_3$}
\psfrag{p1}{\footnotesize $\fp(e_1)$}
\psfrag{p2}{\footnotesize $\fp(e_2)$}
\psfrag{p3}{\footnotesize $\fp(e_3)$}
\psfrag{f1}{\footnotesize $\ff(e_1)$}
\psfrag{f2}{\footnotesize $\ff(e_2)$}
\psfrag{f3}{\footnotesize $\ff(e_3)$}
\includegraphics[scale=0.5]{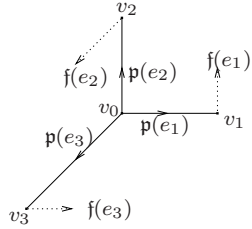}
\end{center}
\caption{A Topological vertex with standard framing}
\end{figure}

In this section, we shall prove

\begin{theo}[invariance of the topological vertex] \label{thm:tri-inv}
Theorem \ref{thm:R2} holds for any topological vertex with standard
framing.
\end{theo}

We fix a topological vertex with standard framing $\Gamma$ once and
for all in this section; we let $\Yr=(\hat Y,\hat D)$ be its
associated FTCY threefold. As before, we continue to denote $T$ the
group $(\CC\sta)^2$. We abbreviate $\Mfml$ to $\Mdy$.

Our first step to prove the invariance of
\begin{equation}
\lab{A}F^{\bu\Gamma}_{\chi,\vd,\vmu}(u_1,u_2)=\frac{1}{|\Aut(\vmu)|}
\int_{[\Mdy^\te]\virt}\frac{e^T(\cT^{1,m})}{e^T(\cT^{2,m})} \in
\QQ(u_1/u_2)
\end{equation}
is to construct a new pair of a nonsingular projective toric
threefold $W$ with a relative divisor $D\subset W$ and a subdivisor
$L\sub D$ so that the similarly defined moduli $\Mw$ (abbreviated to
$\Mdw$) of relative stable morphisms with boundary constraint $L$
admits a morphism
$$
\Phi: \Mdy \  \bl=\Mfml\br    \lra\Mdw\ \ \bl=\Mw\br
$$
so that the induced map on the $T$-fixed loci
$$
\Phi^T: \Mdy^T\lra \Mdw^T
$$
has the property that it is an open and closed embedding and the
obstruction theories of $\Mdw$ along its fixed loci is identical to
that of $\Mdy$ via $\Phi^T$.
Because the equality of two obstruction theories, we have the
identity:
\begin{equation}\lab{I1} 
\int_{[\Mdy^\te]\virt}\frac{e^T(\cT^{1,m})}{e^T(\cT^{2,m})}
=\int_{[\Mdw^\te_\Phi]\virt}\frac{e^T(\cT^{1,m})}{e^T(\cT^{2,m})},
\end{equation}
where $\Mdw^T_\Phi$ is the image of the fixed loci of $\Mdy$. (Here
by abuse of notation, we denote by $\cT^{i,m}$ the moving parts of
the obstruction complex $[\cT^1\to\cT^2]$ of $\Mdy$ as well as
$\Mdw$ along their fixed loci.)

To prove the invariance of the right hand side, we shall devise a
local contribution of $\deg[\Mdw]\virt$ along $\Mdy$; this local
contribution is a sum of the desired term (\ref{A}) with {\sl some
other terms}; we will show that this {\sl some other terms} vanish
completely. This will settle the invariance of the topological
vertex $\Ga$.

\subsection{The Relative Calabi-Yau Manifold $W\urel$ and the Morphism $\Phi$}
\label{sec:targetW}

We begin with constructing the toric variety $W\urel$ as promised.
Looking at the graph $\Gamma$ that we chose, the obvious choice of
$W$ is the toric blowup of $\Po\times\Po\times\Po$ along three
disjoint lines
\begin{equation}\label{3.7}
\ell_1=\infty\times \Po \times 0,\qquad \ell_2=0\times \infty \times
\Po \and \ell_3=\Po\times 0 \times \infty.
\end{equation}

The moment polytope of $W$, which is the image of the moment map
$$
\Upsilon: W\lra \RR^3
$$
of the $U(1)^3$-action on $W$, is shown in Figure 8. It is
diffeomorphic to the quotient $W/U(1)^3$. Here we follow the
convention that $(z_1,z_2,z_3)$ is the point
$([z_1,1],[z_2,1],[z_3,1])$ in $(\Po)^{3}$. We let $D\sub W$ be the
exceptional divisor and let $D_i\sub D$ be its connected component
lying over $\ell_i$. Each $D_i$ is isomorphic to $\Po\times\Po$.
We then let $C_1, C_2$ and $C_3$ be the proper transforms of
$$
\Po \times 0 \times 0,\qquad 0 \times \Po \times 0 \and 0\times
0\times \Po,
$$
and let $L_i\sub D_i$, $i=1,2$ and $3$, be the preimage of
$$
(\infty,0,0)\in\ell_1, \qquad (0,\infty,0)\in\ell_2 \and
(0,0,\infty)\in\ell_3.
$$
Clearly, restricting to $C_i$ the log-canonical sheaf
\begin{equation}\lab{10.1} 
\wedge^3 \Omega_W(\log D)|_{C_i}\cong\cO_{C_i}.
\end{equation}
Hence to the curves $C_i$ the relative pair $W\urel=(W,D)$ is
practically a relative Calabi-Yau threefold.

\begin{figure}[h] \label{figure 8}
\begin{center}
\psfrag{z1}{\footnotesize $z_1$} 
\psfrag{z2}{\footnotesize $z_2$}
\psfrag{z3}{\footnotesize $z_3$} 
\psfrag{p0}{\footnotesize $p_0$}
\psfrag{p1}{\footnotesize $p_1$} 
\psfrag{p2}{\footnotesize $p_2$}
\psfrag{p3}{\footnotesize $p_3$} 
\psfrag{q1}{\footnotesize $q_1$}
\psfrag{q2}{\footnotesize $q_2$} 
\psfrag{q3}{\footnotesize $q_3$}
\psfrag{L1}{\scriptsize $L_1$} 
\psfrag{L2}{\scriptsize $L_2$}
\psfrag{L3}{\scriptsize $L_3$} 
\psfrag{D1}{\scriptsize $D_1$}
\psfrag{D2}{\scriptsize $D_2$} 
\psfrag{D3}{\scriptsize $D_3$}
\psfrag{C1}{\scriptsize $C_1$} 
\psfrag{C2}{\scriptsize $C_2$}
\psfrag{C3}{\scriptsize $C_3$}
\includegraphics[scale=0.4]{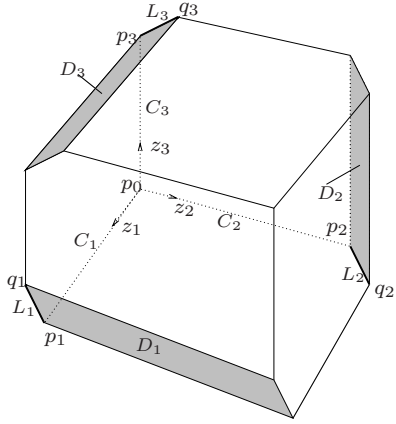}
\end{center}
\caption{Moment Polytope of $W$.}
\flushleft{
\begin{quote} 
\small{All faces of this polytope 
represent the $(\CC^*)^3$ invariant
divisors of $W$. The point ${p}_0$ is the image of the point
$(0,0,0)\in W$. The line $\overline{p_0 p_i}$ is the image of the
curve $C_i\cong\Po$, and the thickened line $\overline{p_i q_i}$ is the image
of the curve $L_i\cong \Po$. The rectangle face containing the edge
$\overline{p_i q_i}$ is the image of the relative divisor $D_i\cong\Po\times\Po$.} 
\end{quote}}
\end{figure}

For later discussion, we agree that under the isomorphisms
$D_i\cong\Po\times\Po$ and $\ell_i\cong\Po$, the tautological
projection $D_i\to \ell_i$ is the first projection. Under this
convention, the line $L_i\sub D_i$ is the line $0\times\Po$ and the
intersection $p_i=C_i\cap D_i$ is the point $(0,0)$.

As to the torus action, we pick the obvious one on $(\Po)^3$ via
\begin{equation}\lab{11.1}
(z_1,z_2,z_3)^{\three{t}}= (t_1 z_1,t_2 z_2,t_3z_3),\qquad
\three{t}\in (\CC^*)^3.
\end{equation}
It lifts to a $(\CC^*)^3$-action on $W$ that leaves $D_i$ and $L_i$
invariant. We let $T\sub(\CC^*)^3$ be the subgroup defined by
$t_1t_2t_3=1$; it is isomorphic to $(\CC\sta)^2$ and is the subgroup
that leaves (\ref{10.1}) invariant. In the following, we shall
view $W\urel=(W,D)$ as a relative Calabi-Yau $T$-manifold to the
curves $C_i$.

Next we will define the moduli space $\Mw$. Clearly, each $C_i$
induces a homology class $[C_i]\in H_2(W;\ZZ)$. For
$$
\vmu=(\up{\mu})\in \cP_+^3,
$$
we let $\vd$ be the homology class
$$
\vd=|\mu^1|[C_1]+|\mu^2|[C_2]+ |\mu^3|[C_3]\in H_2(W;\ZZ).
$$
The pair $(\vd,\vmu)$ is an effective class of $\Gamma$:
$$
\vd(\ee_i)=|\mu^i|,\ \ \vmu(v_i)=\mu^i,\ \ i=1,2,3.
$$
We then let
$$
\Mw,\quad \text{abbreviate to}\ \Mdw,
$$
be the moduli of relative stable morphisms
$$
u:(X;R_1,R_2,R_3)\lra W\urel=(W,D_1,D_2,D_3)
$$
having fundamental classes $\vd$, having ramification patterns
$\mu^i$ along $D_i$, and satisfying $u(R_i)\sub L_i$, modulo the
equivalence relation introduced in \cite{Li2}. It is a proper,
separated DM-stack; it has a perfect obstruction theory \cite{Li1,
Li2}, and thus admits a virtual cycle.

It follows from our construction that the scheme $Y$, which is the
the closure of the three one-dimensional orbits in $\hat Y$, can be
identified with the union $C_1\cup C_2\cup C_3$ in $W$; the formal
scheme $\hat Y$ is then the formal completion of $W$ along $Y$.
Further, the relative divisor $\hat D$ of $\hat Y$ (resp. teh
subdivisor $\hat L\sub\hat D$) is the preimage of the relative
divisor $D\sub W$ (resp. the subdivisor $L\sub D$); the induced
morphism
\begin{equation}\lab{10.2}
\phi: (\hat Y,\hat D, \hat L)\lra (W,D,L)
\end{equation}
is $T$-equivariant; and the two effective classes $(\vd,\vmu)$ are
consistent under the map $\phi$. Therefore, it induces a
$T$-equivariant morphism of the moduli spaces
\begin{equation}\lab{10.3}
\Phi:\Mdy\lra\Mdw,
\end{equation}
which induces a morphism
$$
\Phi^\te: \Mdy^\te\lra \Mdw^\te
$$
between their respective fixed loci.

\begin{lemm}
The morphism $\Phi^\te$ is an open and closed embedding; the
obstruction theories of $\Mdy$ and $\Mdw$ are identical under $\Phi$
along the fixed loci $\Mdy^\te$ and its image in $\Mdw$.
\end{lemm}

\begin{proof}
This follows immediately from that $C_1$, $C_2$ and $C_3$ are the
closures of three one-dimensional orbit, that $Y=C_1\cup C_2\cup
C_3$ and that $\hat Y$ is the formal completion of $W$ along $Y$.
\end{proof}



Later, we shall work with the moduli of relative stable morphisms
$\cM_{\chi,\vd,\vmu}^\bu(Y\urel,L)$, similarly defined as that of
$\Mdw$. Again, for notational simplicity, we shall abbreviate it to
$\Mdyy$.

\subsection{Invariant Relative Stable Morphisms}\label{sectionA2}

Let $a_1,a_2,a_3\in\ZZ$ with $a_1+a_2+a_3=0$ be three relatively
prime integers; $\eta=(a_1,a_2,a_3)$ defines a subgroup
$$
T_\eta=\{(t^{a_1},t^{a_2},t^{a_3})
\mid t\in U(1)\}\sub T.
$$
Our next task is to characterize those stable relative morphisms
that are invariant under $T_\eta\sub T$ and are small deformations
of elements in $\Mdyy$.

To investigate relative stable morphisms to $W$, we need the
expanded relative pair $(\wm,\dm)$, $\bm=(m_1,m_2,m_3)$
(see Figure 9). Let
$\Delta$ be the projective bundle $\mathbb
P\bl\cO_{\Po\times\Po}\oplus\cO_{\Po\times\Po}(0,1)\br$ with two
sections
$$
D_+=\mathbb P\bl\cO_{\Po\times\Po}\oplus0\br \and D_-=\mathbb
P\bl0\oplus\cO_{\Po\times\Po}(0,1)\br;
$$
we form an $m_i$-chain of $\Delta$ by gluing $m_i$ copies of
$\Delta$ via identifying the $D_-$ of one $\Delta$ to the $D_+$ of
the next $\Delta$ using the canonical isomorphism $\text{pr}\mh
D_\pm\to\Po\times\Po$; we then attach this chain to $D_i$ by
identifying the $D_+$ of the first $\Delta$ in the chain with $D_i$
and declaring the $D_-$ of the last $\Delta$ be $D[\bm]_i$. The
scheme $\wm$ is the result after attaching such three chains, of
length $m_1$, $m_2$ and $m_3$ respectively, to $D_1$, $D_2$ and
$D_3$ in $W$. The union
$$
\dm=\dm_1\cup\dm_2\cup\dm_3
$$
is the new relative divisor of $W[\bm]$. Note that our construction
is consistent with that the normal bundle of $D_i$ in $W$ has degree
$-1$ along $L_i$.

\begin{figure}[h]\label{figure9}
\begin{center}
\psfrag{z1}{\footnotesize $z_1$} \psfrag{z2}{\footnotesize $z_2$}
\psfrag{z3}{\footnotesize $z_3$} \psfrag{p0}{\footnotesize $p_0$}
\psfrag{p1}{\footnotesize $p_1$} \psfrag{p2}{\footnotesize $p_2$}
\psfrag{p3}{\footnotesize $p_3$} \psfrag{q1}{\footnotesize $q_1$}
\psfrag{q2}{\footnotesize $q_2$} \psfrag{q3}{\footnotesize $q_3$}
\psfrag{P2}{\footnotesize $p_2'$}
\psfrag{L1}{\scriptsize $L_1$} \psfrag{L2}{\scriptsize $L_2$}
\psfrag{L3}{\scriptsize $L_3$} \psfrag{D1}{\scriptsize $D_1$}
\psfrag{D2}{\scriptsize $D_2$} \psfrag{D3}{\scriptsize $D_3$}
\psfrag{C1}{\scriptsize $C_1$} \psfrag{C2}{\scriptsize $C_2$}
\psfrag{C3}{\scriptsize $C_3$} \psfrag{De}{\footnotesize $\Delta$}
\psfrag{Dm2}{\scriptsize $D[\bm]_2$}
\includegraphics[scale=0.4]{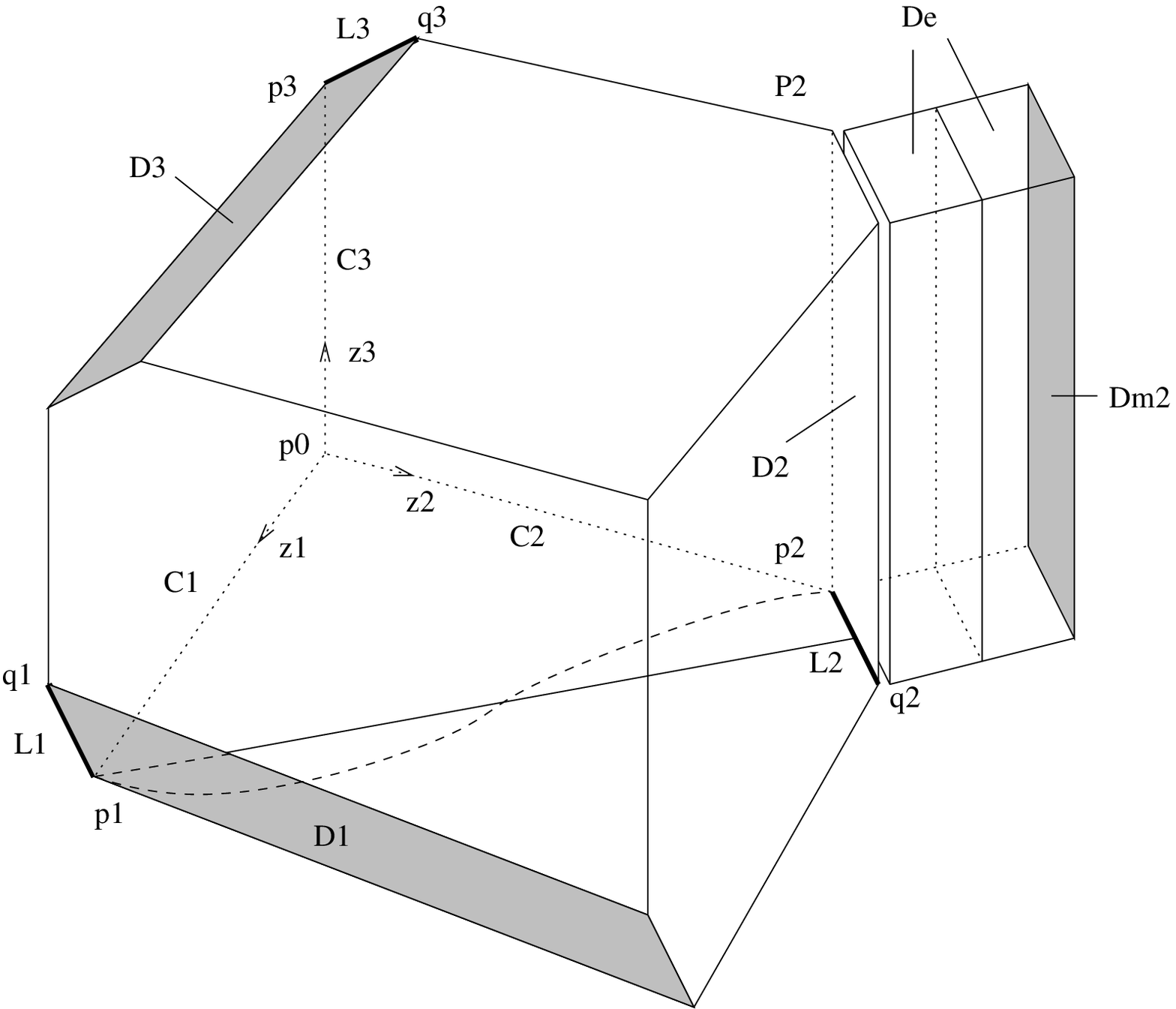}
\end{center}
\caption{ }
\flushleft{
\begin{quote}
\small{This is a sketch of the scheme 
$W[\bm]$ for $\bm=(0,2,0)$.
The main part is the moment polytope of $W$ (see Figure 8).
The added two solids to the right are the
two $\Delta$'s attached to $D_2$, resulting the scheme $W[\bm]$ with
$\bm=(0,2,0)$. The shaded faces are the relative divisor $D[\bm]$ of
$W[\bm]$.
The straight diagonal line contained in the bottom face
indicates the image of $\phi_{k,c}$ in case $\eta=(1,-1,0)$; the
curved line indicates the image of $\phi_{k,c}$ in the other case.}
\end{quote}}
\end{figure}

For future convenience, we denote by $\Delta[m_i]$ the chain of
$\Delta$'s that is attached to $D_i$; we denote by $L[\bm]_i\sub
D[\bm]_i$ the same line as $L_i\sub D_i$. The new scheme $W[\bm]$
contains $W$ as its main irreducible component; it also admits a
stable contraction $\wm\to W$. Also, $(\CC^*)^3$ acts on
$(W[\bm],D[\bm])$ since the $(\CC^*)^3$-action on $N_{D_i/ W}$
induces a $(\CC^*)^3$-action on each $\Delta$ attached to $D_i$.
Therefore $(\CC^*)^3$ and $T$ act on $\Mdw$. Unless otherwise
mentioned, the maps $W\to W[\bm]$ and $W[\bm]\to W$ are these
inclusion and projection; these maps are
$(\CC^*)^3$-equivariant.

The pair $(W,D)$ contains $(Y,p)$, $p=p_1+p_2+p_3$, as its subpair.
Accordingly, the pair $(\wm,\dm)$ contains a subpair
$(Y[\bm],p[\bm])$ whose main part $Y[\bm]$ is the preimage of $Y$
under the contraction $\wm\to W$. The relative divisor $p[\bm]$ is
the intersection $Y[\bm]\cap\dm$. It is the embedding
$Y[\bm]\sub\wm$ that induces the embedding $\Mdyy\sub\Mdw$.

Now let $u_0$ be a relative stable morphism in $\Mdyy$, considered
as an element in $\Mdw$; let $u_s$ be a small deformation of $u_0$
in $\Mdw^{T_\eta}$ that is not entirely contained in $\Mdyy$. Each
$u_s$ is a morphism from its domain $X_s$ to $W[\bm]$ for some
triple $\bm$ possibly depending on $s$. We let $\tilde u_s\mh X_s\to
W$ be the composite of $u_s$ with the contraction $\wm\to W$. Then
$\tilde u_s$ form a flat family of morphisms; it specializes to
$\tilde u_0$ as $s$ specializes to $0$.
Hence as sets, $\tilde u_s(X_s)$ specializes to $\tilde u_0(X_0)$
as $s$ specializes to $0$.
Because $\tilde u_s(X_s)$ are union of algebraic curves in $W$ and
$\tilde u_0(X_0)$ is contained in $C_1\cup C_2\cup C_3$, for general
$s$ the intersection $\tilde u_s(X_s)\cap D$ is discrete. Hence
every irreducible component $Z\sub\tilde u_s\upmo(D_i)$ must be
mapped to a fiber of $\Delta[m_i]/D_i$.

Now suppose there is such a connected component $Z$ with $u_s(Z)$
lies in the fiber of $\Delta[m_i]$ over $q\in D_i$, then by the
pre-deformable requirement on relative stable morphisms forces the
same $q$ in $D[\bm]_i$ to lie in $u_s(X_s)$. Because of the
requirement $u_s(Z)\cap\dm_i\sub L[\bm]_i$ we impose on the $\Mdw$,
we have
\begin{equation}\lab{10.7}
\tilde u_s(X_s)\cap D_i\sub L_i.
\end{equation}

Our primary interest is to those $u$ that are $\soe$-invariant and
are small deformations of elements in $\Mdyy$. This leads to the
following definition.

\begin{defi}\label{def:deform}
We let $\mwtdef$ be the union of all connected components of
$$\{[u,X]\in\Mdw^\soe\mid \tilde u(X)\cap D\ \text{is
finite} \}
$$
that intersect $\Mdyy$ but are not entirely contained in it.
\end{defi}

Following the discussion before Definition \ref{def:deform}, all $u$
in $\mwtdef$ satisfies (\ref{10.7}). In case $a_{i+1}\ne 0$ (we
agree $a_4=a_1$), the only $\soe$-fixed points of $L_i$ are $p_i$
and $q_i$; hence all $u$ in $\mwtdef$ satisfies a strengthened
version to (\ref{10.7}):
\begin{equation}\lab{10.71}
\tilde u_s(X)\cap D_i\sub p_i,\qquad\text{when \ } a_{i+1}\ne0.
\end{equation}
Here $q_i$ is ruled out because each connected component of
$\mwtdef$ intersects $\Mdyy$.

We now characterize elements in $\mwtdef$. We comment that we shall
reserve $a_1$, $a_2$ and $a_3$ for the three components of $\eta$;
we always assume the three $a_i$'s are relatively prime and that
$a_1+a_2+a_3=0$. In this and the next two Subsections, we shall
workout the case $a_1>0$ and $a_2,a_3<0$; the case $\eta=(1,-1,0)$
will be considered in Subsection \ref{later}. Now let
$[u,X]\in\mwtdef$ and let $V\sub\tilde u(X)$ be any irreducible
component. Since $u$ is $\soe$-invariant, $V$ is $\soe$-invariant.
Hence $V$ must be the lift of the set
$$\bar V=\{(c_1 t^{a_1}, c_2 t_2^{a_2}, c_3 t_3^{a_3})\mid
t\in\CC\cup\{\infty\}\ \}\sub(\Po)^3
$$
for some $(c_1,c_2,c_3)$. Thus we see that the following cases
are impossible:
\begin{enumerate}\item all
$c_i$ are non-zero: should this hold, then $\bar{V}\cap
\ell_2=(0,\infty,\infty)$, which violates \eqref{10.7}; \item
$c_1=0$ but the other two are non-zero: should this hold, then
either $V\cap D=V\cap D_1=q_1$, or $V\cap D=V\cap D_2=p_2'$ (see
Figure 9 for the location of $p_2'$), so \eqref{10.71} is violated;
\item $c_2=0$ but the others are non-zero: should this hold, then
$V\cap D_1=q_1$ since $a_1>|a_3|$, which violates \eqref{10.71}.
\end{enumerate}

This leaves us with the only two possibilities: when only one of
$c_i$ is non-zero or $c_3=0$ but the other two are non-zero. In the
first case we have $V=C_i$ for some $i$; in the later case $V$ is
the image of the map
\begin{equation}\lab{10.10}
\phi_{k,c}: \Po\lra W,\qquad k\in \ZZ^+,\ c\in\CC\sta
\end{equation}
that is the lift of $\Po\to(\Po)^{\times 3}$ defined by $\xi\mapsto
(\xi^{k a_1},c^{-k a_2}\xi^{k a_2},0)$. Clearly, $\phi_{k,c}$ is
$T_\eta$-invariant. It is easy to see that these are the only
$\soe$-equivariant maps $Z\to W$ from irreducible $Z$ whose images
are not entirely lie in $C_1\cup C_2\cup C_3$ and the divisor $D$.
This proves

\begin{lemm}\lab{10.6}
Suppose $a_1>0$ and $a_2$, $a_3<0$. Then any $(u, X)\in\mwtdef$
not entirely contained in $Y$ has at least one irreducible
component $Z\sub X$ and a pair $(k,c)$ so that $u|_{Z}\cong
\phi_{k,c}$.
\end{lemm}

Here by $u|_Z\cong\phi_{k,c}$ we mean that there is an isomorphism
$Z\cong\Po$ so that under this isomorphism $u|_Z\equiv\phi_{k,c}$.

When $c$ specialize to $0$, the map $\phi_{k,c}$ specializes to
$$\phi_{k,0}:\Po\sqcup\Po\lra W
$$
defined as follows. We endow the first copy (of $\Po\sqcup \Po$)
with the coordinate $\xi_1$ and the second copy with $\xi_2$; we
then form the nodal curve $\Po\sqcup\Po$ by identifying $0$ of the
first $\Po$ with $0$ of the second $\Po$; we define $\phi_{k,0}$ to
be the lift of the maps
$$
\xi_1\mapsto (\xi_1^{ka_1},0,0) \and \xi_2\mapsto
(0,\xi_2^{-ka_2},0).
$$
Since $\xi_1=0$ and $\xi_2=0$ are both mapped to the origin in
$(\Po)^{\times 3}$, they glue together to form a morphism
$\phi_{k,0}\mh\Po\sqcup \Po\to W$.

This leads to the following definition.
\begin{defi}
A deformable part of a $(u,X)\in\mwtdef$ consists of a curve $Z\sub
X$ and an isomorphism $u|_Z\cong\phi_{k,c}$ for some $(k,c)$.
\end{defi}

Suppose $(u,X)$ has at least two deformable parts, say
$(Z_1,\phi_{k_1,c_1})$ and $(Z_2,\phi_{k_1,c_2})$, then the explicit
expression of $\phi_{k,c}$ ensures that $Z_1$ and $Z_2$ share no
common irreducible components. Should $Z_1\cap Z_2\ne\emptyset$,
their intersection would be a nodal point of $X$ that could only be
mapped to either $D_1$ or $D_2$ of $W$ under $u$. (Note that it
could not be mapped to $p_0$ since then both $c_1$ and $c_2=0$, and
that node would be in more than two irreducible components of $X$.)
However, the case that the node is mapped to $D_1$ or $D_2$ can also
be ruled out because it violates the pre-deformable requirement of
relative stable morphisms \cite{Li1}. Hence $Z_1$ and $Z_2$ are
disjoint. This way, we can talk about the maximal collection of
deformable parts of $(u,X)$; let it be
$$
(Z_1,\phi_{k_1,c_1}),\cdots,(Z_l,\phi_{k_l,c_l}).
$$
For convenience, we order it so that $k_i$ is increasing.

\begin{defi} We define the deformation type of $(u,X)\in\mwtdef$
be
$$(k_i)_l=(k_1\leq k_2\leq\cdots\leq k_l)
$$
It defines a function on $\mwtdef$, called the deformation type
function.
\end{defi}

Let $(u,X)$ be an element in $\mwtdef$ of type $(k_i)_l$.
Intuitively, we should be able to deform $u$ within $\mwtdef$ by
varying $u|_{Z_i}$ using $\phi_{k,t}$ to generate an $\A^l$-family
in $\mwtdef$. It is our next goal to make this precise.

To proceed, we need to show how to put $\pkt$ into a family. We
first blow up $\Po\times\Ao$ at $(0,0)$ to form a family of curves
$\fY$ over $\Ao$. The complement of the exceptional divisor
$\fY-E=\Po\times\Ao-(0,0)$ comes with an induced coordinate
$(\xi,t)$. We define
$$\Phi_k|_{\fY-E}:\fY-E\lra W;\quad
(\xi,t)\mapsto (\xi^{ka_1},t^{-ka_2}\xi^{ka_2},0).
$$
We claim that $\Phi_k|_{\fY-E}$ extends to a $\Phi_k\mh \fY\to W$.
Indeed, if we pick a local coordinate chart near $E$, which is
$(\xi,v)$ with $t=\xi v$, then
$$\Phi_k|_{\fY-E}: (\xi,v)\mapsto (\xi^{ka_1},
(\xi v)^{-ka_2}\xi^{ka_2},0)= (\xi,v)\mapsto (\xi^{ka_1},
v^{-ka_2},0),
$$
which extends to a regular
$$
\Phi_k: \fY\lra W.
$$
Note that for $c\in\Ao$, the fiber of $(\Phi_k,\fY)$ over $c$ is
exactly the $\phi_{k,c}$ we defined earlier. Henceforth, we will
call $(\Phi_k,\fY)$ the standard model of the family $\pkt$; we will
use $\fY_c$ to denote the fiber of $\fY$ over $c\in\Ao$.

To deform $u$ using the family $\Phi_k$, we need to glue $\fY$ onto
the domain $X$. We let $\fD_1$ be the proper transform of
$0\times\Ao\sub\Po\times\Ao$ and let $\fD_2=\infty\times\Ao$ in
$\cY$. Both $\fD_1$ and $\fD_2$ are canonically isomorphic to $\Ao$
via the second projection. For $Z\sub X$, we fix an isomorphism
$Z\cong\fY_c$ so that $u|_Z\cong\phi_{k,c}$; we specify $v_1$,
$v_2\in Z$ so that $u(v_i)\in D_i$; we let $X_0$ be the closure of
$X-Z$ in $X$.

We now glue $\fY$ onto $X_0\times\Ao$. In case both $v_1$ and $v_2$
are nodes of $X$, we glue $\fY$ onto $X_0\times\Ao$ by identifying
$\fD_1$ with $v_1\times\Ao$ and $\fD_2$ with $v_2\times\Ao$, using
the their standard isomorphisms with $\Ao$; in case $v_1$ is a
marked point of $X$ and $v_2$ is a node, we glue $\fY$ onto
$X_0\times\Ao$ by identifying $\fD_2$ with $v_2\times\Ao$ and
declaring $\fD_1$ to be the new marked points, replacing $v_1$; in
case $v_1$ is a node and $v_2$ is a marked points, we repeat the
same procedure with the role of $v_1$ and $v_2$ and of $\fD_1$ and
$\fD_2$ exchanged; finally in case both $v_1$ and $v_2$ are marked
points, we simply replace $Z\times \Ao$ in $X\times\Ao$ by $\fY$
while declaring that $\fD_1$ and $\fD_2$ are the two marked points
replacing $v_1$ and $v_2$. We let $\cX\to\Ao$ be the resulting
family.

The morphisms
$$
X_0\times\Ao\twomapright{\rm{pr}} X_0\twomapright{u|_{X_0}}\wm
\and \Phi_k: \fY\twomapright{} W
$$
glue together to form a morphism
$$\cU:\cX\lra\wm.
$$
The pair $(\cU,\cX)$ is the family in $\mwtdef$ that keeps
$u|_{X_0}$ fixed.

More generally, we can deform $u$ inside $\mwtdef$ by identifying
and altering its restriction to the deformable parts of $X$
simultaneously. This way, any $u\in\mwtdef$ of type $(k_i)_l$
generates an $\A^l$ family of elements in $\mwtdef$.

\subsection{Global Structure of the Loci of Invariant Relative Morphisms}
\lab{global section}

In this subsection, we shall prove that any connected component of
$\mwtdef$ is an $\A^l$-bundle.

We begin with a technical Lemma

\begin{lemm}\lab{intact} Let $v\in X_{\text{node}}$ be a node in $u\upmo(D)$.
Then $v$ remains a node when $u$ deforms infinitesimally in
$\mwtdef$.
\end{lemm}

The key to the proof relies on the fact that the restrictions of the
automorphisms induced by $\soe$ on one irreducible component of $X$
that contains $v$ are infinite while on the other irreducible
component are finite.

Recall that there is an natural inclusion $h:T\to \Aut(\wm)$ induced
by the $T$-actions on $W$ and on $N_{D_i/ W}$. There is a unique
homomorphism $h':T\to \Aut(\wm/W)$ such that\footnote{The
automorphisms $\zeta\in\Aut(\wm)$ that preserve the fibers of the
map $\wm\to W$ are called relative automorphisms of $\wm/W$; the
group of all such automorphisms is denoted by $\Aut(\wm/W)$. If
$\bm=(m_1,m_2,m_3)$ then $\Aut(\wm/W) \cong (\CC^*)^{m_1+m_2+m_3}$.}
$h'(t)\circ h(t)\in \Aut(\wm)$ acts trivially on $p_i[m_i]$,  the
fiber of $\Delta[m_i]$ over $p_i\in D_i$, for all $t\in T$. Since
$u$ is $T_\eta$-invariant and the image
$u\bigl(u^{-1}(\Delta[m_i])\bigr)$ is entiredly contained in
$p_i[m_i]$, there are group homomorphisms $h_1\mh \soe\to\Aut(X)$
and $h_2\mh\soe\to\Aut(\wm/W)$ so that for all
$\sigma\in\soe$,
\begin{equation}\lab{auto}
\textup{(i) }  h'(\sigma)\circ h(\sigma)\circ u = h_2(\sigma)\circ
u\circ h_1(\sigma)\and \textup{(ii) } h_2(\si) \textup{ acts on
}p_i[m_i] \textup{ trivially}.
\end{equation}
Now let $v\in u\upmo(D_i)$ be a node of $X$ that is
mapped to $D_i$ under $u$; let $V_-$ be the irreducible component of
$X$ that contains $v$ that is mapped to $W$ and let $V_+$ be the
other irreducible component of $X$ that contains $v$. Then $u(V_+)$
must be contained in $\Delta[m_i]$. Since $h_1(\id)=\id$ and that
$\soe$ is connected, $h_1(\sigma)(V_\pm)\sub V_\pm$. Hence
$h_1(\sigma)$ are automorphisms of $V_\pm$ that fixed $v$. We let
$$
\soe|_{V_\pm} \stackrel{\mathrm{def}}{=}
\{h_1(\sigma)\mid\sigma\in\soe\}\sub\Aut(V_\pm,v);
$$
it is a group which is a homomorphism image of $T_\eta\cong U(1)$,
so it is either $U(1)$ or trivial.

\begin{lemm}\lab{lem:auto}
The group $\soe|_{V_-}$ is infinite while the group $\soe|_{V_+}$ is
finite. So $\soe|_{V_-}\cong U(1)$ and $\soe|_{V_+}$ is trivial.
\end{lemm}

\begin{proof}
The group $\soe|_{V_-}$ is infinite is obvious. Since $u(v)=p_i$ and
$u(V_-)\sub W$, $u(V_-)$ is $\soe$ invariant but not $\soe$ fixed.
Since the induced action on $u(V_-)$ is infinite, $\soe|_{V_-}$ must
be infinite because $u$ is $\soe$-invariant.

Note that $u(V_+)$ must be contained in $p_i[m_i]$ for some $i$. By
\eqref{auto},
$$
u|_{V_+}= u\circ h_1(\si)|_{V_+}: V_+\lra p_i[m_i]
\quad \forall \si\in \soe\, .
$$
By stability of the relative morphism $u$, $T_\eta|_{V_+}$ is
finite.
%
%
%
\end{proof}

We now prove Lemma \ref{intact}.
\begin{proof}[Proof of Lemma \ref{intact}]
Suppose $v$ can be smoothed at least of first order within
$\mwtdef$. Since $\mwtdef$ is fixed by $\soe$, $\soe|_{V_-}$ is
infinite forces $\soe|_{V_+}$ to be infinite as well. This violates
the assumption that $\soe|_{V_+}$ is finite. This proves the Lemma.
\end{proof}

\smallskip

As we argued before, each $u\in\mwtdef$ contains a deformable part
that is the union of some $\phi_{k_i,c_i}$. Our next task is to show
the deformable part of $u$ remain the same within a connected
component of $\mwtdef$.

We make it more precise now. We let $(u,X)$ be any element in
$\mwtdef$; let $Y_1,\cdots,Y_l\sub X$ be all its deformable parts so
that $u|_{Y_i}\cong\phi_{k_i,c_i}$; we let $v_{i,1}$ and $v_{i,2}\in
Y_i$ be the marked points so that $u(v_{i,j})=p_j$. Then according
to the discussion in the previous Subsection, by varying $u|_{Y_i}$
using $\Phi_{k_i}$ we get a copy $\Ao$ in $\mwtdef$; together they
provide a copy $\A^l$ in $\mwtdef$. This is one of the fiber of the
fiber bundle structure on $\mwtdef$ we are about to construct.

To extend this $\A^l\sub\mwtdef$ to nearby elements of $[u]$, we
need to extend all $Y_i$ in $X$ to a flat family of subcurves.

\begin{lemm}\lab{A9}
The deformation type function on $\mwtdef$ is locally constant.
\end{lemm}

\begin{proof}
We pick an analytic disk $0\in S$ and an analytic map
$\psi\mh S\to \Mdw^\soe\ldef$ so that $\psi(0)=[u]$. The morphism
$\psi$ pulls back the tautological family on $\mwtdef$ to a family
$\cU\mh \cX\to\cZ$ over $S$. The central fiber $\cX_0$ is $X$ and
thus contains $Y_i$. We let $\cN\sub\cX$ be the subscheme of the
nodes of all fibers of $\cX/S$. Since $v_{i,j}$ is either a marked
point or a node of $X$, $v_{i,j}\in\cN\cup\cR$. Let $\cP_{i,j}$ be
the connected component of $\cN\cup\cR$ that contains $v_{i,j}$. We
claim that $\cP_{i,j}$ is a section of $\cN\cup\cR\to S$. First,
$\cP_{i,j}$ is flat over $S$ at $v_{i,j}$. This is true in case
$v_{i,j}$ is a marked point since $\cR$ is flat over $S$ by
definition; in case $v_{i,j}$ is a node it is true because of Lemma
\ref{intact}. Therefore, $\cP_{i,j}$ dominates over $S$. Then
because $\cN\cup \cR$ is proper and unramified over $S$, dominating
over $S$ guarantees that $\cP_{i,j}$ is finite and \'etale over $S$.
But then since $S$ is a disk, $\cP_{i,j}$ must be isomorphic to $S$
via the projection.

We now pick the desired family of curves $\cY_i$. In case
$\cP_{i,j}$ is one of the section of the marked points of $\cX/S$,
we do nothing; otherwise, we resolve the singularity of the fibers
of $\cX$ along $\cP_{i,j}$. As a result, we obtain a flat family of
subcurves $\cY_i\sub\cX$ that contains $Y_i$ as its central fiber.
We let $\cU_i\mh \cY_i\to\cW$ be the restriction of $\cU$ to
$\cY_i$. Because $\cU_i(Y_i)\sub W\sub W[\bm]$, $\cU_i(\cY_i)\sub
W\times S\sub\cW$ as well.

Since $\cU\mh\cX\to\cW$ is a family of $\soe$-equivariant relative
stable maps, $\cU_i\mh \cY_i\to W$ is also a family of
$\soe$-equivariant stable morphisms. Then because $\cU_i|_{Y_i}$ is
isomorphic to $\phi_{k_i,c_i}$, each member of $\cU_i$ must be an
$\phi_{k_i,c}$ for some $c\in\CC$. This proves that the deformation
type of $\cU|_{\cX_s}$ contains that of $\cU|_{\cX_0}$ as a subset.
Because this holds true with $0$ and $s$ exchanged, it shows that
the deformation type function stay constant over $S$.

Finally, because any two elements in the same connected component of
$\mwtdef$ can be connected by a chain of analytic disks, the
deformation type function does take same values on such component.
This proves the lemma.
\end{proof}

We are now ready to exhibit a fiber bundle structure of any
connected component of $\mwtdef$. Let $\cQ\sub\mwtdef$ be any
connected component. According to the previous subsection, all
elements in $\cQ$ are of the same deformation types, say $(k_i)_l$.
In case $\cQ$ is not entirely contained in $\Mdyy$, then necessarily
$l>0$. To get the fiber structure, we need to take a finite cover of
$\cQ$, which we now construct.

\begin{defi}\lab{d10}
We define the groupoid $\bar\cQ$ over $\cQ$ as follows. For any
scheme $S$ over $\cQ$, we let $\bar\cQ(S)$ be the collection of data
$\{(\cU,\cX,\cW),\rho_i,\cZ_i,\pi_i\mid i=1\cdots,l\}$ of which
\begin{enumerate}
    \item $\cU\mh\cX\to\cW$ is an object\,\footnote{Here
    we consider $\cQ$ as a groupoid and $\cQ(S)$ is the collection
    of objects over $S$.} in $\cQ(S)$
    \item $\rho_i$ are morphisms from $S$ to $\Ao_i$, $\Ao_i\cong\Ao$;
    \item $\cZ_i$ are flat families of subcurves in
    $\cX$ over $S$ with all marked points discarded;
    \item $\pi_i\mh\cZ_i\to \rho_i\sta\fY_{k_i}$ an isomorphism over
    $S$;
\end{enumerate}
together they satisfy
$$
\cU|_{\cZ_i}\equiv\rho_i\sta\Phi_{k_i}\circ\pi_i: \cZ_i\lra W.
$$
Further, an arrow from $\{(\cU,\cX,\cW),\rho_i,\cZ_i,\pi_i\}$ to
$\{(\cU\pri,\cX\pri.\cW\pri),\rho_i\pri,\cZ_i\pri,\pi_i\pri\}$
consists of an isomorphism $h_1\mh \cX\to\cX\pri$ and an isomorphism
$h_2\mh \cW\to\cW\pri$ relative to $W$ so that under these
isomorphisms $\cZ_i=\cZ_i\pri$, $\rho_i=\rho_i\pri$,
$\pi_i=\pi_i\pri$ (for all $i$) and $\cU=\cU\pri$.
\end{defi}

Here we use $\Ao_i$ to denote the target of $\rho_i$, which is
$\Ao$. We are doing this to distinguish them for different $i$.

For a fixed type $(k_i)_l$, we let $G_{(k_i)}$ be the group
$\{\sigma\in S_l\mid k_{\sigma(i)}=k_i\}$.

\begin{prop}\lab{A10}
The groupoid $\bar\cQ$ is a DM-stack acted on by $G_{(k_i)}$ ; it is
finite and \'etale over $\cQ$, and $\bar\cQ/G_{(k_i)}=\cQ$.
The morphisms
$\rho_i$ in each object in $\bar\cQ$ glue to a morphism
$\bar\rho_i\mh \bar\cQ\to\Ao_i$. Let
$\bar\cQ_0=(\bar\rho_1,\cdots,\bar\rho_l)\upmo(0)$. Then there is a
canonical projection $\pi\mh \bar\cQ\to\bar\cQ_0$ making it an
$\A^l$-bundle over $\bar\cQ_0$. Finally, the morphism
$$(\pi,(\bar\rho_1,\cdots,\bar\rho_l)):
\bar\cQ\lra\bar\cQ_0\times\A^l
$$
is an isomorphism of DM-stacks.
\end{prop}

\begin{proof}
The proof is straightforward, following the previous discussion, and
will be omitted.
\end{proof}

\subsection{The Obstruction Sheaves}\lab{Sub4}
In this subsection, we will investigate the obstruction sheaf to
deforming a $[u]$ in $\mwtdef$ for the case $a_2$, $a_3<0$; we will
follow the convention introduced in Subsection \ref{global section}.

We let $S\to\bar\cQ$ be an $T$-equivariant \'etale neighborhood,
we let
$$
\cU: \cX\lra \cW,\quad\cR\sub\cX,\quad \cD\sub\cW\and \cZ_i\sub\cX
$$
be the tautological family of $S$ of $\mwtdef$. Here $\cW$ is an
$S$-family of $W[\bm]$ of possibly varying $\bm$, $\cD$ is the
relative divisor of $\cW$, $\cR\sub\cX$ is the sections marked
points and $\cZ_i\sub\cX$ is the i-th deformable parts of $\cU$.

Let $\cT^2$ be the obstruction sheaf over ${S}$ of the obstruction
theory of $\mwtdef$. According to \cite{Li2}, its $\soe$-invariant
part, indicated by the subscript $(\cdot)_\soe$, fits into the long
exact sequences
\begin{equation}
\lab{11.4}\lra\ext^1_{\cX/S} (\Omega_{\cX/S}(\cR),\cO_{\cX})_\soe
\mapright{\beta} \cA^1_\soe \mapright{\delta} \tilde\cT^2_\soe \lra
0;
\end{equation}
\begin{equation}\lab{11.5}
\lra \cB^0_\soe\mapright{}
R^1\pi\lsta\bl\cU\sta\Omega_{\cW^\dagger/{S}}(\log \cD)\dual\br_\soe
\mapright{\alpha} \cA^1_\soe\mapright{} \cB^1_\soe\lra 0
\end{equation}
and
\begin{equation}\lab{11.20}
\lra\tilde\cT^1_\soe\lra\cH_\soe\mapright{}\cT^2_\soe\lra
\tilde\cT^2_\soe\lra0.
\end{equation}
Within these sequences, $\cB^i=\oplus_{j=1}^3\cB^i_j$; each summand
$\cB^i_j$ is a sheaf that associates to the smoothing of the nodes
of the fibers of $\cX$ that are mapped under $\cU$ to $D$ ($\sub W$)
and the singular loci of $\Delta[m_i]$; the $\cW^\dagger$ is the
scheme $\cW$ with the log structure defined in \cite{Li2} and
$\Omega_{\cW^\dagger}$ is the sheaf of log differentials. In our
case, $\cU\sta\Omega_{\cW^\dagger/{S}}(\cD)=
\tilde\cU\sta\Omega_W(\log D)$, where $\tilde\cU\mh\cX\to W$ is the
obvious induced morphism.

Without taking the $\soe$-invariant part, the top two exact
sequences define the obstruction sheaf $\tilde\cT^2$ to deforming
$[u]$ in ${\cM_{\chi,\vd,\vmu}^\bu(W\urel)}$ --- the moduli of
relative stable morphisms to $W\urel$ without requiring $u(R)\sub
L$. Taking the invariant part and adding the last exact sequence
defines the obstruction sheaf $\cT^2_\soe$ of $\mwtdef$. The sheaf
$\cH$ is the pull back of the normal line bundle to $\cL\sub \cD$.
For the $\eta$ we are interested, $\cH_\soe=0$; hence the last exact
sequence reduces to $\cT^2_\soe\equiv\tilde\cT^2_\soe$.

In the following we shall show that the $l$ families $\cZ_i\sub\cX$
of deformable parts of $(\cU,\cX)$ each contributes to a weight zero
trivial quotient sheaf of $\cT^2_\soe$.

We begin with the sheaf
\begin{equation}\lab{15.1}
R^1\pi_{k\ast}\bl\Phi_{k}\sta\Omega_W(\log D)\dual\br,
\end{equation}
where $\Phi_{k_i}\mh\fY\to W$ is the family constructed before and
$\pi_{k}\mh\fY\to\Ao$ is the projection. We let ${D}_{12}$ (resp.
${D}_{31}$) be the $T$-invariant divisor of $W$ that contains $C_1$
and $C_2$ (resp. $C_1$ and $C_3$); it is also the proper transform
of the product of the first and second (resp. the first and third)
copies of $\Po$ in $(\Po)^3$. We let $\pi_{12}\mh W\to {D}_{12}$,
$\pi_{31}\mh W\to {D}_{31}$ and $\pi_1\mh W\to C_1$ be the obvious
projections. We claim that $\Omega_W(\log D)|_{D_{12}}$ has a direct
summand $\pi_1\sta\cN_{C_1/{D}_{31}}\dual$. Indeed, via the
projection $\pi_{31}$ we have homomorphism
$$\pi_{31}\sta\cN\dual_{C_1\backslash D_{31}}|_{D_{12}}\lra
\Omega_W|_{D_{12}}\lra\Omega_W(\log D)|_{D_{12}}.
$$
Also via the projection $\pi_{12}$ we have homomorphism
$$\pi_{12}\sta\Omega_{D_{12}}(\log E_{12})|_{D_{12}}\lra
\Omega_W(\log D)|_{D_{12}},\quad E_{12}=D_{12}\cap W.
$$
Combined, we have
$$\pi_{31}\sta\cN\dual_{C_1\backslash D_{31}}|_{D_{12}}\oplus
\pi_{12}\sta\Omega_{D_{12}}(\log E_{12})|_{D_{12}}\lra \Omega_W(\log
D)|_{D_{12}},
$$
which can easily be shown an isomorphism. This proves that
$\Omega_W(\log D)|_{D_{12}}$ has a direct summand
$\pi_1\sta\cN_{C_1/{D}_{31}}\dual$. Consequently,
$\Phi_{k}\sta\Omega_W(\log D)\dual$ has a direct summand
$\Phi_{k}\sta(\pi_1\sta\cN_{C_1/ {D}_{31}})$.

Because of our choice, the weight of $dz_i$ is $a_i$; the weight of
$T_0\dual \fY_c$ at $0$ is $1/k_i$ and the weight of
$\Phi_{k}\sta(\pi\sta_1\cN_{C_1/{D}_{31}})$ at $0\times\Ao\sub\fY$
is $-a_3$. Hence, the sheaf \eqref{15.1} splits to line bundles of
weights
$$-a_3-a_1+\frac{1}{k},\,-a_3-a_1+\frac{2}{k},\cdots\cdots,
-a_3-\frac{1}{k}.
$$
Since all $a_i$ are integers, and $a_3\leq -1$ and $-a_3-a_1=a_2\leq
-1$, within the above list there is exactly one that is zero. Hence
\begin{equation}\lab{12.1}
R^1\pi_{k\ast}\bl\Phi_{k}\sta\Omega_W(\log D)\dual\br_\soe\cong
\cO_{\Ao}.
\end{equation}

We now let $\rho_i\mh{S}\to \Ao_i$ be so that
$\cU|_{\cZ_i}\cong\rho_i\sta\Phi_{k_i}$.
 Since $\cZ_i\sub\cX$ is a flat
family of subcurves,
$$
R^1\pi\lsta\bl\cU\sta\Omega_{\cW^\dagger/{S}}(\log
\cD)\dual\br_\soe
\lra R^1\pi\lsta\bl\cU\sta\Omega_{\cW^\dagger/{S}}(\log
\cD)\dual|_{\cZ_i}\br_\soe
$$
is surjective; but the last term is isomorphic to the pull back
$\rho_i\sta$ of (\ref{12.1}); hence we obtain a quotient sheaf
\begin{equation}\lab{11.15}
\varphi_i: R^1\pi\lsta\bl\cU\sta\Omega_{\cW^\dagger/S}(\log
\cD)\dual\br_\soe \lra\rho_i\sta \cO_{\Ao_i}.
\end{equation}

\begin{lemm} \lab{lift-1}
The homomorphism $\varphi_i$ canonically lifts to surjective
\begin{equation}\lab{11.6}\hat\varphi_i:
\cT^2_\soe\lra \rho_i\sta \cO_{\Ao_i}.
\end{equation}
\end{lemm}

The default proof is to follow the construction of the sheaves and
the exact sequences in (\ref{11.4}-\ref{11.20}); once it is done,
the required vanishing will follow immediately. However, to follow
this strategy, we need to set up the notation as in \cite{Li2} that
itself requires a lot of efforts. Instead, we will utilize the
decomposition of $S$ to give a more conceptual argument; bypassing
some straightforward but tedious checking.

We first decompose $\cU$ into four subfamilies. Since $\cW/S$ is a
family of expanded pairs of $(W,D)$, $\cW^{[0]}=W\times S$ is a
closed subscheme of $\cW$. We then let
$\cX^{[0]}=\cU\upmo(\cW^{[0]})$. Because of Lemma \ref{intact},
$$
\cU^{[0]}=\cU|_{\cX^{[0]}}:\cX^{[0]}\lra \cW^{[0]}
$$
is an $S$-family of relative stable morphisms relative to
$\cD^{[0]}=D\times S\sub \cW^{[0]}$. Next we consider the composite
$$\tilde \cU: \cX\lra\cW\lra \cW^{[0]}
$$
and the preimage $\tilde\cU\upmo(D_i)$. Because of the same reason,
either this preimage is a flat family of nodes over $S$ or is a flat
family of curves over $S$. In the former case we agree
$\cX^{[i]}=\emptyset$ and in the later case we define
$$\cX^{[i]}=\tilde\cU\upmo(D_i\times S),\quad
\cU^{[i]}=\cU|_{\cX^{[i]}}: \cX^{[i]}\lra \cW^{[i]},
$$
where the last term $\cW^{[i]}$ is the $S$-family of $\Delta[m_i]$'s
that were attached to $\cW^{[0]}$ along $D_i\times S$ to form $\cW$.

Since $\cU^{[0]}\mh \cX^{[0]}\to \cW^{[0]}$ is a family of
$\soe$-equivariant relative stable morphisms, and since $\cU^{[i]}$
is a family of $\soe$-equivariant relative stable morphisms to
$\Delta$ relative to $D_-$ and $D_+$, modulo an additional
equivalence induced by the $\CC\sta$ action on $\Delta$, the
obstruction sheaves $\cT^{[i],2}$ over $S$ to deforming $\cU^{[i]}$
as $\soe$-equivariant maps fit into similar exact sequences
\begin{equation}
\lab{11.8}\lra\ext^1_{\cX^{[i]}/S} (\Omega_{\cX^{[i]}/\hat
S}(\cR^{[i]}),\cO_{\cX^{[i]}})_\soe \mapright{\beta^{[i]}}
\cA^{[i],1}_\soe \mapright{\delta^{[i]}} \cT^{[i],2}_\soe \lra 0
\end{equation}
and
\begin{equation}\lab{11.9}
\lra \cB^{[i],0}_\soe\mapright{}
R^1\pi\lsta\bl\cU^{[i]\ast}\Omega_{\cW^{[i]\dagger}/S}(\log \cD)
\dual\br_\soe \mapright{\alpha^{[i]}} \cA^{[i],1}_\soe\mapright{}
\cB^{[i],1}_\soe\lra 0.
\end{equation}
Here we have already used the observation that
$\tilde\cT^{[i],2}_\soe=\cT^{[i],2}_\soe$.

Now let $\cN_{\text{sp}}\sub\cU\upmo(D_i\times S)$ be any section of
nodes of $\cX$ that separates $\cX^{[0]}$ and $\cX^{[i]}$. By Lemma
\ref{lem:auto}, the induced $\soe$-automorphisms on the connected
component of $\cX^{[0]}$ adjacent to $\cN_{\text{sp}}$ is
infinite and on $\cX^{[i]}$ is finite. Therefore the
$\soe$-invariant parts
$$
\ext^1_{\cX/S}
(\Omega_{\cX/S}(\cR),\cO_{\cX})_\soe=\Oplus_{i=0}^3\ext^1_{\cX^{[i]}/S}
(\Omega_{\cX^{[i]}/ S}(\cR^{[i]}),\cO_{\cX^{[i]}})_\soe.
$$
For the similar reason, because the tangent bundle $T_{p_i}W$ has no
weight $0$ non-trivial $\soe$-invariant subspaces,
\begin{equation}\label{11.10}
R^1\pi\lsta\bl\cU^{\ast}\Omega_{\cW^\dagger/S}(\log
\cD)\dual\br_\soe= \Oplus_{i=0}^3
R^1\pi\lsta\bl\cU^{[i]\ast}\Omega_{\cW^{[i]\dagger}/S}(\log
\cD^{[i]})\dual\br_\soe,
\end{equation}
where $\cD^{[i]}$ is the relative divisor of $\cW^{[i]}$. Further,
if we follow the definition of the sheaves $\cB^i$ and $\cA^i$, we
can prove that
\begin{equation}\lab{11.21}
\Oplus_{i=0}^3\cA^{[i],j}_\soe=\cA^j_\soe\and
\Oplus_{i=0}^3\cB^{[i],j}_\soe=\cB^j_\soe;
\end{equation}
that under these isomorphisms,
\begin{equation}\lab{11.22}
\Oplus_{i=0}^3\alpha^{[i]}=\alpha,\quad
\Oplus_{i=0}^3\beta^{[i]}=\beta\and
\Oplus_{i=0}^3\delta^{[i]}=\delta;
\end{equation}
and
\begin{equation}\lab{11.23}
\Oplus_{i=0}^3\cT^{[i],2}_\soe=\cT^2_\soe.
\end{equation}
Finally, the exact sequences (\ref{11.4}) and (\ref{11.5}) become
the direct sums of the exact sequences (\ref{11.8}) and
(\ref{11.9}).

Now we come back to the weight zero quotient $\varphi_i$ in
(\ref{11.15}). By its construction, $\varphi_i$ is merely the
canonical quotient homomorphism
\begin{equation}\lab{11.24}
R^1\pi\lsta\bl\cU^{[0]\ast} \Omega_{\cW^{[0]}/S}(\log
\cD^{[0]})\dual\br_\soe \lra\qquad\qquad
\end{equation}
$$
\qquad\qquad\qquad\qquad \lra R^1\pi\lsta\bl\cU^{[0]\ast}
\Omega_{\cW^{[0]}/S}(\log\cD^{[0]})\dual|_{\cZ_i}\br_\soe
=\rho\sta_i\cO_{\Ao_i}.
$$
under the isomorphism (\ref{11.10}). Because of (\ref{11.23}), to
lift $\varphi_i$ to $\hat\varphi_i$ we only need to lift
(\ref{11.24}) to $\cT^{[0],2}_\soe\to\rho\sta_i\cO_{\Ao_i}$.

For this, we need to look at the exact sequence (\ref{11.9}) for
$\cX^{[0]}$. Since $\cU^{[0]}$ is a relative stable map to $(W,D)$
--- namely no $\Delta$'s has been attached to $W$ --- the sheaf
$\cB^{[0],j}=0$. Therefore the sequence (\ref{11.9}) reduces to
$\alpha^{[0]}=\text{id}$. On the other hand, $\cT^{[0],2}_\soe$ is
the obstruction sheaf on $S$ to deformations of $\cU^{[0]}$. Since
$\cZ_i$ is a family of connected components of $\cX^{[0]}/S$, the
exact sequence (\ref{11.8}) decomposes into direct sum of individual
exact sequences that contains
\begin{equation}\lab{11.25}
\lra \ext^1_{\cZ_i/S}\bl\Omega_{\cZ_i/\hat
S}(\cR^{[0]}),\cO_{\cZ_i}\br_\soe \mapright{\beta^{[\hat \cZ_i]}}
\qquad\qquad\qquad
\end{equation}
$$\qquad\qquad\qquad \lra R^1\pi\lsta\bl\cU^{[0]\ast} \Omega_{\cW^\dagger/\hat
S}(\log\cD)\dual|_{\cZ_i}\br_\soe
\mapright{\delta^{[\cZ_i]}}\cT^{[\cZ_i],2}_\soe\lra0,
$$
as its factors.

For $\cZ_i$, since it is smooth, it has expected dimension zero and
has actual dimension one, the obstruction sheaf
$\cT^{[\cZ_i],2}_\soe$ must be a rank one locally free sheaf on $S$.
Then because the middle term in (\ref{11.25}) is
$\rho_i\sta\cO_{\Ao_i}$, which is a rank one locally free sheaf, the
arrow $\delta^{[\cZ_i]}$ must be an isomorphism while
$\beta^{[\cZ_i]}=0$. Hence $\varphi_i$ lifts to
$$\cT^{[0],2}_\soe\equiv
\Oplus_j\cT^{[\cZ_j],2}_\soe\lra
\cT^{[\cZ_i],2}_\soe\equiv\rho_i\sta\cO_{\Ao_i},
$$
and lifts to $\hat\varphi_i\mh\cT^2_\soe\to \rho_i\sta\cO_{\Ao_i}$,
thanks to (\ref{11.23}).

\subsection{The Case for $\eta=(1,-1,0)$}\lab{later}

We now investigate the structures of maps $[u]\in\mwtdef$ in case
$\eta=(1,-1,0)$. Let $(u,X)$ be any such map, let $R$ be the marked
points and let $\tilde u$ be the contraction $X\to W$. Because
$a_3=0$, $\tilde u(X)$ intersects $D_1$ at $p_1$; intersects $D_3$
at $p_3$ while intersects $D_2$ can be any point in $L_2$. Thus
being $\soe$-equivariant forces $\tilde u(X)$ to be a finite union
of a subset of $C_1$, $C_2$, $C_3$ and the lifts of the sets
$\{z_1z_2=c,z_3=0\}\sub(\Po)^3$.

In case all irreducible components are mapped to $\cup\, C_i$ under
$\tilde u$, $[u]\in\Mdyy$. For those that are not in $\Mdyy$, there
bound to be some $Y\sub X$ so that $\tilde u(Y)$ is the lifts of
$\{z_1z_2=c,z_3=0\}$. Such $u|_Y$ are realized by the morphism
$\phi_{k,c}\mh\Po\to W$ that are the lifts of
\begin{equation}
\xi\mapsto(c^k\xi^k,\xi^{-k},0)\in(\Po)^3.
\end{equation}

When $c$ specializes to $0$, the map $\phi_{k,c}$ specializes to
$\phi_{k,0}\mh \Po\sqcup\Po\to W$ that is the lift of
$\xi_1\mapsto(\xi_1^k,0,0)$ and $\xi_2\mapsto(0,\xi_2^{-k},0)$.
Indeed, there is a family $\fY\to\Ao$ and a morphism $\Phi_k\mh
\fY\to W$ so that its fiber over $c\in\Ao$ is the $\phi_{k,c}$
defined; also this is a complete list $\soe$-equivariant
deformations of $\phi_{k,c}$. Since the argument is exactly the same
as in the case studied, we shall not repeat it here.

Here comes the main difference between this and the case studied
earlier. In the previous case, $\image\phi_{k,c}\cap D_i=p_i$ for
both $i=1$ and $2$; hence we can deform each $u|_Y\cong\phi_{k,c}$
to produce an $\Ao$ family in $\mwtdef$. In the case under
consideration, though $\image\phi_{k,c}\cap D_1=p_1$, if we fix an
embedding $\Ao\sub L_2$ so that $0\in\Ao$ is the $p_2\in L_2$, then
$\image\phi_{k,c}\cap D_2=c^k\in L_2$. In other words, if we deform
$u|_F\cong \phi_{k,c}$, we need to move the connected component of
$X^{[2]}$ that is connected to $Y$.

This leads to the following definition.

\begin{defi} We say that a connected component
$Y\sub X^{[0]}$ is subordinated to a connected component $E\sub
X^{[2]}$ if $Y\cap E\ne\emptyset$; we say a connected component
$E\sub X^{[2]}$ is deformable if every connected component of
$X^{[0]}$ that is subordinate to $E$ is of the form $\phi_{k,c}$ for
some pair $(k,c)$. We say $u$ has deformation type $(k_i)_l=(k_1\leq
\cdots\leq k_l)$ if it has exactly $l$ deformable connected
components $\phi_{k_1,c_1},\cdots,\phi_{k_l,c_l}$ in $X^{[2]}$.
\end{defi}

 The deformation types
define a function on $\mwtdef$.

\begin{lemm} \lab{A14}The deformation type function is locally constant on
$\mwtdef$.
\end{lemm}

\begin{proof}
The proof is parallel to the case studied, and will be omitted.
\end{proof}

As in the previous case, any $[u]\in\mwtdef$ of deformation type
$(k_i)_l$ generates an $\A^l$ in $\mwtdef$ so that its origin lies
in $\Mdyy$. Let $E_1,\cdots,E_l\sub X^{[2]}$ be the complete set of
deformable parts of $u$; let $Y_{i,j}$, $j=1,\cdots,n_i$ be the
complete set of connected components in $X^{[0]}$ that are
subordinate to $E_i$. By definition, each $u|_{Y_{i,j}}\cong
\phi_{k_{i,j},c_{i,j}}$. To deform $u$, we shall vary the $c_{i,j}$
in each $\phi_{k_{i,j},c_{i,j}}$ and move $E_i$ accordingly to get a
new map.

In accordance, we shall divide $X$ into three parts. We let $X_0$ be
the union of irreducible components of $X$ other than the $E_i$'s
and $Y_{i,j}$'s. The variation of $u$ will remain unchanged over
this part of the curve. The second part is the moving part $E_i$'s.
Recall that each $u|_{E_i}$ is a morphism to $\Delta[m_2]$. Suppose
it maps to the fiber $\Delta[m_2]_c$ of $\Delta[m_2]$ over $c\in
L_2\sub D_2$. To deform $u$, we need to make the new map maps $E_i$
to $\Delta[m_2]_{c\pri}$. Since the total space of $\Delta[m_2]$
over $L_2$ is a trivial $\Po[m_2]$ bundle, there is a canonical way
to do this. We let
$$
\varphi_{c,c\pri}: \Delta[m_2]_c\lmapright{\cong}\Delta[m_2]_{c\pri}
$$
be the isomorphism of the two fibers of $\Delta[m_2]$ over $c$ and
$c\pri\in L_2$ induced by the projection $\Delta[m_2]\to\Po[m_2]$
that is induced by the product structure on $\Delta[m_2]$ over
$L_2$. The third parts are those $Y_{i,j}$ that are subordinate to
$E_i$.

We now deform the map $u$ using the parameter space $\A^l$. We let
$K_i$ be the least common multiple of $(k_{i,1},\cdots,k_{i,n_i})$;
we let $e_{i,j}=K_i/k_{i,j}$. Since $Y_{i,j}$ and $Y_{i,j\pri}$ are
connected to the same connected component $E_i\sub X^{[2]}$,
$c_{i,j}^{k_{i,j}}=c_{i,j\pri}^{k_{i,j\pri}}$; we let it be $c_i$.
For $\bt=(t_1,\cdots,t_l)\in\A^l$, we define
\begin{equation}\lab{14.3}
u^\bt|_{X_0}=u|_{X_0},\quad u^\bt|_{E_i}=\varphi_{c_i, t_i^{K_i}}
\circ u|_{E_i}
\and u^\bt|_{Y_{i,j}}=\phi_{k_{i,j},t_i^{e_{i,j}}}.
\end{equation}
Here in case $Y_{i,j}\cong \Po$, which is the case when $c_{i,j}\ne
0$, by $u^\bt|_{Y_{i,j}}=\phi_{k_{i,j},0}$ we mean that we will
replace $Y_{i,j}$ by $\Po\sqcup\Po$ with necessarily gluing if
required; and vice versa.

The $\A^l$ family $u^\bt$ is a family of $\soe$-equivariant relative
stable morphisms in $\mwtdef$; the map $u^0$ associated to
$0\in\A^l$ lies in $\Mdyy$; the induced morphism $\A^l\to\mwtdef$ is
an embedding up to a finite quotient.

By extending this to any connected component $\cQ$ of $\mwtdef$, we
obtain

\begin{prop}\lab{A15}
Let $\cQ$ be any connected component of $\mwtdef$ that is not
entirely contained in $\Mdyy$. Suppose elements of $\cQ$ has
deformation type $(k_i)_l$. Then there is a stack $\bar\cQ$, a
finite quotient morphism $\bar\cQ/G_{(k_i)}\to \cQ$, a closed
substack $\bar\cQ_0\sub\bar\cQ$, $l$ projections
$\rho_i\mh\bar\cQ\to\Ao_i$ and a projection $\pi\mh\bar\cQ\to
\bar\cQ_0$ so that
$$
\bl\pi,(\rho_1,\cdots,\rho_l)\br:
\bar\cQ\lmapright{\cong}\bar\cQ_0\times\A^l
$$
is an isomorphism. Further, given a $[u]\in\cQ$, the fiber $\A^l$ in
$\bar\cQ$ that contains a lift of $[u]\in\cQ$ is the $\A^l$ family
$\{u^\bt\mid \bt\in\A^l\}$; its intersection with the zero section
$\bar\cQ_0$ is $u^0$. Finally, the intersection $\cQ\cap\Mdyy$ is
the image of $\bar\cQ_0$.
\end{prop}

\begin{proof}
Let $\cU\mh\cX\to\cW$ be the tautological family over $\bar\cQ$. We
choose $\bar\cQ$ so that there are families of subcurves
$\cE_1,\cdots,\cE_l\sub\cX$ so that for each $z\in\bar\cQ$,
$\cE_1\cap\cX_z,\cdots,\cE_l\cap\cX_z$ are exactly the $l$
deformable parts of $\cX_z$. Then the composite $\cE_i\to\cW\lra W$
factor through $L_2\sub W$, and the resulting morphism $\cE_i\to
L_2$ factor through $\bar\cQ\to L_2$. Because each $\cE_i\cap\cX_z$
has a $\phi_{k,c}$ connected to it, the image of $\bar\cQ\to L_2$
lies in $L_2-q_2$. We then fix an isomorphism $\Ao\cong L_2-q_2$
with $0$ corresponding to $p_2$. This way we obtain the desired
morphism
$$
\rho_i: \bar\cQ\lra \Ao_i\cong L_2-q_2.
$$
The proof of the remainder part of the Proposition is exactly the
same as the case studied; we shall not repeat it here.
\end{proof}

The last step is to investigate the obstruction sheaf over $\cQ$, or
its lift to $\bar\cQ$.

Let $\cR\sub\cX$ the divisor of marked points. By passing to an
\'etale covering of $\bar\cQ$, we can assume that $\cR\to\cQ$ is a
union of sections; in other words, we can index the marked points of
$[u]$ in $\bar\cQ$ globally. We then pick an indexing so that for
$i\leq l$ the $i$-th section of the marked points $\cR_i$ lies in
$\cE_i$. For $i=1,\ldots, n$, where $n$ is the number of marked points,
we let $\cU_i\mh\bar\cQ\to \cL_2$ be
$$
\cU_i\stackrel{\mathrm{def}}{=}\cU|_{\cR_i}: \cR_i\cong\bar\cQ\lra\cL_2\sub\cW.
$$
Since $\cL_2\sub\cD_2$ is isomorphic to $L_2\times\bar\cQ\sub
D_2\times\bar\cQ$ under the contraction $\cW\to W\times\bar\cQ$ and
since $\cR_i$ lies in $\cE_i$, for $i\leq l$ the morphism $\cU_i$ is
exactly the $\rho_i$ under the isomorphism $\Ao_i\cong L_2-q_2$, and
$\cU_i\sta\cN_{\cL/\cD}$ is canonically isomorphic to $\rho_i\sta
N_{L_2/D_2}$. Because $D_2$ is fixed by $\soe$, $N_{L_2/D_2}$ is
fixed as well, and hence $\rho_i\sta N_{L_2/D_2}$ is a trivial line
bundle on $\bar\cQ$ with trivial $\soe$-linearization.

Because $\cH=\oplus_{i=1}^n \cU_i\sta\cN_{\cL/\cD}$ (see Section
\ref{Sub4}), $\oplus_{i=1}^l\rho\sta_i N_{L_2/D_2}$ becomes a direct summand of
$\cH$. Because it has weight zero, it induces a canonical
homomorphism
$$
\oplus_{i=1}^l\rho\sta_i N_{L_2/D_2}\lra\cT^2_\soe,
$$
a weight zero subsheaf of $\cT^2_\soe$.

\begin{lemm}\lab{A16}
The homomorphism $\oplus_{i=1}^l \rho_i\sta N_{L_2/D_2}\to\cT^2_\soe$ in
(\ref{11.20}) is injective; thus $\cT^2_\soe$ contains
$\oplus_{i=1}^l \rho_i\sta N_{L_2/D_2}$ as its subsheaf. Indeed, this
subsheaf is canonically a direct summand of $\cT^2_\soe$.
\end{lemm}

\begin{proof}
First the first $l$ marked points lie in the connected components of
$X^{[2]}$ that are connected to the domain of at least one
$\phi_{k,c}$ in $W$. Because all deformations of $\phi_{k,c}$ as
$\soe$-invariant maps are $\phi_{k,c\pri}$, and they intersect $D_2$
in $L_2$ only; hence for these $i$ even if we do not impose the
condition $\cU(\cR_i)\sub\cL_2$ the condition will be satisfied
automatically. In short, the arrow $\tilde\cT^1_\soe\to\cH_\soe$ has
image lies in the summand $\oplus_{i=l+1}^n\cU\sta_i\cN_{\cL/\cD}$. This
proves that the homomorphism $\oplus_{i=1}^l\rho_i\sta
N_{L_2/D_2}\to\cT^2_\soe$ is injective.

We now show that this subsheaf is canonically a summand of the
obstruction sheaf. The ordinary moduli of stable relative morphisms
$\cM_{\chi,\vd,\vmu}^{\bullet}(W\urel)$ requires that the marked
points be sent to the relative divisor. The moduli space
$\Mdw=\cM_{\chi,\vd,\vmu}^{\bullet}(W\urel,L)$ we worked on imposes
one more restriction: the marked points be sent to $L\sub D$. The
obstruction sheaves of the two moduli spaces are related by the
exact sequence (\ref{11.20}) because of the exact sequences
$$
0\lra N_{L_i/D_i}\lra N_{L_i/W}\lra N_{D_i/W}|_{L_i}\lra 0.
$$
In our case, $L_i$ is a $\Po$ and the above exact sequence splits
$T$-equivariantly. Hence the sheaf $\cT^2_\soe$ splits off a factor
that is the kernel of $\tilde\cT^1_\soe\to\cH_\soe$. Therefore
$\oplus^l_{i=1}\rho_i\sta N_{L_2/D_2}$, which is a summand of
$\cH_\soe$ and a subsheaf of $\cT^2_\soe$, becomes a summand
$\cT^2_\soe$.
\end{proof}

\subsection{A Criterion for Constancy}
We first remark that due to the nature of the discussion in the
remainder of this section, we shall use analytic topology of instead
of Zariski topology. In particular, unless otherwise specified, all
open sets will be analytic open sets. Accordingly, we shall view
DM-stack as orbifold with analytic topology.

Before we proceed to the details of the proof, a quick review of the
construction of the virtual cycles of the moduli stack is in order.
Let $T=(\CC\sta)^2$ as before and let $M$ be a proper
Deligne-Mumford stack acted on by $T$ and endowed with a
$T$-equivariant perfect obstruction theory. As shown in
\cite{Beh,Beh-Fan,Li-Tia}, the virtual cycle $[M]\virt$ is
constructed by
\begin{enumerate}
\item[{\sl a}.] identifying the perfect obstruction theory of $M$;
\item[{\sl b}.] picking a vector bundle\footnote{It was shown in
\cite{Li2} the existence of a global vector bundle $E$ can be
replaced by that $M$ is dominated by a quasi-projective scheme.}
(locally free sheaf) $E$ on $M$ so that it surjects onto the
obstruction sheaf of $M$ and
\item[{\sl c}.] constructing an associated cone $C\sub E$ of pure
dimension $\rank E$.
\end{enumerate}
The virtual cycle $[M]\virt$ is the image of the cycle $[C]\in
H_*(E,E-M)$ under the Thom isomorphism
$$\varphi_E: H_*(E,E-M)\lra H_{*-2r}(M),
\quad r=\rank E.
$$
Here as usual, we denote by $E$ the total space of $E$ and denote by
$M \sub E$ its zero section that is isomorphic to $M$. Also, all
(co)homologies are taken with $\QQ$ coefficient.

Following \cite{Gra-Pan}, we can make the above construction
$T$-equivariant. We choose $E$ be a $T$-equivariant vector bundle.
Then the cone $C$ alluded before is a $T$-invariant subcone of
${E}$. Because $C\sub {E}$ is $T$-equivariant, the composite
$$T\times C\mapright{\text{pr}_2}C\lra {E}
$$
defines a $T$-equivariant $[C]^\te\in H^\te\lsta({E},{E}-M)$; its
image under the $T$-equivariant Thom isomorphism $\varphi_{E}$ is
the $T$-equivariant virtual moduli cycle
$$
\varphi_{E}\bl[C]^\te\br=[M]\virtt\in H^T\lsta(M).
$$
Note that the equivariant homologies are $H\sta_T(pt)=\QQ[u_1,u_2]$
modules. As before, we denote by $\fm=(u_1,u_2)$ the maximal ideal
and $H\lsta^T(\cdot)_\fm$ the localization at $\fm$.

Next, we apply the localization theorem to the class $[M]\virtt$.
Let
$$
M^T=\coprod\bigl\{M_a\mid a\in A\bigr\}
$$
be the decomposition of the $T$-fixed loci into connected
components; let
$$
\tau_a: H\lsta^T({M} _a)_\fm \lra H\lsta^T(M)_\fm
$$
be induced by the inclusion. According to \cite{Gra-Pan}, to each
$M_a$ there is a canonically defined virtual cycle $[M_a]\virtt\in
H^\te\lsta(M_a)$ and a virtual $T$-equivariant normal bundle
$\cN_a\virt=\cT^{1,m}-\cT^{2,m}$ so that, after localization,
$$
[M]\virtt_\fm=\sum_{a\in A}\tau\las\bbl
\frac{[M_a]\virtt}{e^\te(\cN_a\virt)}\bbr \in \bl
H^\te\lsta(M)_\fm\br_0.
$$
Here $(\,\cdot\,)_0$ is the degree zero part of the graded ring
inside the parentheses.

Now suppose we have a $T$-invariant close substack $N \sub M$; we
assume that for any $a\in A$, either $M_a\cap N=\emptyset$ or
$M_a\sub N$. Such $N$ divides $A$ into a part of those that are
contained in $N$ and a part that are not.
Summing over those inside $N$, we obtain
$$
[M]\virtt_{\fm,N}=\sum_{M_a\sub
N}\tau\las\bbl\frac{[M_a]\virtt}{e^\te(\cN_a\virt)}\bbr \in \bl
H^\te\lsta(M)\otimes \QQ[u_1,u_2]_\fm\br_0.
$$
Now suppose $M$ has virtual dimension zero, then taking degree
(defined in Section \ref{sec:formalGW})
$$\deg_\fm: \bl
H^\te\lsta(M)_\fm\br_0\lra \QQ(u_1/u_2),
$$
we obtain a rational function in $u_1/u_2$:
$$\deg_\fm [M]\virtt_{\fm,N} \in \QQ(u_2/u_2).
$$

To find a sufficient condition that makes the above quantity
independent of $u_1/u_2$, we shall devise a way to show that this
partial sum is the localization of a purely topological quantity. To
achieve this, we let $T_\RR\sub T$ be the compact maximal real
subgroup and we pick a $T_\RR$-equivariant metric on $M$;
we can do this by viewing $M$ as a possibly singular orbifold.
We then form a ${T_\RR}$-invariant closed tubular neighborhood of
${N}$:
$$
\Sigma=\left\{u\in M\mid \text{dist}(u,{N})\leq \ep\right\}\sub M .
$$
We then close the boundary
$$\pa\Sigma=\left\{u\in M\mid \text{dist}(u,{N})= \ep\right\}\sub M .
$$
of $\Sigma$ by picking a one-parameter subgroup ${T}_\nu\sub T$ and
contract individual $S^1_\nu\sub {T}_\nu$ orbits in $\pa\Sigma$.
Namely, for $S^1_\nu$ the maximal compact subgroup of ${T}_\nu$ and
for $z\in \pa\Sigma$ we shall contract the whole orbit $S^1_\nu\cdot
z\sub\pa \Sigma$ to a single point $[S^1_\nu\cdot z]$. To do so, we
will pick the subgroup ${T}_\nu\sub {T}$ so that
${M}^{{T}_\nu}={M}^T$. We denote the resulting space by $\tilde M$.
Note that because $M$ is algebraic and proper, $\tilde M$ is a
compact singular orbifold.
For later study, we denote $\partial\Sigma/S^1_\nu$ by
$\tilde{M}_\infty$.

We next define the virtual cycle of $\tilcM$. Because ${E}$ is a
${T_\RR}$-equivariant vector bundle, and because points in
$\pa\Sigma$ have finite $S^1_\nu$-stabilizers, ${E}$ descends to a
vector bundle $\tilde {E}$ on $\tilde M$. For the same reason, the
cone $C$ descends to a cone $\tilde C\sub \tilde {E}$. Because $C$
is algebraic and because $\ep$ is sufficiently small, $\tilde C$
defines an element $[\tilde C]\in H\lsta^\te(\tilde {E},\tilde
{E}-\tilde M)$; since all data are ${T_\RR}$-equivariant, it also
defines an equivariant class
$$[\tilde C]^\te\in H\lsta^\te(\tilde {E},\tilde {E}-\tilde M).
$$
We define their respective images under the obvious Thom
isomorphisms their virtual and equivariant virtual cycles:
$$[\tilde M]\virt\in H\lsta(\tilde{M})\and
[\tilde M]\virtt\in H\lsta^\te(\tilde{M}).
$$

We now localize the class $[\tilde M]\virtt$. We let
$$
\tilde M^{T_\RR}=\coprod\bigl\{\tilde M_a\mid a\in \tilde A\bigr\}
$$
be the decomposition of the ${T_\RR}$-fixed loci into connected
components. Because either $M_a\sub N$ or $M_a\cap N=\emptyset$,
$\tilde M^{T_\RR}$ is the disjoint union of $\{M_a\mid M_a\sub N\}$
and $(\tilde M_\infty)^{T_\RR}$.
We now look at the arrows between localized modules:
$$H\lsta^\tz(\tilde E,\tilde E-\tilde M)_\fm \lmapright{\text{Thom}}
H^\tz\lsta(\tilde M)_\fm\lmapright{\iota} \Oplus_{a\in\tilde
A}H^\tz\lsta(\tilde M,\tilde M-\tilde M_a)_\fm.
$$
Here $\iota=\oplus\iota_a$ and
$$\iota_a:H^\tz\lsta(\tilde M)_\fm\lmapright{\iota}
H^\tz\lsta(\tilde M,\tilde M-\tilde M_a)_\fm.
$$
By localization theorem, $\iota$ is an isomorphism.

Since $M$ has virtual dimension zero, $\tilde M$ has virtual
dimension zero as well. Applying the degree map $\deg_\fm$ from
$H^\tz\lsta(\tilde M)_\fm$ and $H^\tz\lsta(\tilde M,\tilde M-\tilde
M_a)_\fm$ to $\QQ(u_1/u_2)$, we obtain
$$\deg_\fm[\tilde M]\virtt=\sum_{a\in\tilde A} \deg_\fm
\iota_a\bl[\tilde M]\virtt\br\in\QQ(u_1/u_2).
$$
Since $\tilde M$ is compact, the left hand side is independent of
$u_1/u_2$. This immediately proves a useful criterion

\begin{lemm}\lab{lem:const}
Let the notation be as before and suppose $\iota_a\bl[\tilde
M]\virtt\br=0$ for all $\tilde M_a\sub(\tilde M_\infty)^{T_\RR}$.
Then the sum
$$\sum_{M_a\sub N}\deg_\fm
\iota_a\bl[\tilde M]\virtt\br\in\QQ(u_1/u_2)
$$
is a constant.
\end{lemm}

\subsection{Proof of Theorem \ref{thm:R2}}

We shall verify this criterion for $M=\Mdw$ and $N=\Mdyy$; this will
complete the proof of Theorem \ref{thm:R2}.

As instructed by the criterion, our first step is to classify the
connected components of $(\tilde M_\infty)^{T_\RR}$. Let $\tilde
z\in\tilde M_\infty$ be any closed point fixed by ${T_\RR}$. By the
construction of $\tilde M_\infty$, $\tilde z$ is an
${S^1_\nu}$-orbit $[{S^1_\nu}\cdot z]$ of some $z\in\pa\Sigma$;
therefore $\tilde z$ is fixed by ${T_\RR}$ if and only if the
${S^1_\nu}$-orbit ${S^1_\nu}\cdot z$ is identical to the
${T_\RR}$-orbit ${T_\RR}\cdot z$, which is possible only if
$\dim_\RR\stab_{T_\RR}(z)\geq 1$. Because $\tilde z\in\pa\Sigma$, it
is not in $M^{T}$; hence there must be a subgroup $\soe\sub {T}$ so
that $\tilde z\in M^\soe$. Finally, because $\ep$ is sufficiently
small, $z\in\mwtdef$. This shows that
$$
\coprod\bigl\{\tilde M_a\mid \tilde M_a\sub(\tilde
M_\infty)^T\bigr\}=\coprod_{\soe\sub
{T}}\bbl\mwtdef\cap\pa\Sigma\bbr/{S^1_\nu}.
$$


We now analyze in more details the individual connected components
appear in the right hand side of the above decomposition. Before we
move on, we remark that we only need to consider the case
$\eta=(a_1,a_2,a_3)$, for $a_1>0$ and $a_2, a_3<0$, and the case
$\eta=(1,-1,0)$. Indeed, since the symmetry of $(\Po)^3$ defined by
$(z_1,z_2,z_3)\mapsto(z_2,z_3,z_1)$ lifts to a symmetry of $W$, any
statement that holds true for the $\eta=(a_1,a_2,a_3)$ holds true
for $\eta\pri=(a_2,a_3,a_1)$. Consequently, we only need to work
with those $\eta$ so that $|a_1|\geq |a_2|$ and $|a_3|$. Then
because ${T}_\eta={T}_{-\eta}$, we can assume further that $a_1>0$.
Hence either $a_2$ and $a_3<0$ or one of them is zero. For former is
the {\sl case one}; in the later case, by applying the $S_3$
symmetry we can reduce it to the case $\eta=(1,-1,0)$.

We fix a $\soe\sub {T}$ belongs to the two classes just mentioned.
We let $\cQ_a$ be a connected component of $\mwtdef$ associated to
$\tilde M_a$. According to Proposition \ref{A10} and \ref{A15},
after a finite branched covering $\pi_a:\bar\cQ_a\to\cQ_a$,
$\bar\cQ_a$ is isomorphic to $\cQ_{a,0}\times\A^l$ for some integer
$l>0$; the ${S^1_\nu}$-action on $\bar\cQ_a$ is the product of the
action on $\bar\cQ_{a,0}$ induced by that on $\Mdyy$ and the action
\begin{equation}\lab{15.3}
(u_1,\cdots,u_l)^\sigma=(\sigma^{w_1}u_1,\cdots,\sigma^{w_l}u_l)
\in\A^l
\end{equation}
for some $\bw=(w_1,\cdots,w_l)$ where $w_i$ are nonzero
rational numbers. In case some $w_i$ are non-integers, we let $d$ be the least common multiple of 
the denominators of all $w_i$ and replace the $S^1$ action by the composing it with the 
degree $d$ homomorphism $S^1\to S^1$. This way the new exponents are  $d w_i$, which are integrals. 
Thus without loss of generality, we can assume that all $w_i$ are integers in the first place. 
Hence if we let
$P_a\mh\bar\cQ_a\to\A^l$ be the projection, which is
$(\rho_1,\cdots,\rho_l)$ by our convention, and if we endow $\A^l$
with the ${S^1_\nu}$-action (\ref{15.3}), then $\bar\cQ_a\to\A^l$ is
${S^1_\nu}$-equivariant.

We now pick an ${S^1_\nu}$-invariant Riemannian metric on $\A^l$; we
let $S^{2l-1}_\ep\sub\A^l$ be the $\ep$-sphere under this metric.
$P_a\upmo(S_\ep^{2l-1})$. Without lose of generality, we can assume
that the metric on $\A^l$ and on $M$ are chosen so that
$P_a\upmo(S_\ep^{2l-1})\sub\bar\cQ_a$ is the preimage of
$\cQ_a\cap\pa\Sigma$ in $\bar\cQ_a$. Hence $P_a$ induces an
${S^1_\nu}$-equivariant map
\begin{equation}\lab{15.5}
P_a\upmo(S_\ep^{2l-1})\lra S_\ep^{2l-1}
\end{equation}
and thus induces a map between their quotients
$$
\mu_a:\bar M_a\stackrel{\mathrm{def}}{=} P_a\upmo(S_\ep^{2l-1})/{S^1_\nu}\lra
S_\ep^{2l-1}/{S^1_\nu}=\PP^{l-1}_\bw.
$$
Here we use the subscript $\bw$ to indicate the weighted and the
superscript $l-1$ to denote the dimension of weighted projective
space; to be precise, we shall view the weighted projective spaces
as DM-stack. Since the specific weight is irrelevant to our study,
we shall not keep track of it in our study. We let
$$\bar\pi_a: \bar M_a\lra \tilde M_a
$$
the projection induced by $\bar\cQ_a\to\cQ_a$.

We next put our prior knowledge of the invariant part of the
obstruction sheaf of $\cQ_a$ in this setting. We let $\cT^2_{a}$ be
the obstruction sheaf on $\bar\cQ$ and let $\cT^2_{a,\soe}$ be its
invariant part. By Lemma \ref{lift-1} and \ref{A16}, there is a
canonical quotient sheaf homomorphism
\begin{equation}\lab{16.4}
\cT^2_{a,\soe}\lra\Oplus_{i=1}^l\rho_i\sta\cO_{\Ao_i},
\end{equation}
both with trivial $\soe$-actions.

A direct check shows that to each $i$ there is a
${S^1_\nu}$-linearization on $\cO_{\Ao_i}$ so that the above
homomorphism is ${S^1_\nu}$-equivariant. Because ${S^1_\nu}\cdot
S^1_\eta={T_\RR}$, the adopted ${S^1_\nu}$-linearization and the
trivial $S^1_\eta$-linearization on $\cO_{\Ao_i}$ makes (\ref{16.4})
${T_\RR}$-equivariant.

Since the obstruction sheaf $\cT^2$ on ${M}$ is a
${T_\RR}$-equivariant quotient sheaf of ${E}|_{{M}}$, pull back to
$\bar\cQ_a$, denoted by ${E}|_{\bar\cQ_a}$, and then composing with
(\ref{16.4}) give us a ${T_\RR}$-equivariant quotient sheaf
\begin{equation}\lab{16.7}
{E}|_{\bar\cQ_a}\lra\Oplus_{i=1}^l\rho_i\sta\cO_{\Ao_i}.
\end{equation}
Their descents to $\bar{M}_a$ then give rise to a quotient
homomorphism
$$
\bar{E}|_{\bar{M}_a}\lra\mu_a\sta \cV_a.
$$
Here $\cV_a$ is the descent (or the ${S^1_\nu}$-quotient) of
$\oplus_{i=1}^l\cO_{\Ao_i}|_{S^{2l-1}}$ --- a {\sl rank} $l$ vector
bundle on $\PP^{l-1}_\bw$ with trivial ${T_\RR}$-action;
$\bar{E}|_{\bar{M}_a}$ is the lift of $\tilde E|_{\tilde{M}_a}$ to
$\bar{M}_a$.

We need a key technical Lemma recently proved in \cite[Lemma
2.6]{KL} concerning the cone $C\sub E$ and its restriction to
$\cQ_a$.

\begin{lemm}[\cite{KL}]\lab{KL} Let $C|_{\cQ_a}\sub E|_{\cQ_a}$ be the
restriction of $C\sub E$ to $\cQ_a$; let $C|_{\bar\cQ_a}\sub
 E|_{\bar\cQ_a}$ be the pull back of $C|_{\cQ_a}$ to
$\bar\cQ_a$.
Then $C|_{\bar\cQ_a}$ lies in the kernel bundle of the homomorphism
\eqref{16.7}.
\end{lemm}

We are now ready to prove Theorem \ref{thm:R2}.

\begin{proof}[Proof of Theorem \ref{thm:R2}]
According the the constancy criterion Lemma \ref{lem:const}, we only
need to check that for any connected component $\tilde M_a\sub\tmt$
and for $\iota_a$ the localization homomorphism
$$\iota_a:H^\tz\lsta(\tilde M)_\fm\lmapright{\iota}
H^\tz\lsta(\tilde M,\tilde M-\tilde M_a)_\fm,
$$
we have $\iota_a [\tilde M]\virtt=0$.

To prove this vanishing, we shall first construct a
$T_\RR$-invariant finite covering
$$\pi: \bar U \to U
$$
of $\cQ_a\sub U\sub M$ that extends the existing covering
$\bar\cQ_a\to\cQ_a$; we shall then extend the existing
$\bar\cQ_a\to\A^l$ to a $T_\RR$-equivariant $\varphi\mh \bar
U\to\A^l$ and extend \eqref{16.7} to a $T_\RR$-equivariant
$$\Phi:\pi\sta E\lra \varphi\sta\bone_{\A^l}^{\oplus l},
$$
where $\bone_{\A^l}$ is the trivial line bundle on $\A^l$. Once
these are done, we then pick a $T_\RR$-invariants multi-section
$\xi$ of $\bone_{\A^l}^{\oplus l}$ on $\A^l$. We shall argue that it
can be chosen to be nowhere vanishing over $\A^l-0$.
We shall also pick a $T_\RR$-invariant $C^0$-splitting
$$
\Phi\pri: \varphi\sta\bone_{\A^l}^{\oplus l}\lra\pi\sta E
$$
of $\Phi$. Then the section $\Phi\pri(\xi)$ of $\pi\sta E$, by Lemma
\ref{KL}, is disjoint from the pull back cone $\pi\sta C$.

To define the localization $\iota_a[M]\virtt$, we shall construct an
open neighborhood $\tilde U$ of $\tilde M_a\sub\tilde M$. For this,
we form $U\cap \Sigma$; we then close the $S^1_\nu$-orbit of $U\cap
\pa\Sigma$ to form $\tilde U$. Since $U$ is a $T_\RR$-invariant
neighborhood of $\cQ_a$, $\tilde U$ is a $T_\RR$-invariant
neighborhood of $\tilde M_a\sub\tilde M$. Next we look at the
restriction to $\pi\upmo(\Sigma)\sub \bar U$ of $\Phi\pri(\xi)$,
$\Phi\pri(\xi)|_{\pi\upmo(\Sigma)}$, and its image in
$E|_{U\cap\Sigma}$ under the projection
$$\pi\sta
E|_{\pi\upmo(\Sigma)}\lra E|_{U\cap\Sigma}.
$$
Since $\Phi\pri(\xi)$ is
$T_\RR$-invariant, this restriction descends to a multi-section of
$\tilde E|_{\tilde U}$. We denote the resulting multi-section by
$\tilde \xi$. Because $\Phi\pri(\xi)\cap \pi\sta C=\emptyset$,
$\tilde\xi\cap \tilde C|_{\tilde U}=\emptyset$. Therefore,
$$0=\iota_a[M]\virtt=[\tilde\xi\cap \tilde C|_{\tilde U}]\in
H\lsta^T(\tilde U,\tilde U-\tilde M_a)_\fm= H\lsta^T(\tilde M,\tilde
M-\tilde M_a)_\fm.
$$
This will prove the vanishing $\iota_a[M]\virtt=0$.
\smallskip


We now provide the details. First, we remark that since $M$ is a
moduli of relative stable morphisms, its coarse moduli exists and is
compact. We denote the coarse moduli of $M$ by $M_-$; for
consistency, we will denote the coarse moduli of other DM-stacks
using the same subscript ``${}_-$''. Also, for any $z\in M$, we can
find a finite group $G_z$ and a $G$ space $V$ such that that stack
quotient $[V/G]$ is an open neighborhood of $z\in M$ while the
quotient $V/G$ is an open neighborhood of $z\in M_-$.

\smallskip

{\sl 1. Extending a neighborhood of $\bar\cQ_a$}. We first find an
open $\cQ_a\sub U\sub M$ and a $\Gk\times T_\RR$-orbifold $\bar U$
containing $\bar\cQ_a$ as sub-$\Gk\times T_\RR$-orbifold so that
$\bar U/\Gk=U$ and that the induced $T_\RR$-structure on $U$
coincides with the one induced by $M$.


Indeed, we first pick a $T_\RR$-invariants open $U\sub M$ containing
$\cQ_a$. we can cover $U$ by $U_\alpha$, $\alpha\in\Lambda$, so that
each $U_\alpha$ is a $G\lalp$ orbifold quotient of a $G\lalp$-space
$V\lalp$; namely $U\lalp=[V\lalp/G\lalp]$. We let $p\lalp\mh
V\lalp\to U\lalp$ be the quotient map. We next construct $\bar
U\lalp$, a $\Gk$-space whose quotient by $\Gk$ will be $U\lalp$. To
do this, we first construct a $\Gk$-space quotient $\bar U\lzalp\to
U\lzalp=U\lalp\cap \cQ_a$.

Since $M$ is the moduli of stable morphism, we have a tautological
$G\lalp$-family of stable morphisms $[\cU,\cX]$ over
$V\lzalp=p\lalp\upmo(\cQ_a)$. Using this family, we can construct a
$\Gk$-quotient $\bar V\lzalp\to V\lzalp$ similar to our construction
of $\bar\cQ$ from $\cQ$ in Lemma \ref{A10}. We claim that $G\lalp$
acts on $\bar V\lzalp$ naturally. Let $g\in G\lalp$ and $z\in
V\lalp$; then $g$ induces an isomorphism $\bar g\mh
\cX_{z}\to\cX_{gz}$. Because a point $w\in \bar U\lzalp$ over $z$
are $[\cU_z,\cX_z]$ together with $Z_i\sub \cX_z$ and
$\cU|_{Z_i}\cong\phi_{k_i,c_i}$, the isomorphism $\bar g\mh
\cX_{z}\to\cX_{gz}$ defines $\bar g(Z_i)\sub \bar g(\cX_z)$ and
$\cU|_{\bar g(Z_i)}\cong\phi_{k_i,c_i}$. This way we see that
$G\lalp$ acts on $\bar V\lzalp$ naturally.

We now construct $G\lalp$ pair $\bar V\lalp\to V\lalp$ extending
$\bar V\lzalp\to V\lzalp$. This time, since $V\lalp$ is an analytic
space, by shrinking $V\lalp$ while keep $V\lzalp$ unchanged if
necessarily, we can extend the $\Gk$-covering $\bar V\lzalp\to
V\lzalp$ to a covering $\bar V\lalp\to V\lalp$. Because this
construction is unique, the $G\lalp$ action lifts to a $G\lalp$
action on $\bar V\lalp$. We define $\bar U\lalp$ be the orbifold
quotient $[\bar V\lalp/G\lalp]$. As to the $T_\RR$-action, since it
is compact, we can make $U\lalp$ $T_\RR$-invariant; thus $\bar
V\lalp\to V\lalp$ is $T_\RR$-equivariant. Finally, since this
construction is canonical once we have the chart $V\lalp\to U\lalp$,
the orbifold quotient $\bar U\lalp$ patch together to form an
$\Gk\times T_\RR$-orbifold $\pi\mh\bar U\to U$. Lastly, we comment
that $\bar U$ has a coarse moduli the the coarse moduli $U_-$ is the
quotient $\bar U_-/\Gk$.

\smallskip

{\sl 2. Extending the map $\bar\cQ_a\to \A^l$}. We now show that the
map $\rho\mh\bar\cQ_a\to \A^l$ constructed in Proposition \ref{A15}
extends to a $T_\RR$-equivariant $C^0$-map $\varphi\mh\bar U\lra
\A^l$. Here the $T_\RR$ action on $\A^l$ is induced by the trivial
$S^1_\nu$ action and the $S^1_\eta$ action defined in \eqref{15.3}.
Since $\rho$ factors through its course moduli $\rho_-\mh
\bar\cQ_{a-}\to\A^l$, to extend $\rho$, we only need to extend
$\rho_-$ to a $T_\RR$-equivariant continuous $\varphi_-\mh \bar
U_-\to \A^l$. This is possible because $T_\RR$ is compact and that
$T_\RR$ acts on $\A^l$ linearly. Once we have $\varphi_-$, we define
$$\varphi: \bar U\lra \A^l
$$
be the composite of $\bar U\to \bar U_-$ with $\varphi_-$.

\smallskip
{\sl 3. Extending the homomorphism $E|_{\bar\cQ_a}\to
\rho\sta\bone_{\A^l}^{\oplus l}$}. Lastly, we shall extend the sheaf
homomorphism \eqref{16.7} to a $T_\RR$-equivariant $\Phi\mh\pi\sta
E\lra  \varphi\sta\bone_{\A^l}^{\oplus l}$. This is possible for the
same as the previous extension problem.

To complete the proof, we need to construct an invariant hermitian
metric on $\pi\sta E$ that will provide us a splitting $\Phi\pri$.
This is simple using partition of unity over $U_-$. First, we find
hermitian metric $h\lalp$ on $p\lalp\sta E$ over $V\lalp$. Since
both $T_\RR$ and $G\lalp$ are compact, we can assume $h\lalp$ is
$G\lalp\times T_\RR$-invariant. we then use a partition of unity
$\zeta\lalp$ subordinate to the covering $U_{\alpha-}\sub U_i$ and
define the hermitian metric on $\pi\sta E$ be the sum of the pull
back of $\zeta\lalp h\lalp$. The resulting metric $h$ satisfies the
desired property.

Finally, we need to find a $T_\RR$-invariant multisection of
$\bone_{\A^l}^{\oplus l}$ which is nowhere vanishing over $\A^l-0$.
Since the $T_\RR$ action is via the trivial $S^1_\eta$ and the $S^1_\nu$
action specified in \eqref{15.3}, we only need to find an
$S_\nu^1$-invariant multiple section.
The action of $S^1_\nu$
on $\bone_{\A^l}^{\oplus l}$ is given by
$$
(z_1,\ldots,z_l)^\sigma = (\sigma^{v_1} z_1,\ldots, \sigma^{v_l} z_l).
$$
for some $(v_1,\ldots,v_l)$. For $i=1,\ldots,l$, let
\begin{equation}\label{eqn:s}
s_i(u_1,\ldots,u_l)=\begin{cases}
u_i^{v_i/w_i} & v_i/w_i>0,\\
\bar u_i^{-v_i/w_i} & v_i/w_i<0,\\
1 & v_i=0.
\end{cases}
\end{equation}
Then $(s_1,\ldots,s_l)$ is an $S_\nu^1$-invariant multiple section of
$\bone_{\A^l}^{\oplus l}$ which is nonwhere vanishing over $\A^l-0$.
\end{proof}

%% file: sec6.tex
\section{Topological Vertex, Hodge Integrals and Double Hurwitz numbers}
\label{sec:hodge}

Let $\Ga_{\bn;w_1,w_2}$ be the FTCY in Figure 10, where 
\begin{equation}\label{eqn:fw}
f_1=w_2-n_1 w_1,\quad 
f_2=w_3-n_2 w_2,\quad
f_3=w_1-n_3 w_3,\quad  
w_3=-w_1-w_2,
\end{equation}
and $\bn=\three{n}\in\ZZ^{\oplus 3}$.
Any topological vertex (defined in
Definition \ref{def:vertex}) is of this form. 

\begin{figure}[h] \label{figure10}
\begin{center}
\psfrag{v0}{\footnotesize $v_0$} 
\psfrag{v1}{\footnotesize $v_1$} 
\psfrag{v2}{\footnotesize $v_2$}
\psfrag{v3}{\footnotesize $v_3$} 
\psfrag{p1}{\footnotesize $\fp(e_1)=w_1$}
\psfrag{p2}{\footnotesize $\fp(e_2)=w_2$} 
\psfrag{p3}{\footnotesize $\fp(e_3)=w_3$}
\psfrag{f1}{\footnotesize $\ff(e_1)=f_1$} 
\psfrag{f2}{\footnotesize $\ff(e_2)=f_2$}
\psfrag{f3}{\footnotesize $\ff(e_3)=f_3$}
\includegraphics[scale=0.7]{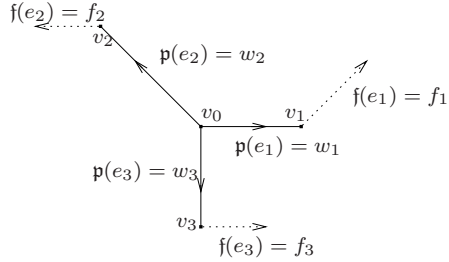}
\end{center}
\caption{The graph of a topological vertex }
\end{figure}

In this section, we will compute
\begin{equation}\label{eqn:Fn}
F^\bu_{\chi,\vmu}(\bn;w_1,w_2)\stackrel{\mathrm{def}}{=}F^{\bu\Ga_{\bn;w_1,w_2}}_{\chi,\vd,\vmu}(u_1,u_2)
\end{equation}
where the RHS is defined by \eqref{eqn:FN}.
To simplify the notation, we will fix $\bn=\three{n}$ and $(w_1,w_2)$ and
write $\Ga$ instead of $\Ga_{\bn;w_1,w_2}$.

\subsection{Torus Fixed Points and Label Notation}
\label{sec:vertex-loci}
In this subsection, we describe the $T$-fixed points in
$\Mmu\stackrel{\mathrm{def}}{=}\Mfml$, and introduce the label notation. 
To each label corresponds a disjoint union of connected components of
$$
\Mmu^T=\Mmu^{T_\RR},
$$
or equivalently, a collection of graphs in the graph notation.

Let $\Yr=(\hY,\hD)$ be the FTCY associated to $\Ga$, and
let 
$$
\hD^i=\hD^{v_i},\quad C^i=C^{\ee_i}
$$
for $i=1,2,3$. Given $u:(X, \bq)\lra \bl\hY_\bm, \hD_\bm \br$
which represents a point in $\Mmu^T$, let $\tu=\pi_\bm\circ u: X\to
\hatYrel$, where $\pi_\bm:\hY_\bm\to \hY$ be the projection defined
in Section \ref{sec:moduli}. Then $\tu(X)\subset C^1\cup C^2 \cup C^3$.
Let $z^0$ and $z^i$ be the two T fixed points on $C^i$, and let
$$
V^i=\tu^{-1}(z^i)
$$
for $i=0,1,2,3$. Let $E^i$ be the closure of
$\tu^{-1}(C^i\setminus\{z^0,z^i\})$ for $i=1,2,3$.
Then $E^i$ is a union of projective lines, and
$u|_{E^i}:E^i\to C^i$ is a degree $d^i=|\mu^i|$ cover
fully ramified over $z^0$ and $z^i$.

Define
$$
\PP^i(m^i)=\pi_\bm^{-1}(z^i)
$$
which is a point if $m^i=0$, and is a chain of $m^i$ copies of $\Po$
if $m^i>0$.

For $i=1,2,3$, let
$$
\hat{u}^i = u|_{V^i}:V^i \lra \PP^i(m^i),\quad
\tu^i = u|_{E^i}:E^i \lra C^i.
$$
The degrees of $\tu^i$ restricted to connected components
of $E^i$ determine a partition $\nu^i$ of $d^i$.

For $i=0,1,2,3$, let
$V_1^i,\ldots, V_{k^i}^i$ be the connected components of $V^i$, and
let $g_j^i$ be the arithmetic genus of $V_j^i$. (We define
$g_j^i=0$ if $V_j^i$ is a point.) Define
$$
\chi^i=\sum_{j=1}^{k^i}(2-2g^i_j).
$$

Then
$$
-\sum_{i=0}^3\chi^i + 2\sum_{i=1}^3 \ell(\nu^i) = -\chi.
$$

Note that $\chi^i \leq 2\min \{\ell(\mu^i),\ell(\nu^i)\}$
for $i=1,2,3$, so
$$
\xnm\geq 0
$$
and the equality holds if and only if $m^i=0$. In this case, we have
$\nu^i=\mu^i$, $\chi^i= 2\ell(\mu^i)$.

We introduce moduli spaces of relative stable maps to the non-rigid $(\Po,\{0,\infty\})$
(called {\em rubber} in \cite{Oko-Pan3} etc.):
$$
\tMP\stackrel{\mathrm{def}}{=}\overline{\cM}^{\bu}_{\chi}(\Po,\nu,\mu)//\CC^*
$$
where $\overline{\cM}^\bu_\chi(\PP^1,\nu,\mu)//\CC^*$ is defined
as in \cite[Section 5]{LLZ2}.  

For each $i\in \{1,2,3\}$, there are two cases:
\paragraph{Case 1: $m^i=0$.} Then $\hu^i$ is a constant map
from $\ell(\mu^i)$ points to $p^i$.
\paragraph{Case 2: $m^i>0$.} Then $\hu^i$ represents a point in
$\tMi$.

\begin{defi}
An {\em admissible label} of $\Mmu$ is
a pair $\pair$ such that
\begin{enumerate}
\item $\vx=(\chi^0,\up{\chi})$, where $\chi^i\in 2\ZZ$.
\item $\vnu=(\up{\nu})$, where $\nu^i$ is a partition such
      that $|\nu^i|= |\mu^i|$.
\item $\chi^0\leq 2\sum_{i=1}^3 \ell(\nu^i)$.
\item $\chi^i\leq 2\min\{\ell(\mu^i), \ell(\nu^i)\}$ for
      $i=1,2,3$.
\item $-\sum_{i=0}^3\chi^i + 2\sum_{i=1}^3 \ell(\nu^i)=-\chi$.
\end{enumerate}
Let $\GYfm$ denote the set of all admissible labels of $\Mmu$.
 \end{defi}
For a nonnegative integer $g$ and a positive integer $h$,
let $\Mbar_{g,h}$ be the moduli space of
stable curves of genus $g$ with $h$ marked points.
Although $\Mbar_{g,h}$ is empty for $(g,h)=(0,1), (0,2)$,
for simiplicity of notation we will
formally assume the following integrals exist:
$$
\int_{\Mbar_{0,1}}\frac{1}{1-d\psi}=\frac{1}{d^2}\ ,\quad\quad
\int_{\Mbar_{0,2}}\frac{1}{(1-\mu_1\psi_1)(1-\mu_2\psi_2)}
= \frac{1}{\mu_1+\mu_2}.
$$
This convention will give the correct final results.

For a nonnegative integer $g$ and a positive integer $h$,
let $\Mbar^\bu_{\chi,h}$ be the moduli of possibly
disconnected  stable curves $C$  with $h$ marked points such that
\begin{itemize}
\item If $C_1,\ldots,C_k$ are connected components of $C$, and $g_i$ is the
  arithmetic genus of $C_i$, then
$$
\sum_{i=1}^k (2-2g_i)=\chi.
$$
\item Each connected component contains at least one marked point.
\end{itemize}

The connected components of $\Mbar^\bu_{\chi,h}$ are of the form
$$
\Mbar_{g_1,h_1}\times\cdots\times\Mbar_{g_k,h_k}.
$$
where
$$
\sum_{i=1}^k(2-2g_i)=\chi,\ \ \ \sum_{i=1}^k h_i=h.
$$
The restriction of the Hodge bundle $\EE\to \Mbar^\bu_{\chi,h}$ to
the above connected component is the direct sum of the Hodge bundles
on each factor, and
$$
\Lambda^\vee(u)=\prod_{i=1}^k\Lambda_{g_i}^\vee(u).
$$

We define 
$$
\xn{\Mbar} = \prod_{i=0}^3\xn{\Mbar}^i
$$
where $\xn{\Mbar}^0=\Mbar^\bu_{\chi^0,\ell(\vnu)}$, and 
for $i\in\{ 1,2,3 \}$,
$$ 
\xn{\Mbar}^i = \begin{cases}
\{\mathrm{pt} \}, & \xnm=0,\\
\tMi, & \xnm>0.
\end{cases}
$$

For each $\pair \in \GYfm$, there is a morphism
$\xn{i}:\xn{\Mbar}\to \Mmu^T$,
whose image $\xn{\cF}$ is a union of connected components of $\Mmu^T$.
The morphism $\xn{i}$ induces an isomorphism
$$
\xn{\Mbar}\biggl/\biggl(\, \prod_{i=1}^3\xn{A}^i \biggr)\cong \xn{\cF}
$$
where $\xn{A}^i$ is the automorphism group associated to the edge $e_i$:
\begin{eqnarray*}
\xn{A}^i=\prod_{j=1}^{\ell(\nu^i)}\ZZ_{\nu^i_j}, 
&& \xnm=0;\\
1\to \prod_{j=1}^{\ell(\nu^i)}\ZZ_{\nu^i_j} \to \xn{A}^i \to \Aut(\nu^i) \to 1,
&& \xnm>0.
\end{eqnarray*}

The fixed points set $\Mmu^T$ is a disjoint union of
$$
\{ \xn{\cF} \mid \pair\in\GYfm \}.
$$

\begin{rema}
There are two perfect obstruction theories on $\xn{\cF}$: one
is the fixed part $[\cT^{1,f}\to \cT^{2,f}]$ of the restriction
of the perfect obstruction theory on $\Mmu$; the other comes from the perfect
obstruction theory on the moduli spaces $\Mbar^\bu_{\chi^0,\ell(\vnu)}$ and $\tMi$. 
It is straightforward to check that they coincide.
\end{rema}

\subsection{Contribution from Each Label}
\label{sec:vertex-label}
We view $w_i$ and $f_i$ in
Equation \eqref{eqn:fw} as elements in 
$$
\ZZ u_1 \oplus \ZZ u_2 = \Lambda_T \cong H^2_T(\mathrm{pt},\QQ).
$$ 
Recall that $H^*_T(\mathrm{pt};\QQ)=\QQ[u_1,u_2]$.
The results of localization calculations will involve rational functions of
$w_i$ and $f_i$ which are elements in  $\QQ(u_1,u_2)$.

If $m^i>0$, let $\psi^0_i,\psi^\infty_i$ denote the target $\psi$ class of
$\Mbar^i_{\vx,\vnu}$ (see e.g. \cite[Section 5]{LLZ2} for definitions).
Let $N\virt_{\vx,\vnu}$ denote
the virtual bundle on $\Mbar_{\vx,\vnu}$ which is
the pull back of $\cT^{1,m}-\cT^{2,m}$ under $i_{\vx,\vmu}$.

With the above notation and the explicit description of 
$[\cT^1\to \cT^2]$ in Section \ref{sec:tan-obs}, calculations similar to those 
in \cite[Appenix A]{LLZ1} show that
$$
\frac{1}{e_T(N_{\vx,\vnu}\virt)}
= \prod_{i=0}^3 B_{v_i}\prod_{i=1}^3 B_{e_i},
$$
where
$$
B_{v_0}=\prod_{i=1}^3
\frac{a_{\nu^i}\Lambda^\vee(w_i)w_i^{\ell(\vnu)-1}}
{\prod_{j=1}^{\ell(\nu^i)}
(w_i(w_i-\nu^i_j\psi^i_j))},
$$
and for $i\in \{ 1,2,3 \}$, 
\begin{eqnarray*}
B_{v_i}&=&
\begin{cases}
 1,& \xnm=0\\
{\displaystyle (-1)^{\ell(\nu^i)-\chi^i/2} a_{\nu^i}
\frac{f_i^{\xnm}}{-w_i-\psi^0_i}, } & \xnm>0 
\end{cases}\\
B_{e_i}&=& (-1)^{|\nu^i|n^i+\ell(\nu^i)-|\nu^i|}\prod_{j=1}^{\ell(\nu^i)}
\frac{\prod_{a=1}^{\nu^i_j-1}
(w_{i+1}\nu^i_j+ a w_i)}{(\nu^i_j -1)! w_i^{\nu^i_j -1} }
\end{eqnarray*}
The disconnected double Hurwitz numbers $H^\bu_{\chi,\nu,\mu}$
(see Section \ref{sec:H}) can be related to intersection
of the target $\psi$ class (see \cite[Section 5]{LLZ2} for
a derivation): 
\begin{equation}\label{eqn:H}
H^\bu_{\chi,\nu,\mu}=\frac{(-\chi+\ell(\nu)+\ell(\mu))!}{|\Aut(\nu)\times\Aut(\mu)|}
\int_{[\tMP]\virt}(\psi^0)^{-\chi+\ell(\nu)+\ell(\mu)-1}
\end{equation}
The three-partition Hodge integral $G^\bu_{\chi,\vmu}(\bw)$ defined
by \eqref{eqn:Gchimu} in Section \ref{sec:G} can be expressed as
\begin{equation}\label{eqn:G}
G^\bu_{\chi,\vnu}(\bw)=(-\sqrt{-1})^{\ell(\vnu)} V_{\chi,\vnu}(\bw)
\prod_{i=1}^3 E_{\nu^i}(w_{i+1},w_i)
\end{equation}
where
\begin{equation}\label{eqn:V}
V_{\chi,\vnu}(\bw)
=\frac{1}{|\Aut(\vnu)|}\int_{\Mbar^\bu_{\chi,\ell(\vnu)}}
\prod_{i=1}^3\frac{\Lambda^\vee(w_i)w_i^{\ell(\vnu)-1}}
{\prod_{j=1}^{\ell(\nu^i)}(w_i(w_i-\nu^i_j \psi^i_j)) }
\end{equation}
\begin{equation}\label{eqn:E}
E_\nu(x,y)=\prod_{j=1}^{\ell(\nu)}
\frac{\prod_{a=1}^{\nu_j-1}(y\nu_j +ax)}{(\nu_j-1)!x^{\nu_j-1} }
\end{equation}

Set
$$
\xn{I}(\bn;\bw)=
\int_{[\xn{\cF}]\virt}\frac{1}{ e_T(N_{\vx,\vnu}\virt)} 
$$
Then
\begin{eqnarray*}
&& \xn{I}(\bn;\bw)= \frac{1}{\prod_{i=1}^3|\xn{A}^i|}\int_{[\xn{\Mbar}]\virt}
\frac{1}{e_T(N_{\vx,\vnu}\virt)}\\
&=& |\Aut(\vmu)|(-1)^{\sum_{i=1}^3 (n_i-1)|\mu^i|}
\bigl(-\sqrt{-1}\bigr)^{\ell(\vmu)+\ell(\vnu)}
   V_{\chi^0,\vnu}(\bw)\\
&& \cdot \prod_{i=1}^3 E_{\nu^i}(w_i,w_{i+1}) z_{\nu^i}
   \Bigl(-\sqrt{-1}\frac{f_i}{w_i}\Bigr)^{\xnm}
   \frac{H^\bu_{\chi^i,\nu^i,\mu^i} }{(\xnm)!}
\end{eqnarray*}
So
\begin{equation}\label{eqn:IG}
\begin{aligned}
& \xn{I}(\bn;\bw)=
|\Aut(\vmu)|(-1)^{\sum_{i=1}^3 (n_i-1)|\mu^i|} \bigl(-\sqrt{-1}\bigr)^{\ell(\vmu)}
   G^\bu_{\chi^0,\vnu}(\bw)\\
&\, \cdot\prod_{i=1}^3 z_{\nu^i}
  \Bigl(\sqrt{-1}\Bigl(n_i-\frac{w_{i+1}}{w_i}\Bigr)\Bigr)^{-\chi^i+\ell(\nu^i)+\ell(\mu^i)}
   \frac{H^\bu_{\chi^i,\nu^i,\mu^i}}{(-\chi^i+\ell(\nu^i)+\ell(\mu^i))!}.
\end{aligned}
\end{equation}

\subsection{Sum over Labels} 
\label{sec:vertex-sum}
We have
$$
F^\bu_{\chi,\vmu}(\bn;w_1,w_2)=\frac{1}{|\Aut(\vmu)|}
\sum_{\pair\in\GYfm}\xn{I}(\bn;\bw).
$$
Define generating functions
\begin{equation}\label{eqn:F}
F^\bu_\vmu(\lam;\bn;w_1,w_2)=
\sum_{\chi\in 2\ZZ,\chi\leq \ell(\vmu)} \lam^{-\chi+\ell(\vmu)}
F^\bu_{\chi,\vmu}(\bn;w_1,w_2)
\end{equation}
\begin{equation}\label{eqn:tF}
\tF^\bu_\vmu(\lam;\bn;w_1,w_2)=
(-1)^{\sum_{i=1}^3 (n_i-1)|\mu^i|}\sqrt{-1}^{\ell(\vmu)}
F^\bu_{\vmu}(\lam;\bn;w_1,w_2)
\end{equation}
Then relation \eqref{eqn:IG} becomes
\begin{equation} \label{eqn:tFG}
\tF^\bu_\vmu(\lam;\bn;w_1,w_2)
=\sum_{|\nu^i|=|\mu^i|}G^\bu_\vnu(\lam;\bw)
\prod_{i=1}^3 z_{\nu^i}\Phi^\bu_{\nu^i,\mu^i}
\Bigl(\sqrt{-1}\Bigl(n_i-\frac{w_{i+1}}{w_i}\Bigr)\lam\Bigr),
\end{equation}
where $G^\bu_\vmu(\lam;\bw)$ is defined by \eqref{eqn:Gmu} in Section \ref{sec:G},
and $\Phi^\bu_{\nu,\mu}(\lam)$ is the generating function of 
disconnected double Hurwitz numbers defined in Section \ref{sec:H}.

Equations (\ref{eqn:tFG}), (\ref{eqn:Hsum}), and (\ref{eqn:Hinitial}) imply that
\begin{equation} \label{eqn:tFn}
\tF^\bu_\vmu(\lam;\bn;w_1,w_2)
=\sum_{|\nu^i|=|\mu^i|}\tF^\bu_\vnu(\lam;\mathbf{0},w_1,w_2)
\prod_{i=1}^3 z_{\nu^i}\Phi^\bu_{\nu^i,\mu^i}
\left(\sqrt{-1}n_i\lam\right)
\end{equation}
By Theorem \ref{thm:tri-inv},
$$
F^\bu_\vnu(\lam;\mathbf{0};w_1,w_2)=\sum_\chi \lam^{-\chi+\ell(\vnu)}F^{\bu\Ga^0}_{\chi,\vd,\vmu}(w_1,w_2)
$$
does not depend on $w_1,w_2$. So by (\ref{eqn:tF})
and (\ref{eqn:tFn}), $F^\bu_\vmu(\lam;\bn;w_1,w_2)$ and
$\tF^\bu_\vmu(\lam;\bn;w_1,w_2)$ do not depend on $w_1,w_2$.
From now on, we will write
$$
F^\bu_\vmu(\lam;\bn),\quad
\tF^\bu_\vmu(\lam;\bn)
$$
instead of $F^\bu_\vmu(\lam;\bn;w_1,w_2)$,
$\tF^\bu_\vmu(\lam;\bn;w_1,w_2)$.
To summarize, for each $\vmu\in\cP_+^3$ and each $\bn\in\ZZ^3$, we have
defined an generating function $F^\bu_\vmu(\lam;\bn)$ which can be
expressed in terms of Hodge integrals and double Hurwitz numbers
as follows.
\begin{prop}\label{thm:FG}
\begin{equation} \label{eqn:FG}
\tF^\bu_\vmu(\lam;\bn)
=\sum_{|\nu^i|=|\mu^i|}G^\bu_\vnu(\lam;\bw)
\prod_{i=1}^3 z_{\nu^i}\Phi^\bu_{\nu^i,\mu^i}
\Bigl(\sqrt{-1}\Bigl(n_i-\frac{w_{i+1}}{w_i}\Bigr)\lam\Bigr),
\end{equation}
\end{prop}

Proposition \ref{thm:FG} and the sum formula (\ref{eqn:Hsum}) of double
Hurwitz numbers imply:
\begin{coro}[framing dependence in winding basis]
\begin{equation} \label{eqn:tFframe}
\tF^\bu_\vmu(\lam;\bn)
=\sum_{|\nu^i|=|\mu^i|}\tF^\bu_\vnu(\lam;\mathrm{0})
\prod_{i=1}^3 z_{\nu^i}\Phi^\bu_{\nu^i,\mu^i}
\left(\sqrt{-1}n_i\lam\right)
\end{equation}
\end{coro}
Note that (\ref{eqn:tFframe}) is valid for any three complex numbers
$\lo{n}$.

\subsection{Representation Basis}
The framing dependence (\ref{eqn:tFframe}) is particularly simple
in the representation basis used in \cite{AKMV}. We use the notation
in Section \ref{sec:W}. Define
$\tC_{\mu}(\lam;\bn)$ by
\begin{equation}\label{eqn:tCtF}
\tC_{\vmu}(\lam;\bn)=\sum_{|\nu^i|=|\mu^i|}\tF^\bu_\vnu(\lam;\bn)\prod_{i=1}^3
\chi_{\mu^i}(\nu^i)
\end{equation}
which is equivalent to
\begin{equation}\label{eqn:tFtC}
\tF_{\vmu}(\lam;\bn)=\sum_{|\nu^i|=|\mu^i|}\tC^\bu_\vnu(\lam;\bn)\prod_{i=1}^3
\frac{\chi_{\nu^i}(\mu^i)}{z_{\mu^i}}.
\end{equation}
Then (\ref{eqn:tFframe}) is equivalent to
\begin{prop}[framing dependence in representation basis]\label{thm:tCframe}
\begin{equation}\label{eqn:tCframe}
\tC_\vmu(\lam;\bn)=
e^{\frac{1}{2}\sqrt{-1}(\sum_{i=1}^3 \kappa_{\mu^i}n_i)\lam} \tC_\vmu(\lam;\mathbf{0}).
\end{equation}
\end{prop}
Define $\tC_\vmu(\lam)=\tC_\vmu(\lam;\mathbf{0})$, and let $q=e^{\sqrt{-1}\lam}$.
Then (\ref{eqn:tCtF}), (\ref{eqn:FG}), and the Burnside formula
(\ref{eqn:burnside}) of double Hurwitz numbers imply
\begin{prop}\label{thm:GtC}
We have
\begin{equation}\label{eqn:tCG}
\tC_{\vmu}(\lam)=q^{-\frac{1}{2}(\sum_{i=1}^3\kappa_{\mu^i}\frac{w_{i+1}}{w_i})}
\sum_{|\nu^i|=|\mu^i|}G^\bu_\vnu(\lam;\bw)\prod_{i=1}^3\chi_{\mu^i}(\nu^i).
\end{equation}
\begin{equation}\label{eqn:GtC}
G^\bu_\vmu(\lam;\bw)=\sum_{|\nu^i|=|\mu^i|}\prod_{i=1}^3\frac{\chi_{\nu^i}(\mu^i)}{z_{\mu^i}}
q^{\frac{1}{2}(\sum_{i=1}^3 \kappa_{\nu^i}\frac{w_{i+1}}{w_i})}\tC_\vnu(\lam).
\end{equation}
\end{prop}

%% file: sec7.tex
\section{Gluing Formulae of Formal Relative  Gromov-Witten Invariants}\label{sec:glue}

Let $\Ga$ be a FTCY graph (see Definition \ref{def:FTCYgraph}),
and let $(\vd,\vmu)$ be an effective class 
of $\Gamma$ (defined in Definition \ref{def:effective}).
In this section, we will calculate the formal relative Gromov-Witten
invariant
\begin{equation}
F^{\bu\Ga}_{\chi,\vd,\vmu}(u_1,u_2)\in\QQ(u_2/u_1).
\end{equation}
We will reduce the invariance of
$F^{\Gamma\bu}_{\chi,\vd,\vmu}$ (Theorem \ref{thm:R2}) to
the invariance of the topological vertex at the standard framing
(Theorem \ref{thm:tri-inv}). We will
derive gluing formulae for such invariants.

As in Definition \ref{def:effective}, we will use the abbreviation
$$
d^{\ee}=\vd(\ee),\quad \ee\in E(\Ga); \
\mu^v=\vmu(v),\quad v\in V_1(\Ga).
$$
\subsection{Torus Fixed Points and Label Notation}\label{sec:cy-label}
In this subsection, we describe the $T$-fixed points in
$\Mfml$, and introduce the label notation.
This is a generalization of Section \ref{sec:vertex-loci}.

Given a morphism
$$
u:(X,\bq)\lra (\hY_\bm,\hD_\bm)
$$
which represents a point in $\Mfml^T$, let $\tu=\pi_\bm\circ u: X\to \hY$,
as before. Then
$$
\tu(X)\subset \bigcup_{\ee\in E(\Ga)}C^\ee.
$$
where $C^\ee$ is defined as in Section \ref{sec:threefold}.

Let $z^v$ be the $T$ fixed point associated to $v\in V(\Gamma)$, as in Section
\ref{sec:threefold}, and let
$$
V^v=\tu^{-1}(z^v).
$$
Let $E^\ee$ be the closure of
$\tu^{-1}(C^\ee\setminus \{z^{\fv_0(e)},z^{\fv_1(e)} \} )$ for $\ee=\{e,-e\}\in E(\Gamma)$.
Then $E^\ee$ is a union of projective lines, and
$u|_{E^\ee}:E^\ee\to C^\ee$ is a degree $d^{\ee}$ cover
fully ramified over $\fv_0(e)$ and $\fv_1(e)$.

For $v\in V_1(\Ga)\cup V_2(\Ga)$, define
$$
\PP^v(m^v)=\pi_\bm^{-1}(z^v)
$$
which is a point if $m^v=0$, and is a chain of $m^v$ copies of $\Po$ if
$m^v>0$. Let
$$
\hat{u}^v = u|_{V^v}:V^v \to \PP^v(m^v).
$$

For $\ee\in E(\Gamma)$, define
$$
\tu^\ee = u|_{E^\ee}:E^\ee \to C^\ee.
$$
The degrees of $\tu^\ee$ restricted to connected components
of $E^\ee$ determine a partition $\nu^e=\nu^{-e}$ of $d^{\ee}$.

For $v\in V(\Gamma)$, let
 $V_1^v,\ldots, V_{k^v}^v$ be the connected components of $V^v$, and
let $g_j^v$ be the arithmetic genus of $V_j^v$. (We define
$g_j^v=0$ if $V_j^v$ is a point.) Define
$$
\chi^v=\sum_{j=1}^{k^v}(2-2g^v_j).
$$
Then
$$
-\sum_{v\in V(\Gamma)} \chi^v + \sum_{e\in E\ori(\Gamma)} \ell(\nu^e) = -\chi.
$$

Given $v\in V_1(\Ga)$ with $\fv_1^{-1}(v)=\{e\}$, we have 
$\chi^v \leq 2\min \{\ell(\nu^e),\ell(\mu^v)\}$, so 
\begin{equation}\label{eqn:rVone}
r^v \stackrel{\mathrm{def}}{=} -\chi^v +\ell(\nu^e)+\ell(\mu^v) \geq 0
\end{equation}
and the equality holds if and only if $m^v=0$. In this case, we have
$\nu^e=\mu^v$, $\chi^v= 2\ell(\mu^v)$.
For each $v\in V_1(\Ga)$, there are two cases:
\paragraph{Case 1: $m^v=0$.} Then $\hu^v$ is a constant map
from $\ell(\mu^v)$ points to $z^v$.
\paragraph{Case 2: $m^v>0$.} Then $\hu^v$ represents a point in
$\Mbar^{\bu\sim}_{\chi^v,\nu^e,\mu^v}$.

Given $v\in V_2(\Ga)$ with $\fv_1^{-1}(v)=\{e,e'\}$, we 
have $\chi^v \leq 2\min \{\ell(\nu^e,\ell(\nu^{e'})\}$, so 
\begin{equation}\label{eqn:rVtwo}
r^v\stackrel{\mathrm{def}}{=} -\chi^v+ \ell(\nu^e)+\ell(\nu^{e'})\geq 0
\end{equation}
and the equality holds if and only if $m^v=0$. In this case, we have
$\nu^e=\nu^{e'}$, $\chi^v= 2\ell(\nu^e)$.
For each $v\in V_2(\Ga)$, there are two cases:
\paragraph{Case 1': $m^v=0$.} Then $\hu^v$ is a constant map
from $\ell(\mu^v)$ points to $z^v$.
\paragraph{Case 2': $m^v>0$.} Then $\hu^v$ represents a point in
$\Mbar^{\bu\sim}_{\chi^v,\nu^e,\nu^{e'}}$.

\begin{defi}
An {\em admissible label} of $\Mfml$ is
a pair $\pair$ such that
\begin{enumerate}
\item $\vx: V(\Ga)\lra 2\ZZ$. Let $\chi^v$ denote $\vx(v)$.
\item $\vnu: E\ori(\Ga)\lra \cP$, where $\vnu(e)=\vnu(-e)$ and $|\vnu(e)|=d^{\ee}$.
      We write $\nu^e$ for $\vnu(e)$.
\item For $v\in V_1(\Ga)$ with $\fv_1^{-1}(v)=\{e\}$, we have
      $\chi^v\leq 2\min\{\ell(\nu^e), \ell(\mu^v)\}$.
\item For $v\in V_2(\Ga)$ with $\fv_1^{-1}(v)=\{e,e'\}$, we have
      $\chi^v\leq 2\min\{\ell(\nu^e), \ell(\nu^{e'})\}$.
\item For $v\in V_3(\Gamma)$, define
      $\ell_{\vnu}(v)=\sum_{e\in\fv_0^{-1}(v)}\ell(\nu^e)$. Then $\chi^v\leq 2 \ell_{\vnu}(v)$.
\item $-\sum_{v\in V(\Ga)}\chi^v +2\sum_{e\in E(\Gamma)}\ell(\nu^e)=-\chi$.
\end{enumerate}
Let $\Gdmu$ denote the set of all admissible labels of $\Mfml$.
 \end{defi}

Given $\pair\in\Gdmu$, define $r^v$ as in (\ref{eqn:rVone}) and
(\ref{eqn:rVtwo}) for $v\in V_1(\Ga)$ and $v\in V_2(\Ga)$, respectively. 
We define
$$
\xn{\Mbar} = \prod_{v\in V(\Gamma)}\xn{\Mbar}^v
$$
where
$$
\xn{\Mbar}^v =\begin{cases}
\{\mathrm{pt} \}, & v\in V_1(\Ga)\cup V_2,\ r^v=0,\\
\Mbar^{\bu\sim}_{\chi^v,\nu^e,\mu^v}, 
& v\in V_1(\Ga),\ \fv_1^{-1}(v)=\{ e\},\ r^v>0,\\
\Mbar^{\bu\sim}_{\chi^v,\nu^e,\nu^{e'}}, 
& v\in V_2(\Ga),\ \fv_1^{-1}(v)=\{e,e'\},\ r^v>0,\\
\Mbar^\bu_{\chi^v,\ell_\vnu(v)}, & v\in V_3(\Ga) .
\end{cases}
$$

For each $\pair \in \Gdmu$, there is a morphism
$\xn{i}:\xn{\Mbar}\to \Mdmu^T$,
whose image $\xn{\cF}$ is a union of connected components of
$\Mdmu^T$. The morphism $\xn{i}$ induces an isomorphism
$$
\xn{\Mbar}\biggl/\biggl(\ \prod_{\ee\in E(\Gamma)}\xn{A}^\ee\ \biggr)\cong \xn{\cF},
$$
where $\xn{A}^\ee$ is the automorphism group associated to the edge $\ee$:
\begin{eqnarray*}
&& \xn{A}^\ee=\prod_{j=1}^{\ell(\nu^e)}\ZZ_{\nu^e_j}\;, \quad
\{\fv_0(e),\fv_1(e)\}\cap V_1(\Ga)=\{v\}\neq \emptyset \quad \textup{and} \quad r^v=0;\\
&& 1\to \prod_{j=1}^{\ell(\nu^e)}\ZZ_{\nu^e_j} \to \xn{A}^\ee \to \Aut(\nu^e) \to 1\; ,\quad
\textup{otherwise.}
\end{eqnarray*}

The fixed points set $\Mfml^T$ is a disjoint union of
$$
\{ \xn{\cF} \mid \pair\in\Gdmu \}.
$$

\subsection{Perfect Obstruction Theory on Fixed Points Set}

There are two perfect obstruction theories on $\xn{\cF}$:
one is the fixing part $[\cT^{1,f}\to \cT^{2,f}]$ of
the restriction of the perfect obstruction theory on $\Mfml$; 
the other comes from the perfect obstruction theory on the moduli spaces
$\Mbar^\bu_{\chi^v,\ell_\vnu(v)}$ and  $\tMP$.
Let $[\xn{\Mbar}]\virt$ denote the virtual cycle defined by
$[\cT^{1,f}\to \cT^{2,f}]$. By inspecting the $T$-action
on the perfect obstruction theory on $\Mfml$ (see \cite{Li2}
and the description in Section \ref{sec:invariants}), we get
$$
[\xn{\Mbar}]\virt=\prod_{v\in V(\Ga)}[\xn{\Mbar}^v]\virt
$$
where
$$
[\xn{\Mbar}^v]\virt =\begin{cases}
[\{\mathrm{pt} \}], & v\in V_1(\Ga)\cup V_2(\Ga),\ r^v=0,\\
{[\Mbar^{\bu\sim}_{\chi^v,\nu^e,\mu^v}]}\virt, 
& v\in V_1(\Ga),\ \fv_1^{-1}(v)=\{ e\},\ r^v>0,\\
c_1(\LL)\cap {[\Mbar^{\bu\sim}_{\chi^v,\nu^e,\nu^{e'}}]}\virt, 
& v\in V_2(\Ga),\ \fv_1^{-1}(v)=\{e,e'\},\ r^v>0,\\
[\Mbar^\bu_{\chi^v,\ell_\vnu(v)}] , & v\in V_3(\Ga),
\end{cases}
$$ 
where $\LL$ is a line bundle on $\xn{\Mbar}^v$ coming
from the restriction of the line bundle $\bL^v$ on $\MY$
(see Section \ref{sec:tan-obs}). 

We now give a more explicit description of $\LL$. 
Let 
$$
u:(X,\bq)\lra (\Po(m),p_0, p_{m})
$$ 
represent a point in $\Mbar^{\bu\sim}_{\chi^v,\nu^+,\nu^-}$, where
$\Po(m)$ is a chain of $m>0$ copies of $\Po$
with two relative divisors $p_0$ and $p_m$. 
Let $\Delta_l$ be the $l$-th irreducible component of $\Po(m)$ so
that $\Delta_l \cap \Delta_{l+1}=\{p_l\}$. 
The complex lines  
$$
\LL^0_u=T_{p_0}\Delta_1,\ \
\LL^1_u= \bigotimes_{l=1}^{m-1} T_{p_l}\Delta_l\otimes
T_{p_l}\Delta_{l+1},\ \
\LL^\infty_u= T_{p_m} \Delta_m.
$$ 
form line bundles $\LL^0$, $\LL^1$, $\LL^\infty$ on
$\Mbar^{\bu\sim}_{\chi^v,\nu^+,\nu^-}$
when we vary $u$ in $\Mbar^{\bu\sim}_{\chi^v,\nu^+,\nu^-}$.
The line bundle $\LL$ is given by
$$
\LL=\LL^0\otimes \LL^1\otimes \LL^\infty.
$$
Note that
$$
c_1(\LL^0)= -\psi^0,\ \ c_1(\LL^\infty) = -\psi^\infty,
$$
where $\psi^0, \psi^\infty$ are target $\psi$ classes (see e.g. \cite[Section 5]{LLZ2}).

Let $\cD$ be the divisor in $\Mbar^{\bu\sim}_{\chi^v,\nu^+,\nu^-}$
which corresponds to morphisms with target $\Po(m)$, $m>1$.
Then $\LL^1=\cO(\cD)$. Let 
\begin{equation}
\begin{aligned}
J'_{\chi^v,\nu^+,\nu^-}
 =  &\bigl\{ (\chi^+,\chi^-,\si)\mid
\chi^+,\chi^-\in 2\ZZ,\
\si\in \cP,\ 
|\si|=|\nu^+|=|\nu^-|,\\
&\quad  -\chi^+ + 2\ell(\si) -\chi^-=-\chi^v,
\quad  -\chi^\pm +\ell(\nu^\pm)+\ell(\si)>0\bigr\}.
\end{aligned}
\end{equation}
For each $(\chi^+,\chi^-,\si)\in J'_{\chi,\nu^+,\nu^-}$, there is a morphism
$$
\pi_{\chi^+,\chi^-,\si}:
\Mbar^{\bu\sim}_{\chi^+,\nu^+,\si} \times \Mbar^{\bu\sim}_{\chi^-,\si,\nu^-}
\lra \Mbar^{\bu\sim}_{\chi^v,\nu^+,\nu^-}
$$
with image contained in $\cD$. Moreover,
$$
[\Mbar^{\bu\sim}_{\chi^v,\nu^+,\nu^-}]\virt\cap c_1(\LL^1)
=\sum_{\tiny \begin{array}{c}(\chi^+,\chi^-,\si)\\ \in J'_{\chi^v,\nu^+,\nu^-}\end{array}}
\frac{a_\si}{|\Aut(\si)|} 
(\pi_{\chi^+,\chi^-,\si})_*\left( 
[\Mbar^{\bu\sim}_{\chi^+,\nu^+,\si}]\virt
\times[\Mbar^{\bu\sim}_{\chi^-,\si,\nu^-}]\virt\right)
$$
where $a_\si$ and $\Aut(\si)$ are defined in Section \ref{sec:W}.

\subsection{Contribution from Each Label}\label{sec:cy-each}
We use the definitions in Section \ref{sec:W} and Section \ref{sec:graph}.

In this subsection, we view the position $\fp(e)$ and 
the framing $\ff(e)$ as elements in
$$
\ZZ u_1\oplus \ZZ u_2=\Lambda_T\cong H^2_T(\mathrm{pt},\QQ).
$$
Recall that $H^*_T(\mathrm{pt};\QQ)=\QQ[u_1,u_2]$. The results
of localization calculations will involve rational functions of
$\fp(e)$ and $\ff(e)$.

Let $\xn{N}\virt$ denote the pull back of $\cT^{1,m}-\cT^{2,m}$
of $\xn{\cF}$ under $\xn{i}$. Let $r^v$
be defined as (\ref{eqn:rVone}) and (\ref{eqn:rVtwo}).
For $e\in E\ori(\Ga)$, let $\ee=\{e,-e\}\in E(\Ga)$ as before.

With the above notation and the explicit description of $[T^1\to T^2]$
in Section \ref{sec:tan-obs}, 
calculations similar to those in \cite[Appendix A]{LLZ1} show that
$$
\frac{1}{e_T(\xn{N}\virt)}
= \prod_{v\in V(\Ga)} B_v \prod_{\ee\in E(\Gamma)} B_\ee,
$$
where
$$
B_v=\begin{cases} 
1 &  v\in V_1(\Ga)\cup V_2(\Ga),\ r^v=0,\\
{\displaystyle (-1)^{\ell(\nu^e)-\chi^v/2} a_{\nu^e}
\frac{\ff(e)^{r^v}}{-\fp(e)-\psi^0} }
& v\in V_1(\Ga),\fv_1^{-1}(v)=\{e\},r^v>0, \\
{ \displaystyle (-1)^{\ell(\nu^e)-\chi^v/2}
\frac{ a_{\nu^e}a_{\nu^{e'}}\ff(e)^{r^v}}{(-\fp(e)-\psi^0)(-\fp(e')-\psi^\infty)} } 
& v\in V_2(\Ga),\fv_1^{-1}(v)=\{e,e'\},r^v>0,\\
{\displaystyle \prod_{e\in \fv_0^{-1}(v)}
\frac{a_{\nu^e }\Lambda^\vee(\fp(e))\fp(e)^{\ell_\vnu(v)-1}}
{\prod_{j=1}^{\ell(\nu^e)}
\left(\fp(e)(\fp(e)-\nu^e_j\psi^e_j)\right)} } & 
v\in V_3(\Gamma);
\end{cases}
$$

$$
B_\ee=
(-1)^{n^e d^{\ee}}\cdot
\begin{cases}
 E_{\nu^e}\bigl(\fp(e),\fl_0(e)\bigr)\cdot
E_{\nu^e}\bigl(\fp(-e)),\fl_0(-e)\bigr)
& \fv_0(e),\fv_1(e)\in V_3(\Ga),\\
(-1)^{\ell(\nu^e)-d^{\ee}} E_{\nu^e}(\fp(e),\fl_0(e))
& \fv_0(e)\in V_3(\Ga), \fv_1(e)\notin V_3(\Ga),\\
1 & \fv_0(e)\notin V_3(\Ga),\fv_1(e)\notin V_3(\Ga).
\end{cases}
$$
Recall that $n^e$ is defined in Definition \ref{def:vn}
and $E_\nu(x,y)$ is defined by (\ref{eqn:E}).

For $v\in V_2(\Ga)$, we have
\begin{eqnarray*}
&&\int_{[\xn{\Mbar}^v]\virt}\frac{\ff(e)^{r^v} }{(-\fp(e)-\psi^0)(-\fp(e')-\psi^\infty)}
=\int_{[\Mbar^{\bu\sim}_{\chi^v,\nu^e,\nu^{e'} } ]\virt}
\frac{\ff(e)^{r^v}c_1(\LL)}{(-\fp(e)-\psi^0)(\fp(e)-\psi^\infty)}\\
&=&\int_{[\Mbar^{\bu\sim}_{\chi^v,\nu^e,\nu^{e'} }]\virt}
\frac{\ff(e)^{r^v}(-\fp(e)-\psi^0
  +\fp(e)-\psi^\infty + c_1(\LL^1))}{(-\fp(e)-\psi^0)(\fp(e)-\psi^\infty)}\\
&=&\int_{[\Mbar^{\bu\sim}_{\chi^v,\nu^e,\nu^{e'} } ]\virt}
\frac{\ff(e)^{r^v}}{\fp(e)-\psi^\infty}
+\int_{[\Mbar^{\bu\sim}_{\chi^v,\nu^e,\nu^{e'} }]\virt}
\frac{\ff(e)^{r^v}}{-\fp(e)-\psi^0}\\
&& +\sum_{(\chi^+,\chi^-,\si)\in I_{\chi^v,\nu^e,\nu^{e'}}}
\frac{a_\si}{|\Aut(\si)|}
\int_{[\Mbar^{\bu\sim}_{\chi^+,\nu^e,\si}]\virt}
\frac{\ff(e)^{r^+_{\chi^+,\si}}}{-\fp(e)-\psi^0}
\int_{[\Mbar^{\bu\sim}_{\chi^-,\si,\nu^{e'}}]\virt}
\frac{\ff(e)^{r^-_{\chi^-,\si}}}{\fp(e)-\psi^\infty}\\
&=&|\Aut(\nu^e)\times\Aut(\nu^{e'})|\Bigl(\frac{\ff(e)}{\fp(e)}\Bigr)^{r^v}
\sum_{(\chi^+,\chi^-,\si)\in J_{\chi^v,\nu^e,\nu^{e'}}}(-1)^{r^+_{\chi^+,\si}}
  \frac{H^\bu_{\chi^+,\nu^e,\si}}{r^+_{\chi^+,\si}!}z_\si 
\frac{H^\bu_{\chi^-,\si,\nu^{e'}}}{r^-_{\chi^-,\si}!}
\end{eqnarray*}
where $r^+_{\chi^+,\si}=-\chi^+ +\ell(\nu^e)+\ell(\si)$,
$r^-_{\chi^-,\si}=-\chi^- +\ell(\si)+\ell(\nu^{e'})$, and
$$
J_{\chi^v,\nu^e,\nu^{e'}}=J'_{\chi^v,\nu^e,\nu^{e'}}
\cup\{(2\ell(\nu^e),\chi, \nu^e),(\chi,2\ell(\nu^{e'}),\nu^{e'})  \}.
$$

Given $v\in V_3(\Gamma)$, we have $\fv_0^{-1}(v)=\{ e_1,e_2,e_3\}$,
where $\fp(e_1)\wedge \fp(e_2)=u_1\wedge u_2$. Then $(e_1,e_2,e_3)$
is unique up to cyclic permutation.
Let
\begin{equation}\label{eqn:vnuw}
\vnu^v=(\nu^{e_1},\nu^{e_2},\nu^{e_3}),\quad
\bw^v=(\fp(e_1),\fp(e_2),\fp(e_3)).
\end{equation}
Then $V_{\chi^v,\vnu^v}(\bw^v)$ is independent of choice of cyclic
ordering of $e_1,e_2,e_3$, where $V_{\chi,\vnu}(\bw)$ is
defined by (\ref{eqn:V}). Set
$$
\xn{I}(u_1,u_2)=\int_{[\xn{\cF}]\virt}\frac{1}{ e_T(\xn{N}\virt)}.
$$
Then
\begin{eqnarray*}
&&\xn{I}(u_1,u_2)=\frac{1}{\prod_{\ee\in E(\Gamma)}|\xn{A}^\ee|}\int_{[\xn{\Mbar}]\virt}
\frac{1}{e_T(\xn{N}\virt)}\\
&=&|\Aut(\vmu)|\prod_{\ee\in E(\Gamma)}(-1)^{n^e d^{\ee} }z_{\nu^e}
\prod_{v\in V_3(\Gamma)} V_{\chi^v,\vnu^v}(\bw^v)
\prod_{e\in \fv_0^{-1}(V_3(\Ga))}E_{\nu^e}(\fp(e),\fl_1(e))\\
&& \cdot \prod_{v\in V_1(\Gamma), \fv_1(e)=v}
\sqrt{-1}^{\ell(\nu^e)+\ell(\mu^v)}(-1)^{d^{\ee}}
\cdot \Bigl(\sqrt{-1}\frac{\ff(e)}{\fp(e)}\Bigr)^{r^v}
\frac{H^\bu_{\chi^v,\nu^e,\mu^v} }{r^v!}\\
&& \cdot \prod_{\tiny \begin{array}{cc}v\in V_2(\Ga)\\ \fv_1^{-1}(v)=\{e,e'\}\end{array} }
\biggl(\sqrt{-1}^{\ell(\nu^e)+\ell(\nu^{e'})}
\Bigl(\sqrt{-1}\frac{\ff(e)}{\fp(e)}\Bigr)^{r^v}
\cdot \sum_{\tiny\begin{array}{l}(\chi^+,\chi^-,\si)\\ 
\in J_{\chi^v,\nu^e,\nu^{e'}}\end{array} }
   \frac{H^\bu_{\chi^+,\nu^e,\si} }{r^+_{\chi^+,\si}!} (-1)^{\ell(\si)} z_\si
   \frac{H^\bu_{\chi^-,\nu^{e'},\si}}{r^-_{\chi^-,\si}!}\biggr)
\end{eqnarray*}
So
\begin{equation}\label{eqn:Iuu}
\begin{aligned}
& \xn{I}(u_1,u_2)= |\Aut(\vmu)|
\prod_{\ee\in E(\Gamma)} (-1)^{n^e d^{\ee} }
\prod_{v\in V_3(\Gamma)}\sqrt{-1}^{\ell_\vnu(v)}
G^\bu_{\chi^v,\vnu^v}(\bw^v)\\
&\makebox[3cm]{ } \cdot \prod_{v\in V_1(\Gamma), \fv_1(e)=v }
\sqrt{-1}^{\ell(\mu^v)+\ell(\nu^v)}(-1)^{d^{\ee} }
\Bigl(\sqrt{-1}\frac{\ff(e)}{\fp(e)}\Bigr)^{r^v}
\frac{H^\bu_{\chi^v,\nu^v,\mu^i} }{r^v!}\\
& \cdot \prod_{\tiny\begin{array}{c} v\in V_2(\Ga)\\ \fv_1^{-1}(v)=\{e,e'\}\end{array} }
\biggl(\sqrt{-1}^{\ell(\nu^e)+\ell(\nu^{e'})}
\left(\sqrt{-1}\frac{\ff(e)}{\fp(e)}\right)^{r^v}
\cdot\sum_{\tiny\begin{array}{l}(\chi^+,\chi^-,\si)\\ 
\in J_{\chi^v,\nu^e,\nu^{e'}}\end{array} }
   \frac{H^\bu_{\chi^+,\nu^e,\si} }{r^+_{\chi^+,\si}!}z_\si (-1)^{\ell(\si)}
   \frac{H^\bu_{\chi^-,\nu^{e'},\si}}{r^-_{\chi^-,\si}!}\biggr)
\end{aligned}
\end{equation}

\subsection{Sum over Labels}\label{sec:cy-sum}
Finally, with the notation above, the formal relative
GW invariants of a general FTCY graph $\Ga$ are
$$
F^{\bu\Ga}_{\chi,\vd,\vmu}(u_1,u_2)=\frac{1}{|\Aut(\vmu)|}
\sum_{\pair\in\GX}\xn{I}(u_1,u_2).
$$
Define a generating function
\begin{equation}
F^{\bu\Ga}_{\vd,\vmu}(\lam;u_1,u_2)=
\sum_{\chi\in 2\ZZ,\chi\leq \ell(\vmu)} \lam^{-\chi+\ell(\vmu)}
F^{\bu\Ga}_{\chi,\vd,\vmu}(u_1,u_2).
\end{equation}
Then \eqref{eqn:Iuu} becomes
\begin{equation}\label{eqn:Fuu}
\begin{aligned}
& F^{\bu\Ga}_{\vd,\vmu}(\lam;u_1,u_2)
 =\sum_{|\nu^\ee|=\vd(\ee)}\prod_{\ee\in E(\Gamma)}(-1)^{n^e d^{\ee} }z_{\nu^\ee}
\prod_{v\in V_3(\Gamma)}\sqrt{-1}^{\ell(\vnu^v)}G^\bu_{\vnu^v}(\lam;\bw_v)\\
&\makebox[3cm]{ }\cdot \prod_{v\in V_1(\Gamma),\fv_1(e)=v}
(-1)^{d^{\ee} }\sqrt{-1}^{\ell(\nu^e)+\ell(\mu^v)}\Phi^\bu_{\nu^e,\mu^v}
\Bigl(\sqrt{-1}\frac{\ff(e)}{\fp(e)}\lam\Bigr)\\
& \cdot \prod_{v\in V_2(\Gamma),\fv_1^{-1}(v)=\{e,e'\} }
\sqrt{-1}^{\ell(\nu^e)+\ell(\nu^{e'})}
\Phi^\bu_{\nu^e,\si}\Bigl(\sqrt{-1}\frac{\ff(e)}{\fp(e)}\lam\Bigr)
(-1)^{\ell(\si)}z_\si
\Phi^\bu_{\nu^{e'},\si}\Bigl(\sqrt{-1}\frac{\ff(e')}{\fp(e')}\lam\Bigr)
\end{aligned}
\end{equation}
where $G^\bu_{\vmu}$ is defined by \eqref{eqn:Gmu},
$\bw_v$ is defined in \eqref{eqn:vnuw}, and
$\Phi^\bu_{\nu,\mu}$ is defined in Section \ref{sec:H}. 
Equations \eqref{eqn:FG} and \eqref{eqn:tF} imply
\begin{equation}\label{eqn:GK}
\begin{aligned}
& \sqrt{-1}^{\ell(\vmu)}G^\bu_{\vmu}(\lam;\fp(e_1),\fp(e_2),\fp(e_3))\\
=& \sqrt{-1}^{\ell(\vmu)}\sum_{|\nu^i|=|\mu^i|}\tF^\bu_{\vnu}(\lam;\mathrm{0})
     \prod_{i=1}^3 z_{\nu^i}\Phi^\bu_{\nu^i,\mu^i}
     \Bigl(\sqrt{-1}\frac{\fl_0(e_i)}{\fp(e_i)}\lam\Bigr)\\
=& (-1)^{\sum_{i=1}^3 \vd(\ee_i)}
    \sum_{|\nu^i|=|\mu^i|}F^\bu_{\vnu}(\lam;\mathrm{0})
    \prod_{i=1}^3 \sqrt{-1}^{\ell(\nu^i)-\ell(\mu^i)}z_{\nu^i}\Phi^\bu_{\nu^i,\mu^i}
    \Bigl(-\sqrt{-1}\frac{\fl_0(e_i)}{\fp(e_i)}\lam\Bigr)
\end{aligned}
\end{equation}

\subsection{Invariance}
In this subsection, we prove that formal relative Gromov-Witten
invariants are rational numbers independent of $u_1, u_2$ (Theorem
\ref{thm:R2}). We will use operations on FTCY graphs such as 
smoothing  and normalization (defined in Section \ref{sec:operation})
to reduce this to the invariance of the topological vertex
(Theorem \ref{thm:tri-inv}).
 
Let $\Ga$ be a FTCY graph, and let
$$
\Ga_2=\Gamma_{V_2(\Ga)},\ \ \Ga^2=\Ga^{V_2(\Ga)}.
$$
Then $\Gamma_2$, $\Gamma^2$ are {\em regular} FTCY
graphs. We call $\Gamma_2$ the {\em full smoothing}
of $\Gamma$, and $\Gamma^2$ the full resolution of
$\Gamma$. We have surjective maps
$$
\pi_2=\pi_{V_2(\Gamma)}:E\ori(\Gamma)\to E\ori(\Gamma_2),\ \
\pi^2=\pi^{V_2(\Gamma)}:V(\Gamma^2)\to V(\Gamma).
$$
\begin{defi}
Let $\Ga$ be a FTCY graph, and let $\Ga^2$ be the
full resolution of $\Ga$, 

Let $(\vd,\vmu)$ be an effective class of $\Ga$.
A {\em splitting type} of $(\vd,\vmu)$ is a map
$\vsi: V_2(\Ga)\to \cP$ such that $|\vsi(v)|=\vd(\ee)$ if
$\fv_1(e)=v$.

Given a splitting type $\vsi$ of an effective class $(\vd,\vmu)$
of $\Ga$, let $(\vd, \vmu\sqcup\vsi)$ denote
the effective class of $\Ga^2$ defined by 
$\vd:E(\Ga^2)=E(\Ga)\to \ZZ_{\geq 0}$ and
$$
\vmu\sqcup\vsi(v)=\left\{ \begin{array}{ll}
\vmu (\pi^2(v)), & \pi^2(v) \in V_1(\Ga)\\
\vsi (\pi^2(v)), & \pi^2(v) \in V_2(\Ga)
\end{array}\right.
$$
Let $S_{\vd,\vmu}$ denote the set of all splitting
types of $(\vd,\vmu)$.
\end{defi}

The following is clear from  the expression (\ref{eqn:Fuu}).
\begin{lemm}\label{thm:sing}
Let $\Ga$ be a FTCY graph, and let $(\vd,\vmu)$ be an effective
class of $\Ga$. Then
$$
F^{\bu\Ga}_{\vd,\vmu}(\lam;u_1,u_2)=\sum_{\si\in S_{\vd,\vmu}}z_{\vsi}
F^{\bu\Ga^2}_{\vd,\vmu\sqcup\vsi}(\lam;u_1,u_2)
$$
where $\displaystyle{z_{\vsi}=\prod_{v\in V_2(\Ga)} z_{\vsi(v)}}$.
\end{lemm} 

By Lemma \ref{thm:sing}, it suffices to consider regular FTCY graphs.
For a regular FTCY graph $\Ga$, (\ref{eqn:Fuu}) reduces to
\begin{equation}\label{eqn:Fu}
\begin{aligned}
& F^{\bu\Ga}_{\vd,\vmu}(\lam;u_1,u_2)
=\sum_{|\nu^\ee|=d^{\ee}}\prod_{\ee\in E(\Ga)}(-1)^{n^e d^{\ee} }z_{\nu^\ee}
\prod_{v\in V_3(\Gamma)}\sqrt{-1}^{\ell(\vnu^v)}G^\bu_{\vnu^v}(\lam;\bw_v)\\
&\makebox[2cm]{ }\cdot \prod_{v\in V_1(\Gamma),\fv_1(e)=v}
(-1)^{d^{\ee}}\sqrt{-1}^{\ell(\nu^e)+\ell(\mu^v)}\Phi^\bu_{\nu^e,\mu^v}
\left(\sqrt{-1}\frac{\ff(e)}{\fp(e)}\lam\right)
\end{aligned}
\end{equation}
since $V_2(\Ga)=\emptyset$.

Let $(\vd,\vmu)$ be the effective class of a regular FTCY graph.
Let $P_{\vd,\vmu}$ be the set of all maps $\vnu:E\ori(\Gamma)\to \cP$ such that
\begin{itemize}
\item $|\vnu(e)|=d^{\ee}$.
\item $\vnu(e)=\vmu(v)$ if $\fv_0(e)=v\in V_1(\Gamma)$.
\end{itemize}
Note that we do not require $\vnu(e)=\vnu(-e)$.
Denote $\vnu(e)$ by $\nu^e$.
Given $v\in V_3(\Gamma)$, there exist $e_1,e_2,e_3\in E(\Gamma)$, unique
up to a cyclic permutation, such that  $\fv_0^{-1}(v)=\{e_1,e_2,e_3\}$
and $\fp(e_1)\wedge\fp(e_2)=u_1\wedge u_2$. Define
\begin{equation}\label{eqn:vnuv}
\vnu^v=(\nu^{e_1},\nu^{e_2},\nu^{e_3}),\ \
z_{\vnu^v}=z_{\nu^{e_1}}z_{\nu^{e_2}} z_{\nu^{e_3}}.
\end{equation}
Note that $F^\bu_{\vnu^v}(\lam;\mathbf{0})$ and $z_{\vnu^v}$ are invariant
under cyclic permutation of $e_1,e_2,e_3$, thus well-defined.

Using (\ref{eqn:GK}) and the sum formula (\ref{eqn:Hsum}) of
double Hurwitz numbers, we can rewrite (\ref{eqn:Fu}) as follows:
\begin{equation}\label{eqn:FF}
F^{\bu\Ga}_{\vd,\vmu}(\lam;u_1,u_2)
=\sum_{\vnu\in P_{\vd,\vmu}}
\prod_{v\in V_3(V)} F^\bu_{\vnu^v}(\lam;\mathrm{0})z_{\vnu^v}
\prod_{\ee\in E(\Ga)}\sqrt{-1}^{\ell(\nu^e)-\ell(\nu^{-e})}
(-1)^{n^e d^{\ee}}\Phi^\bu_{\nu^e,\nu^{-e} }(\sqrt{-1}n^e\lam).
\end{equation}
Note that the right hand side of (\ref{eqn:FF}) does not depend
on $u_1,u_2$. This completes the proof of Theorem \ref{thm:R2}.
From now on, we write $F^{\bu\Ga}_{\vd,\vmu}(\lam)$
instead of $F^{\bu\Ga}_{\vd,\vmu}(\lam;u_1,u_2)$. We define
$$
F^{\bu\Ga}_{\chi,\vd,\vmu}=F^{\bu\Ga}_{\chi,\vd,\vmu}(u_1,u_2),
$$
to be {\em formal relative Gromov-Witten invariants} of $\hat{Y}\urel_\Gamma$.

\subsection{Gluing Formulae} \label{sec:cy-glue}

Let $(\vd,\vmu)$ be an effective class of a regular FTCY graph $\Ga$.
Let 
$$
T_{\vd,\vmu}=\left\{ \vnu:E(\Gamma)\to \cP\ \Bigl |\ 
|\vnu(e)|=\vd(\ee),\ \vnu(-e)=\vnu(e)^t \right\}.
$$
Note that we do not require $\vnu(e)=\vmu(v)$ if $\fv_0(e)=v\in V_1(E)$.
We have
\begin{equation}\label{eqn:KtC}
F^\bu_{\vmu}(\lam;\mathbf{0})=\frac{(-1)^{|\mu^1|+|\mu^2|+|\mu^3|} }{\sqrt{-1}^{\ell(\vmu)}}
\sum_{|\nu^i|=|\mu^i|}\tC_{\vnu}(\lam)
\prod_{i=1}^3\frac{\chi_{\nu^i}(\mu^i) }{z_{\mu^i}},
\end{equation}
where $\tC_{\vnu}(\lam)=\tC_{\vnu}(\lam;\mathbf{0})$. Applying
(\ref{eqn:KtC}) and  the Burnside formula (\ref{eqn:burnside}) of
double Hurwitz numbers, we see that (\ref{eqn:FF}) is equivalent
to the following.
\begin{prop}\label{thm:R4}
Let $\Ga$ be a regular FTCY graph. Then
$$
F^{\bu\Ga}_{\vd,\vmu}(\lam)=\sum_{\vnu\in T_{\vd,\vmu}}
\prod_{\ee\in E(\Gamma)}(-1)^{(n^e+1) d^{\ee} }
e^{-\sqrt{-1}\kappa_{\nu^e} n^e\lam/2}
\prod_{v\in V_3(\Ga)}\tC_{\vnu^v}(\lam)
\prod_{\tiny\begin{array}{c}v\in V_1(\Ga)\\ v_0(e)=v\end{array}}
\frac{\chi_{\nu^e}(\mu^v)}{\sqrt{-1}^{\ell(\mu^v)}z_{\mu^v} }.
$$
\end{prop}
Recall that $\kappa_{\nu^e}$ is defined by \eqref{eqn:kappa}, and
we have $n^{-e}=-n^e$, $\kappa_{\nu^t}=-\kappa_\nu$, so  
$$
\kappa_{\nu^{-e} }n^{-e}=\kappa_{(\nu^e)^t}\cdot (-n^e)=\kappa_{\nu^e} n^e.
$$

\begin{theo}[gluing formula]\label{thm:R3}
Let $\Ga$ be a FTCY graph, and let $\Ga_2$ and $\Gamma^2$
be its full smoothing and its full resolution, respectively.
Let $(\vd,\vmu)$ be an effective class of $\Ga$ which
can also be viewed as an effective class of $\Ga_2$. Then
\begin{equation}
F^{\bu\Gamma_2}_{\vd,\vmu}(\lam)=F^{\bu\Ga}_{\vd,\vmu}(\lam)
=\sum_{\vsi\in S_{\vd,\vmu} }z_{\vsi}
F^{\bu\Gamma^2}_{\vd_\Gamma,\vmu\sqcup\vsi}(\lam).
\end{equation}
\end{theo}

\begin{proof}
By Lemma \ref{thm:sing} and Proposition \ref{thm:R4}, it suffices to 
show that if $|\mu|=|\nu|=d$, then
$$
\sum_{|\si|=d}\frac{\chi_\mu(\si)}{\sqrt{-1}^{\ell(\si)}z_\si}
z_\si\frac{\chi_\nu(\si)}{\sqrt{-1}^{\ell(\si)}z_\si}
=(-1)^d\delta_{\mu(\nu^t)},
$$
which is obvious.
\end{proof}

\subsection{Sum over Effective Classes}
Given a regular FTCY graph, let $\mathrm{Eff}(\Gamma)$ denote
the set of effective classes of $\Ga$. Introduce formal K\"{a}hler
parameters
$$
\bt=\{ t^\ee: \ee\in E(\Gamma) \}
$$
and winding parameters
$$
\bp=\{ p^v=(p^v_1,p^v_2,\ldots):v\in V_1(\Gamma)\}
$$
We define the {\em formal relative Gromov-Witten partition function} of $\hat{Y}\urel_\Gamma$
to be
\begin{equation}
Z^\Ga_{\mathrm{rel}}(\lam;\bt;\bp)=\sum_{(\vd,\vmu)\in \mathrm{Eff}(\Gamma)}
F^{\bu\Gamma}_{\vd,\vmu}(\lam)
e^{-\sum_{\ee\in E(\Gamma)}\vd(\ee)t^\ee}
\prod_{v\in V_1(\Gamma)}p^v_{\mu^v}
\end{equation}
where $p^v_\mu=p^v_{\mu_1}\cdots p^v_{\mu_{\ell(\mu)}}$.

Let $T^\Ga$ denote the set of pairs
$(\vnu,\vmu)$ such that
\begin{itemize}
\item $\vnu:E\ori(\Ga)\to \cP$ such that $\vnu(-e)=\vnu(e)^t$. 
\item $\vmu:V_1(\Gamma)\to \cP$. 
\item $|\nu^e|=|\mu^v|$ if $\fv_0(e)=v$.
\end{itemize}
We abbreviate $\vnu(e)$ to $\nu^e$ for $e\in E\ori(\Ga)$, abbreviate
$\vmu(v)$ to $\mu^v$ for $v\in V_1(\Ga)$, and
define $\vmu^v$ by (\ref{eqn:vnuv}) for $v\in V_3(\Ga)$. The
following is a direct consequence of Proposition \ref{thm:R4}.
\begin{coro}\label{thm:ZtC}
\begin{eqnarray*}
Z^\Ga_{\mathrm{rel}}(\lam;\bt;\bp)&=&
\sum_{(\vnu,\vmu)\in T^\Gamma}
\prod_{\ee\in E(\Gamma)}e^{-|\nu^e|t^\ee}
(-1)^{(n^e+1)|\nu^e|} e^{-\sqrt{-1}\kappa_{\nu^e} n^e\lam/2}\\
&&\cdot\prod_{v\in V_3(\Gamma)}\tC_{\vnu^v}(\lam)
\prod_{v\in V_1(\Ga), v_0(e)=v}
\frac{\chi_{\nu^e}(\mu^v)}{\sqrt{-1}^{\ell(\mu^v)}z_{\mu^v} }
\end{eqnarray*}
\end{coro}

%% file: sec8.tex
\section{Combinatorial Expressions for the topological vertex}\label{sec:comb}
We use the notation introduced in Section \ref{sec:W}. The goal of this
section is to derive the following combinatorial expression
for $\tC_{\vmu}(\lam)$:
\begin{theo}\label{thm:R5}
Let $\vmu\in\cP^3_+$. Then
$$
\tC_{\vmu}(\lam)=\tilde{\cW}_\vmu(q),
$$
where $q=e^{\sqrt{-1}\lam}$, and $\tilde{\cW}_\vmu(q)$ is defined by 
\eqref{eqn:tW}.
\end{theo}

We now outline our strategy to prove Theorem \ref{thm:R5}. 
By Proposition \ref{thm:GtC},
$$
\tC_\vmu(\lam)=\sum_{|\nu^i|=|\mu^i|}\prod_{i=1}^3\chi_{\mu^i}(\nu^i)
q^{-\frac{1}{2}\left(\sum_{i=1}^3\kappa_{\nu^i}\frac{w_{i+1}}{w_i}\right)}G^\bu_\vnu(\lam;\bw),
$$
where $\bw$ is as in \eqref{eqn:w-convention}. Since the above sum is independent
of $\bw$, we may take $\bw=(1,1,-2)$ and obtain
\begin{equation}
\tC_\vmu(\lam)=\sum_{|\nu^i|=|\mu^i|}\prod_{i=1}^3\chi_{\mu^i}(\nu^i)
q^{-\frac{1}{2}\kappa_{\nu^1}+\kappa_{\nu^2}+\frac{1}{4}\kappa_{\nu^3}}\cdot G^\bu_\vnu(\lam;1,1,-2).
\end{equation}

In Section \ref{sec:Gonetwo}, we will show that the main result in \cite{LLZ2}
gives a combinatorial expression of $G^\bu_{\mu,\nu,\emptyset}(\lam;\bw)$
(Theorem \ref{thm:Gtwo}).
In Section \ref{sec:reduction}, we will relate $G^\bu_\vmu(\lam;1,1,-2)$ to
$G^\bu_{\emptyset,\mu^1\cup\mu^2,\mu^3}(\lam;1,1,-2)$.
This gives the combinatorial expression $\tilde{\cW}_\vmu(q)$ in Theorem \ref{thm:R5}.
Moreover, (\ref{eqn:GtC}) and Theorem \ref{thm:R5} imply the following formula
of three-partition Hodge integrals.
\begin{theo}[Formula of three-partition Hodge integrals]\label{thm:Gthree}
Let $\bw$ be as in \eqref{eqn:w-convention} and let $\vmu=(\up{\mu})\in\cP^3_+$.
Then
\begin{equation}
G^\bu_\vmu(\lam;\bw)=\sum_{|\nu^i|=|\mu^i|}\prod_{i=1}^3
\frac{\chi_{\nu^i}(\mu^i)}{z_{\mu^i}}
q^{\frac{1}{2}\left(\sum_{i=1}^3\kappa_{\nu^i}\frac{w_{i+1}}{w_i}\right)}
\tilde{\cW}_\vnu(q).
\end{equation}
\end{theo}

The cyclic symmetry of $\tC_{\vmu}(\lam)$ is obvious from definition. By
Theorem \ref{thm:R5} we have the following cyclic symmetry
$$
\tilde{\cW}_{\mu^1,\mu^2,\mu^3}(q)=\tilde{\cW}_{\mu^2,\mu^3,\mu^1}(q)
=\tilde{\cW}_{\mu^3,\mu^1,\mu^2}(q)
$$
which is far from being obvious.

Finally, we conjecture that the combinatorial expression
$\tilde{\cW}_{\vmu}(q)$ coincides with $\cW_\vmu(q)$ predicted
in \cite{AKMV}:
\begin{conj}\label{conj:WW}
Let $\vmu\in\cP^3_+$. Then
$$
\tilde{\cW}_{\vmu}(q)=\cW_\vmu(q),
$$
where $q=e^{\sqrt{-1}\lam}$, and $\cW_\vmu(q)$ is defined by \eqref{eqn:Wthree}.
\end{conj}

We have strong evidence for Conjecture \ref{conj:WW}.   
By Theorem \ref{thm:R5} and Corollary \ref{thm:R5-two}, Conjecture
\ref{conj:WW} holds when one of the three partitions is empty. When none of the partitions
is empty, A. Klemm has checked by computer that Conjecture \ref{conj:WW} holds
in all the cases where
$$
|\mu^i|\leq 6,\ \ i=1,2,3.
$$
We will list some of these cases in Section \ref{sec:lowdegree}.

\subsection{One-Partition and Two-Partition Hodge Integrals}
\label{sec:Gonetwo}
We recall some notation in \cite{LLZ1}.
\begin{equation}\label{eqn:CG}
\cC^\bu_\mu(\lam;\tau)=\sqrt{-1}^{|\mu|}G^\bu_{\mu,\emptyset,\emptyset}(\lam;1,\tau,-\tau-1).
\end{equation}
\begin{equation}\label{eqn:VW}
V_\mu(q)=q^{-\kappa_\mu/4}\sqrt{-1}^{|\mu|}\cW_\mu(q).
\end{equation}
where $\cW_\mu(q)=\cW_{\mu,\emptyset,\emptyset}(q)$ is defined in Section \ref{sec:W}.
The main result of \cite{LLZ1} is the following formula conjectured by
Mari\~{n}o and Vafa \cite{Mar-Vaf} (see \cite{Oko-Pan2} for another proof):
\begin{theo}\label{thm:MV}
\begin{equation}
\cC^\bu_\mu (\lam;\tau)=\sum_{|\nu|=|\mu|}\frac{\chi(\mu)}{z_\mu}
q^{\kappa_\nu(\tau+\frac{1}{2})/2}V_\nu(q)
\end{equation}
\end{theo}

Theorem \ref{thm:MV} can be reformulated in our notation as follows:
\begin{theo}[Formula of one-partition Hodge integrals]\label{thm:Gone}
Let $\bw$ be as in \eqref{eqn:w-convention}, and let $\mu\in\cP_+$.
Then
\begin{equation}\label{eqn:Gone}
G^\bu_{\mu,\emptyset,\emptyset}(\lam;\bw)
=\sum_{|\nu|=|\mu|}\frac{\chi_\nu(\mu)}{z_\mu}
q^{\frac{1}{2}\kappa_\nu\frac{w_2}{w_1}}\cW_{\nu,\emptyset,\emptyset}(q).
\end{equation}
\end{theo}

Let
\begin{equation}
G^\bu_{\mu^+,\mu^-}(\lam;\tau)=(-1)^{|\mu^-|-\ell(\mu^-)}
G^\bu_{\mu^+,\mu^-}(\lam;1,\tau,-1-\tau).
\end{equation}
The main result of \cite{LLZ2}
is the following formula conjectured in \cite{Zho2}:
\begin{theo}\label{thm:Z}
Let $(\mu^+,\mu^-)\in\cP^2_+$. Then
$$
G^\bu_{\mu^+,\mu^-}(\lam;\tau)=\sum_{|\nu^\pm|=|\mu^\pm|}
\frac{\chi_{\nu^+}(\mu^+)}{z_{\mu^+}}
\frac{\chi_{\nu^-}(\mu^-)}{z_{\mu^-}}
q^{\left(\kappa_{\nu^+}\tau+\kappa_{\nu^-}\tau^{-1}\right)/2}
\cW_{\nu^+,\nu^-}(q).$$
\end{theo}

We now reformulate Theorem \ref{thm:Z} in the notation of this paper.
\begin{eqnarray*}
&&G^\bu_{\mu^1,\mu^2,\emptyset}(\lam;1,\tau,-1-\tau)\\
&=&(-1)^{|\mu^2|-\ell(\mu^2)}
\sum_{|\nu^i|=|\mu^i|}
\frac{\chi_{\nu^1}(\mu^1)}{z_{\mu^1}}
\frac{\chi_{\nu^2}(\mu^2)}{z_{\mu^2}}
q^{\left(\kappa_{\nu^1}\tau+\kappa_{\nu^2}\tau^{-1}\right)/2}
\cW_{\nu^1,\nu^2}(q)\\
&=&\sum_{|\nu^i|=|\mu^i|}
\frac{\chi_{\nu^1}(\mu^1)}{z_{\mu^1}}
\frac{\chi_{(\nu^2)^t}(\mu^2)}{z_{\mu^2}}
q^{\left(\kappa_{\nu^1}\tau+\kappa_{\nu^2}\tau^{-1}\right)/2}
q^{\kappa_{\nu^2}/2}\cW_{\nu^1,(\nu^2)^t,\emptyset}(q)\\
&=&\sum_{|\nu^i|=|\mu^i|}
\frac{\chi_{\nu^1}(\mu^1)}{z_{\mu^1}}
\frac{\chi_{\nu^2}(\mu^2)}{z_{\mu^2}}
q^{\left(\kappa_{\nu^1}\tau+\kappa_{\nu^2}\frac{-\tau-1}{\tau}\right)/2}
\cW_{\nu^1,\nu^2,\emptyset}(q)
\end{eqnarray*}

Theorem \ref{thm:Z} is equivalent to the following:
\begin{theo}[Formula of two-partition Hodge integrals]\label{thm:Gtwo}
Let $\bw$ be as in \eqref{eqn:w-convention} and let $(\mu^1,\mu^2)\in\cP^2_+$. Then
\begin{equation}
G^\bu_{\mu^1,\mu^2,\emptyset}(\lam;\bw)=\sum_{|\nu^i|=|\mu^i|}
\sum_{|\nu^i|=|\mu^i|}
\frac{\chi_{\nu^1}(\mu^1)}{z_{\mu^1}}
\frac{\chi_{\nu^2}(\mu^2)}{z_{\mu^2}}
q^{\frac{1}{2}\left(\kappa_{\nu^1}\frac{w_2}{w_1}+\kappa_{\nu^2}\frac{w_3}{w_2}\right)}
\cW_{\nu^1,\nu^2,\emptyset}(q).
\end{equation}
\end{theo}
Note that Theorem \ref{thm:Gone} corresponds the special case
where $(\mu^1,\mu^2)=(\mu,\emptyset)$.
Theorem \ref{thm:Gtwo} and (\ref{eqn:tCG}) imply
\begin{coro}\label{thm:R5-two}
Let $\vmu=(\up{\mu})\in\cP^3_+$, and let $q=e^{\sqrt{-1}\lam}$. Then
$$
\tC_{\vmu}(\lam)=\cW_\vmu(q)
$$
when one of $\mu^1,\mu^2,\mu^3$ is empty.
\end{coro}

\subsection{Reduction}
\label{sec:reduction}

Recall that
$$
G_{g,\vmu}(\tau)=G_{g,\vmu}(1,\tau,-\tau-1).
$$
For two partitions $\mu^1$ and $\mu^2$, let $\mu^1\cup \mu^2$
be the partition with 
$$
m_i(\mu^1\cup \mu^2)=m_i(\mu^1)+m_i(\mu^2),\quad
\forall i\geq 1. 
$$
We have
\begin{lemm}\label{thm:Gg-reduce}
Let $\vmu=(\up{\mu})\in\cP^3_+$. Then
\begin{equation}\label{eqn:Gg-reduce}
\begin{aligned}
G_{g,\vmu}(\lambda;1)
=& (-1)^{|\mu^1|-\ell(\mu^1)}
\frac{z_{\mu^1\cup \mu^2}}{z_{\mu^1}\cdot z_{\mu^2}}
G_{g,\emptyset,\mu^1 \cup \mu^2, \mu^3}(\lambda;1)\\
& +\delta_{g0}\sum_{m\geq 1} \delta_{\mu^1 (m)} \delta_{\mu^2\emptyset}
\delta_{\mu^3(2m)}\frac{(-1)^{m-1}}{m}
\end{aligned}
\end{equation}
\end{lemm}
\begin{proof}
Let
$$
I_{g,\vmu}(\bw)=\int_{\Mbar_{g,\ell(\vmu)} }
\prod_{i=1}^3\frac{\Lambda^\vee(w_i) w_i^{\ell(\vmu)-1}}
{\prod_{j=1}^{\ell(\mu^i)}(w_i(w_i-\mu^i_j\psi_{d_\vmu^i+j}) }
$$
and let $I_{g,\vmu}(\tau)=I_{g,\vmu}(1,\tau,-\tau-1)$.
Then
\begin{equation}
I_{0,\vmu}(\tau)=\frac{(\tau(-\tau-1))^{\ell(\vmu)-1}}{\tau^{2\ell(\mu^2)}(-\tau-1)^{2\ell(\mu^3)}}
\left(|\mu^1|+\frac{|\mu^2|}{\tau} +\frac{|\mu^3|}{-\tau-1} \right)^{\ell(\vmu)-3}
\end{equation}
Note that $I_{g,\vmu}(\tau)$ has a pole at $\tau=1$ only if
\begin{equation}\label{eqn:except}
g=0,\ \ \vmu=((m),\emptyset,(2m))\text{ or } (\emptyset,(m),(2m)),
\end{equation}
where $m>0$. Let
\begin{equation}
E_\mu(\tau)=\prod_{j=1}^{\ell(\mu)}\frac{\prod_{a=1}^{\mu_j-1}(\tau \mu_j+a)}{(\mu_j-1)!}.
\end{equation}
Then $E_\mu(\tau)$ is a polynomial in $\tau$ of degree $|\mu|-\ell(\mu)$, and
$$
E_\mu(-\tau-1)=(-1)^{|\mu|-\ell(\mu)} E_\mu(\tau).
$$
Then
\begin{eqnarray*}
&& G_{g,\vmu}(\tau)
=\frac{(-\sqrt{-1})^{\ell(\vmu)} }{|\Aut(\vmu)|}
E_{\mu^1}(\tau)E_{\mu^2}(-1-\tau^{-1})E_{\mu^3}\left(\frac{1}{-\tau-1}\right)
I_{g,\vmu}(\tau)\\
&&\quad = (-1)^{|\mu^1|-\ell(\mu^1)}\frac{(-\sqrt{-1})^{\ell(\vmu)} }{|\Aut(\vmu)|}
E_{\mu^1}(-1-\tau)
E_{\mu^2}(-1-\tau^{-1})E_{\mu^3}\left(\frac{1}{-\tau-1}\right)
I_{g,\vmu}(\tau)
\end{eqnarray*}
while
\begin{eqnarray*}
G_{g,\emptyset,\mu^1\cup\mu^2,\mu^3}(\tau)
&=&\frac{(-\sqrt{-1})^{\ell(\vmu)} }{|\Aut(\mu^1\cup\mu^2)\times \Aut(\mu^3)|}
E_{\mu^1\cup\mu^2}(-1-\tau^{-1})\\
&&\cdot E_{\mu^3}\left(\frac{1}{-\tau-1}\right)
I_{g,\emptyset,\mu^1\cup\mu^2,\mu^3}(\tau)
\end{eqnarray*}
where $E_{\mu^1\cup\mu^2}(-1-\tau^{-1})
=E_{\mu^1}(-1-\tau^{-1}) E_{\mu^2}(-1-\tau^{-1})$.

Suppose that $(g,\vmu)$ is not the exceptional case listed in
(\ref{eqn:except}). Then neither is $(g,\emptyset,\mu^1\cup\mu^2,\mu^3)$.
It is immediate from the definition that
$$
I_{g,\up{\mu}}(1)=I_{g,\emptyset,\mu^1\cup\mu^2,\mu^3}(1),
$$
so
\begin{equation}\label{eqn:general}
G_{g,\vmu}(1)=(-1)^{|\mu^1|-\ell(\mu^1)}
\frac{|\Aut(\mu^1\cup\mu^2)|}{|\Aut(\mu^1)\times \Aut(\mu^2)|}
G_{g,\emptyset,\mu^1\cup\mu^2,\mu^3}(1)
\end{equation}

For the exceptional case (\ref{eqn:except}), we have
\begin{eqnarray*}
G_{0,(m),\emptyset,(2m)}(\tau)
&=&\frac{\tau}{(\tau+1)(m-1)!(2m-1)!}
   \prod_{a=1}^{m-1}(\tau m+a)\\
&& \cdot \prod_{a=1}^{m-1}(\frac{2m}{-\tau-1}+a)
\prod_{a=m+1}^{2m-1}(\frac{2m}{-\tau-1}+a) 
\end{eqnarray*}
while
\begin{eqnarray*}
G_{0,\emptyset,(m),(2m)}(\tau)
&=&\frac{-1}{(\tau+1)(m-1)!(2m-1)!}
\prod_{a=1}^{m-1}(\frac{-\tau-1}{\tau} m+a)\\
&& \cdot \prod_{a=1}^{m-1}(\frac{2m}{-\tau-1}+a)
\prod_{a=m+1}^{2m-1}(\frac{2m}{-\tau-1}+a)
\end{eqnarray*}
So
\begin{equation}\label{eqn:error}
G_{0,(m),\emptyset,(2m)}(1)=\frac{(-1)^{m-1}}{2m},\quad
G_{0,\emptyset,(m),(2m)}(1)=\frac{-1}{2m}.
\end{equation}
Combining the general case (\ref{eqn:general}) and
the exceptional case (\ref{eqn:error}), we obtain (\ref{eqn:Gg-reduce}).
\end{proof}

Let $\bp$, $p^i$, $p^i_\mu$ be defined as in Section \ref{sec:G},
and let $G^\bu(\lam;\bp;\tau)$ be defined as in \eqref{eqn:Gmu}. We have
\begin{lemm}\label{thm:main-extra}
Let
\begin{equation}\label{eqn:plus}
p^+_i =(-1)^{i-1} p_i^1 +  p_i^2,\quad
p^+_\mu=\prod_{j=1}^{\ell(\mu)}p^+_{\mu_j}. 
\end{equation}
Then
\begin{equation}\label{eqn:Gplus}
G^\bu(\lam;p^1,p^2,p^3;1)=G^\bu(\lam;0,p^+, p^3;1)\exp
\biggl(\, \sum_{m\geq 1} \frac{(-1)^{m-1}}{m}p^1_m p^3_{2m}\biggr).
\end{equation}
\end{lemm}
\begin{proof}
We have
$$
G(\lam;\bp;1)= G(\lam;\bp;1,1,-2)
=\sum_{\vmu\in\cP^3_+}\sum_{g=0}^\infty\lam^{2g-2+\ell(\vmu)}
   G_{g,\vmu}(\lam;1)p^1_{\mu^1}p^2_{\mu^2}p^3_{\mu^3}.
$$
By Lemma \ref{thm:Gg-reduce},
\begin{eqnarray*}
&& G(\lam;\bp;1)\\
&=&\sum_{\vmu\in\cP^3_+}\sum_{g=0}^\infty\lam^{2g-2+\ell(\vmu)}
   G_{g,\emptyset,\mu^1\cup\mu^2,\mu^3}(1)
  \frac{z_{\mu^1\cup\mu^2} }{z_{\mu^1}z_{\mu^2}}
  (-1)^{|\mu^1|-\ell(\mu^1)}p^1_{\mu^1}p^2_{\mu^2}p^3_{\mu^3}
+\sum_{m\geq 1} \frac{(-1)^{m-1}}{m}p^1_m p^3_{2m}\\
&=&\sum_{(\mu^+,\mu^3)\in\cP^2_+}
\sum_{g=0}^\infty\lam^{2g-2+\ell(\mu^+)+\ell(\mu^3)}
   G_{g,\emptyset,\mu^+,\mu^3}(1)
  \left(\sum_{\mu^1\cup\mu^2=\mu^+}\frac{z_{\mu^+} }{z_{\mu^1}z_{\mu^2}}
  (-1)^{|\mu^1|-\ell(\mu^1)}p^1_{\mu^1}p^2_{\mu^2}\right)p^3_{\mu^3}\\
&&+\sum_{m\geq 1}  \frac{(-1)^{m-1}}{m}p^1_m p^3_{2m}
\end{eqnarray*}
It is easy to see that
\begin{equation} \label{eqn:Union}
\sum_{\mu^1\cup\mu^2=\mu^+}\frac{z_{\mu^+}}{z_{\mu^1} z_{\mu^2}}
(-1)^{|\mu^1|-\ell(\mu^1)}p^1_{\mu^1}p^2_{\mu^2} = p_{\mu^+}^+.
\end{equation}
So
\begin{equation}
G(\lam;\up{p};1)=
G(\lam;0,p^+, p^3;1)
+\sum_{m=1}^\infty \frac{(-1)^{m-1}}{m}p^1_m p^3_{2m}
\end{equation}
which is equivalent to \eqref{eqn:Gplus}.
\end{proof}

\subsection{Combinatorial Expression}
\begin{lemm}\label{thm:main}
Let $p^+$ be defined by \eqref{eqn:plus}.
\begin{equation}
G^\bu(\lambda;0,p^+,p^3;1)=
 \sum_{\nu^{+}, \nu^i, \mu^i \in \cP}
c^{\nu^+}_{(\nu^1)^t\nu^2} q^{(-2\kappa_{\nu^+}  - \frac{ \kappa_{\nu^3} }{2})/2}
\cW_{\nu^+, \nu^3}(q)
(-1)^{|\mu^3|-\ell(\mu^3)}\prod_{i=1}^3
\frac{\chi_{\nu^i}(\mu^i)}{z_{\mu^i}}p^i_{\mu^i}.
\end{equation}
\end{lemm}

\begin{proof}
By Theorem \ref{thm:Gtwo},
$$
G^\bu(\lambda;0,p^+,p^3;1)
=  \sum_{\mu^{\pm}, \nu^{\pm}, \mu^3\in \cP}
\frac{\chi_{\nu^+}(\mu^+)}{z_{\mu^+}} \frac{\chi_{\nu^3}(\mu^3)}{z_{\mu^3}}
q^{(-2\kappa_{\nu^+} -\frac{\kappa_{\nu^3}}{2})/2}
\cW_{\emptyset,\nu^+, \nu^3}(q) p^+_{\mu^+}p^3_{\mu^3}.
$$

Recall that
$$
\cW_{\emptyset,\nu^+,\nu^3}(q)
=q^{\kappa_{\nu^3}/2}\cW_{\nu^+,(\nu^3)^t}(q),
$$
$$
p^+_{\mu^+}=\sum_{\mu^1\cup\mu^2=\mu^+}\frac{z_{\mu^+}}{z_{\mu^2}z_{\mu^2}}
(-1)^{|\mu^1|-\ell(\mu^1)} p^1_{\mu^1} p^2_{\mu^2}.
$$
Let ${\displaystyle
s^i_\mu =\sum_{|\nu|=|\mu|} \frac{\chi_\mu(\nu)}{z_\nu} p_\nu }$
be Schur functions. Then
\begin{eqnarray*}
&& G^\bu(\lambda;0,p^+,p^3;1)\\
&=&  \sum_{\mu^{\pm}, \nu^{\pm}, \mu^3\in \cP}
\frac{\chi_{\nu^+}(\mu^+)}{z_{\mu^+}} \frac{\chi_{\nu^3}(\mu^3)}{z_{\mu^3}}
q^{(-2\kappa_{\nu^+}  + \frac{\kappa_{\nu^3}}{2})/2 }
\cW_{\nu^+, (\nu^3)^t}(q) p^+_{\mu^+}p^3_{\mu^3}\\
& = & \sum_{\mu^{\pm}, \nu^{\pm}, \mu^3\in \cP}
\frac{\chi_{\nu^+}(\mu^+)}{z_{\mu^+}} \frac{\chi_{(\nu^3)^t}(\mu^3)}{z_{\mu^3}}
q^{(-2\kappa_{\nu^+}  - \frac{\kappa_{\nu^3}}{2})/2}
\cW_{\nu^+, \nu^3}(q) p^3_{\mu^3} \\
&& \cdot \sum_{\mu^1 \cup \mu^2=\mu^+} \frac{z_{\mu^+}}{z_{\mu^1} \cdot z_{\mu^2}}
(-1)^{|\mu^1|-\ell(\mu^1)}p^1_{\mu^1}p^2_{\mu^2}\\
& = & \sum_{\mu^{i},\nu^{+}, \nu^3 \in \cP}
\frac{\chi_{\nu^+}(\mu^1\cup \mu^2)}{z_{\mu^1} \cdot z_{\mu^2}}
\frac{\chi_{(\nu^3)^t}(\mu^3)}{z_{\mu^3}}
q^{(-2\kappa_{\nu^+}  - \frac{\kappa_{\nu^3}}{2})/2} \cW_{\nu^+, \nu^3}(q)
(-1)^{|\mu^1|-\ell(\mu^1)} p^1_{\mu^1} p^2_{\mu^2}  p^3_{\mu^3}\\
& = & \sum_{\mu^{1},\mu^2,\nu^{+}, \nu^i \in \cP}
\frac{\chi_{\nu^+}(\mu^1\cup \mu^2)\chi_{\nu^1}(\mu^1)\chi_{\nu^2}(\mu^2)}{z_{\mu^1} \cdot z_{\mu^2}}
q^{(-2\kappa_{\nu^+}  - \frac{\kappa_{\nu^3}}{2})/2} \cW_{\nu^+, \nu^3}(q)
(-1)^{|\mu^1|-\ell(\mu^1)} s^1_{\nu^1} s^2_{\nu^2}  s^3_{(\nu^3)^t} \\
& = & \sum_{\mu^{1},\mu^2,\nu^{+}, \nu^i \in \cP}
\frac{\chi_{\nu^+}(\mu^1\cup \mu^2)\chi_{(\nu^1)^t}(\mu^1)\chi_{\nu^2}(\mu^2)}{z_{\mu^1} \cdot z_{\mu^2}}
q^{(-2\kappa_{\nu^+}  - \frac{\kappa_{\nu^3}}{2})/2} \cW_{\nu^+, \nu^3}(q)
s^1_{\nu^1} s^2_{\nu^2}  s^3_{(\nu^3)^t}\\
& = & \sum_{\nu^{+}, \nu^i \in \cP}
c^{\nu^+}_{(\nu^1)^t\nu^2} q^{(-2\kappa_{\nu^+}  - \frac{\kappa_{\nu^3}}{2})/2} \cW_{\nu^+, \nu^3}(q)
s^1_{\nu^1} s^2_{\nu^2}  s^3_{(\nu^3)^t} \\
& = & \sum_{\nu^{+}, \nu^i, \mu^i \in \cP}
c^{\nu^+}_{(\nu^1)^t\nu^2} q^{(-2\kappa_{\nu^+}  - \frac{\kappa_{\nu^3}}{2})/2} \cW_{\nu^+, \nu^3}(q)
(-1)^{|\mu^3|-\ell(\mu^3)}\prod_{i=1}^3
\frac{\chi_{\nu^i}(\mu^i)}{z_{\mu^i}}p^i_{\mu^i}.
\end{eqnarray*}
In the above we have used (\ref{eqn:Union}) and the following identity:
\begin{eqnarray*}
&& c^{\mu}_{\mu^+\mu^-} = \sum_{\nu^+, \nu^-}
\frac{\chi_{\mu^+}(\nu^+)\chi_{\mu^-}(\nu^-)\chi_{\mu}(\nu^+\cup\nu^-)}{z_{\nu^+}z_{\nu^-}}.
\end{eqnarray*}
\end{proof}

\begin{rema}
By the same method we also have
\begin{equation}\label{eqn:anotherW}
G^\bu(\lambda;0,p^+,p^3;1)=
 \sum_{\gamma^{+}, \gamma^i, \mu^i \in \cP}
c^{\gamma^+}_{(\gamma^1)^t\gamma^2} q^{(-2\kappa_{\gamma^+}  +\frac{ \kappa_{\gamma^3} }{2})/2}
\cW_{\gamma^+, (\gamma^3)^t}(q)\prod_{i=1}^3
\frac{\chi_{\gamma^i}(\mu^i)}{z_{\mu^i}}p^i_{\mu^i}.
\end{equation}
\end{rema}

\begin{lemm}\label{thm:extra}
We have
\begin{equation}
\exp \biggl(-\sum_{m\geq 1}\frac{(-1)^{m-1}}{m}p^1_m p^3_{2m} \biggr)
=\sum_{\mu\in\cP}\frac{(-1)^{|\mu|-\ell(\mu)} }{z_\mu} p_\mu^1 p_{2\mu}^3
\end{equation}
where 
$2\mu$ is the partition $(2\mu_1, 2\mu_2, \dots, 2\mu_{\ell(\mu)})$.
\end{lemm}

\begin{proof}
Let $(x^i_1, \dots, x_n^i, \dots)$ be formal variables such that
${ \displaystyle p_m^i = \sum_n (x_n^i)^m}$.
By standard series manipulations,
\begin{eqnarray*}
&& \exp \biggl(-\sum_{m \geq 1} \frac{(-1)^{m-1}}{m}p^2_m p^3_{2m}\biggr)
= \exp \biggl(\, \sum_{m\geq 1} \frac{(-1)^{m-1}}{m}
\sum_{n_1, n_3} (x^1_{n_1})^m (x_{n_3}^3)^{2m} \biggr)\\
&& \quad =  \prod_{n_1, n_3} \exp
\biggl(\, \sum_{m \geq 1} \frac{(-1)^{m-1}}{m} (p_{n_1}^1(p_{n_3}^3)^2)^m\biggr) 
=  \prod_{n_1, n_3} (1 + x_{n_1}^1(x_{n_3}^3)^2).
\end{eqnarray*}
Now recall (cf. \cite{Mac}, p. 65, (4.1')):
$$
\prod_{i,j} (1+x_iy_j) = \sum_{\mu\in \cP} \frac{(-1)^{|\mu|-\ell(\mu)}}{z_{\mu}}
p_{\mu}(x) p_{\mu}(y).
$$
Hence we have
\begin{eqnarray*}
 \exp \biggl(\, \sum_{m \geq 1} \frac{(-1)^{m-1}}{m} p^1_mp^3_{2m}\biggr)
& = &  \sum_{\mu\in \cP}  \frac{(-1)^{|\mu|-\ell(\mu)}}{z_{\mu}}
p_{\mu}(x^1) p_{\mu}((x^3)^2) \\
& = &  \sum_{\mu\in \cP}  \frac{(-1)^{|\mu|-\ell(\mu)}}{z_{\mu}}
p_{\mu}(x^1) p_{2\mu}(x^3).
\end{eqnarray*}
\end{proof}

By Lemma \ref{thm:main-extra}, Lemma \ref{thm:main}, and Lemma
\ref{thm:extra}, we have
\begin{eqnarray*}
&& G^\bu(\lambda;p^1,p^2,p^3;1) \\
& = & G^\bu(\lam;0,p^+,p^3;1)
\exp  \bigg(\,\sum_{m\geq 1} \frac{(-1)^{m-1}}{m} p_m^1p_{2m}^3 \bigg) \\
& = & \sum_{\nu^{+}, \nu^i \in \cP}
c^{\nu^+}_{(\nu^1)^t\nu^2} q^{(-2\kappa_{\nu^+}  - \frac{\kappa_{\nu^3}}{2})/2} \cW_{\nu^+, \nu^3}
s^1_{\nu^1} s^2_{\nu^2}  s^3_{(\nu^3)^t}
\sum_{\mu\in \cP}  \frac{(-1)^{|\mu|-\ell(\mu)}}{z_{\mu}}
p_{\mu}(x^1) p_{2\mu}(x^3) \\
& = & \sum_{\nu^{+}, \nu^i \in \cP}
c^{\nu^+}_{(\nu^1)^t\nu^2} q^{(-2\kappa_{\nu^+}  - \frac{\kappa_{\nu^3}}{2})/2} \cW_{\nu^+, \nu^3}(q)
s^1_{\nu^1} s^2_{\nu^2}  s^3_{(\nu^3)^t} \\
&& \cdot \sum_{\mu, \eta^1, \eta^3 \in \cP}  \frac{(-1)^{|\mu|-\ell(\mu)}}{z_{\mu}}
\chi_{\eta^1}(\mu)\chi_{\eta^3}(2\mu) s^1_{\eta^1} s^3_{\eta^3} \\
& = &  \sum_{\nu^+,\nu^1,\nu^3,\eta^1,\eta^3,\mu \in\cP}
c^{\nu^+}_{(\nu^1)^t\nu^2}
c_{\eta^1\nu^1}^{\rho^1}c_{\eta^3(\nu^3)^t}^{\rho^3}
q^{(-2\kappa_{\nu^+}  - \frac{\kappa_{\nu^3}}{2})/2}
\cW_{\nu^+,\nu^3}(q)\frac{\chi_{(\eta^1)^t}(\mu)\chi_{\eta^3}(2\mu)}{z_\mu}
s^1_{\rho^1} s_{\rho^2}^2 s^3_{\rho^3}\\
& = &  \sum_{\nu^+,\nu^1,\nu^3,\eta^1,\eta^3,\mu \in\cP}
c^{\nu^+}_{(\nu^1)^t\nu^2}
c_{(\eta^1)^t\nu^1}^{\rho^1}c_{\eta^3(\nu^3)^t}^{\rho^3}
q^{(-2\kappa_{\nu^+}  - \frac{\kappa_{\nu^3}}{2})/2}
\cW_{\nu^+,\nu^3}(q)\frac{\chi_{\eta^1}(\mu)\chi_{\eta^3}(2\mu)}{z_\mu}
s^1_{\rho^1} s_{\rho^2}^2 s^3_{\rho^3}.
\end{eqnarray*}

By Proposition \ref{thm:GtC},
\begin{eqnarray*}
G^\bu(\lam;\bp;1)
& = & \sum_{\mu^i, \nu^i\in \cP} \tC_{\vnu}(\lam)
q^{(\kappa_{\nu^1}-2\kappa_{\nu^2}-\frac{1}{2}\kappa_{\nu^3})/2}
\prod_{i=1}^3 \frac{\chi_{\nu^i}(\mu^i)}{z_{\mu^i}}p^i_{\mu^i} \\
& = & \sum_{\nu^i\in \cP} \tC_{\vnu}(\lam)
q^{(\kappa_{\nu^1}-2\kappa_{\nu^2}-\frac{1}{2}\kappa_{\nu^3})/2}
\prod_{i=1}^3 s_{\nu^i}(x^i).
\end{eqnarray*}
By comparing coefficients,
\begin{eqnarray*}
&&  \sum
c^{\nu^+}_{(\nu^1)^t\nu^2}
c_{(\eta^1)^t\nu^1}^{\rho^1}c_{\eta^3(\nu^3)^t}^{\rho^3}
q^{(-2\kappa_{\nu^+}  - \frac{\kappa_{\nu^3}}{2})/2} \cW_{\nu^+, \nu^3}(q)
\frac{\chi_{\eta^1}(\mu)\chi_{\eta^3}(2\mu)}{z_\mu}s^1_{\rho^1} s_{\rho^2}^2 s^3_{\rho^3} \\
& = & \sum_{\mu^i, \nu^i\in \cP} \tC_{\vnu}(\lam)
q^{(\kappa_{\nu^1}-2\kappa_{\nu^2}-\frac{1}{2}\kappa_{\nu^3})/2}
\prod_{i=1}^3 s^i_{\nu^i}.
\end{eqnarray*}
Therefore,
$$
\tC_{\vec{\rho}}(\lam)=\tilde{\cW}_{\vec{\rho}}(q)
$$
where $\tilde{\cW}_{\vec{\rho}}(q)$ is defined by \eqref{eqn:tW}.
This completes the proof of Theorem \ref{thm:R5}.

\begin{rema}
By \eqref{eqn:anotherW} one gets a slightly different expression.
\end{rema}

\subsection{Examples of Conjecture \ref{conj:WW}} \label{sec:lowdegree}

Recall that

\medskip

\paragraph{\bf Conjecture \ref{conj:WW}}
{\em Let $\vmu\in\cP^3_+$. Then
$$
\tilde{\cW}_{\vmu}(q)=\cW_\vmu(q),
$$
where $q=e^{\sqrt{-1}\lam}$, and $\cW_\vmu(q)$ is defined by \eqref{eqn:Wthree}}.
\bigskip

We have seen in Section \ref{sec:comb} that Conjecture \ref{conj:WW} holds
when one of the three partitions is empty. When none of the partitions
is empty, A. Klemm has checked by computer that Conjecture \ref{conj:WW} holds
in all the cases where
$$
|\mu^i|\leq 6,\ \ i=1,2,3.
$$
We list some of these cases here.

$$
 \tilde\cW_{(1),(1),(1)}(q)=\cW_{(1),(1),(1)}(q)
=\frac{q^4-q^3+q^2-q+1}{q^{1/2}(q-1)^3}
$$
$$
 \tilde\cW_{(1),(1),(2)}(q)=\cW_{(1),(1),(2)}(q)
=\frac{q^6-q^5+q^3-q+1}{(q^2-1)(q-1)^3}
$$
$$
 \tilde\cW_{(1),(1),(1,1)}(q)=\cW_{(1),(1),(1,1)}(q)
=\frac{q^6-q^5+q^3-q+1}{q(q^2-1)(q-1)^3}
$$
$$
 \tilde\cW_{(1),(1),(3)}(q)=\cW_{(1),(1),(3)}(q)
=\frac{q^{3/2}(q^8-q^7+q^4-q+1)}{(q^3-1)(q^2-1)(q-1)^3}
$$
$$
 \tilde\cW_{(1),(1),(2,1)}(q)=\cW_{(1),(1),(2,1)}(q)
=\frac{q^8-2q^7+3q^6-3q^5+3q^4-3q^3+3q^2-2q+1}{q^{1/2}(q^3-1)(q-1)^4}
$$
$$
 \tilde\cW_{(1),(1),(1,1,1)}(q)=\cW_{(1),(1),(1,1,1)}(q)
=\frac{q^8-q^7+q^4-q+1}{q^{3/2}(q^3-1)(q^2-1)(q-1)^3}
$$
$$
 \tilde\cW_{(1),(2),(2)}(q)=\cW_{(1),(2),(2)}(q)
=\frac{q^{1/2}(q^8-q^7+q^5-q^4+q^3-q+1)}{(q^2-1)^2(q-1)^3}
$$
$$
 \tilde\cW_{(1),(1,1),(2)}(q)=\cW_{(1),(1,1),(2)}(q)
=\frac{q^9-q^8+q^6-q^5+2q^3-q^2-q+1}{q^{3/2}(q^2-1)^2(q-1)^3}
$$
$$
 \tilde\cW_{(1),(2),(1,1)}(q)=\cW_{(1),(2),(1,1)}(q)
=\frac{q^9-q^8-q^7+2q^6-q^4+q^3-q+1}{q^{1/2}(q^2-1)^2(q-1)^3}
$$
$$
 \tilde\cW_{(1),(1,1),(1,1)}(q)=\cW_{(1),(1,1),(1,1)}(q)
=\frac{q^8-q^7+q^5-q^4+q^3-q+1}{q^{3/2}(q^2-1)^2(q-1)^3}
$$
$$
 \tilde\cW_{(1),(1),(4)}(q)= \cW_{(1),(1),(4)}(q)
=\frac{q^4(q^{10}-q^9+q^5-q+1)}{(q^4-1)(q^3-1)(q^2-1)(q-1)^3}
$$
$$
\tilde\cW_{(1),(1),(3,1)}(q)=\cW_{(1),(1),(3,1)}(q)
=\frac{q(q^{10}-2q^9+2q^8-2q^6+3q^5-2q^4+2q^2-2q+1)}{(q^4-1)(q^2-1)(q-1)^4}
$$
$$
\tilde\cW_{(1),(1),(2,2)}(q)=\cW_{(1),(1),(2,2)}(q)
=\frac{q(q^8-2q^6+q^5+q^4+q^3-2q^2+1)}{(q^3-1)(q^2-1)^2(q-1)^3}
$$
$$
\tilde\cW_{(1),(1),(2,1,1)}(q)=\cW_{(1),(1),(2,1,1)}(q)
=\frac{q^{10}-2q^9+2q^8-2q^6+3q^5-2q^4+2q^2-2q+1}{q(q^4-1)(q^2-1)(q-1)^4}
$$
$$
\tilde\cW_{(1),(1),(1,1,1,1)}(q)=\cW_{(1),(1),(1,1,1,1)}(q)
=\frac{q^{10}-q^9+q^5-q+1}{q^2(q^4-1)(q^3-1)(q^2-1)(q-1)^3}
$$
$$
\tilde\cW_{(1),(2),(3)}(q)=\cW_{(1),(2),(3)}(q)
=\frac{q^2(q^{10}-q^9+q^6-q^4+q^3-q+1)}{(q^3-1)(q^2-1)^2(q-1)^3}
$$
$$
\tilde\cW_{(1),(3),(2)}(q)=\cW_{(1),(3),(2)}(q)
=\frac{q^2(q^{10}-q^9+q^7-q^6+q^4-q+1)}{(q^3-1)(q^2-1)^2(q-1)^3}
$$
$$
\tilde\cW_{(1),(2),(2,1)}(q)=\cW_{(1),(2),(2,1)}(q)
=\frac{q^{11}-2q^{10}+2q^9- q^8 +q^7-q^6+q^4-q+1}{(q^3-1)(q^2-1)(q-1)^4}
$$
$$
\tilde\cW_{(1),(2,1),(2)}(q)=\cW_{(1),(2,1),(2)}(q)
=\frac{q^{11}-q^{10}+q^7-q^5+q^4-q^3+2 q^2 -2q+1}{q(q^3-1)(q^2-1)(q-1)^4}
$$
$$
\tilde\cW_{(1),(2),(1,1,1)}(q)=\cW_{(1),(2),(1,1,1)}(q)
=\frac{q^{12}-q^{11}-q^{10}+q^9 +q^8-q^6+q^4-q+1}{q(q^3-1)(q^2-1)^2(q-1)^3}
$$
$$
\tilde\cW_{(1),(1,1,1),(2)}(q)=\cW_{(1),(1,1,1),(2)}(q)
=\frac{q^{12}-q^{11}+q^8-q^6+q^4+q^3 -q^2-q+1}{q^3(q^3-1)(q^2-1)^2(q-1)^3}
$$
$$
\tilde\cW_{(1),(1,1),(3)}(q)=\cW_{(1),(1,1),(3)}(q)
=\frac{q^{12}-q^{11}+q^8-q^6+q^4+q^3-q^2-q+1}{q(q^3-1)(q^2-1)^2(q-1)^3}
$$
$$
\tilde\cW_{(1),(3),(1,1)}(q)=\cW_{(1),(3),(1,1)}(q)
=\frac{q(q^{12}-q^{11}-q^{10}+q^9+q^8-q^6+q^4-q+1)}{(q^3-1)(q^2-1)^2(q-1)^3}
$$
$$
\tilde\cW_{(1),(1,1),(2,1)}(q)=\cW_{(1),(1,1),(2,1)}(q)
=\frac{q^{11}-q^{10}+q^7-q^5+q^4-q^3+2q^2-2q+1}{q^2(q^3-1)(q^2-1)(q-1)^4}
$$
$$
\tilde\cW_{(1),(2,1),(1,1)}(q)=\cW_{(1),(2,1),(2,1)}(q)
=\frac{q^{11}-2q^{10}+2q^9-q^8+q^7-q^6+q^4-q+1}{q(q^3-1)(q^2-1)(q-1)^4}
$$
$$
\tilde\cW_{(1),(1,1),(1,1,1)}(q)=\cW_{(1),(1,1),(1,1,1)}(q)
=\frac{q^{10}-q^9+q^7-q^6+q^4-q+1}{q^2(q^3-1)(q^2-1)^2(q-1)^3}
$$
$$
\tilde\cW_{(1),(1,1,1),(1,1)}(q)=\cW_{(1),(1,1,1),(1,1)}(q)
=\frac{q^{10}-q^9+q^6-q^4+q^3-q+1}{q^2(q^3-1)(q^2-1)^2(q-1)^3}
$$
$$
\tilde\cW_{(2),(2),(2)}(q)=\cW_{(2),(2),(2)}(q)
=\frac{q(q^{10}-3q^8+3q^7+2q^6-5q^5+2q^4+3q^3-3q^2+1)}{(q^2-1)^3(q-1)^3}
$$
$$
\tilde\cW_{(2),(2),(1,1)}(q)=\cW_{(2),(2),(1,1)}(q)
=\frac{q^{12}-q^{11}-q^{10}+2q^9-q^7+q^6-q^5+2q^3-q^2-q+1}{q(q^2-1)^3(q-1)^3}
$$
$$
\tilde\cW_{(2),(1,1),(1,1)}(q)=\cW_{(2),(1,1),(1,1)}(q)
=\frac{q^{12}-q^{11}-q^{10}+2q^9-q^7+q^6-q^5+2q^3-q^2-q+1}{q^2(q^2-1)^3(q-1)^3}
$$
$$
\tilde\cW_{(1,1),(1,1),(1,1)}(q)=\cW_{(1,1),(1,1),(1,1)}(q)
=\frac{q^{10}-3q^8+3q^7+2q^6-5q^5+2q^4+3q^3-3q^2+1}{q^2(q^2-1)^3(q-1)^3}
$$
$$
\cW_{(1),(2),(3,1)}(q)=\tilde{\cW}_{(1),(2),(3,1)}(q)
=\frac{q^{3/2}(q^{13}-2q^{12}+q^{11}+2q^{10}-3q^9+2q^8-2q^6+2q^5-q+1)}
{(q^4-1)(q^2 -1)^2(q-1)^4}
$$
\begin{eqnarray*}
&& \cW_{(1,1),(2,1),(3)}(q)=\tilde{\cW}_{(1,1),(2,1),(3)}(q)
=\bigl(q^{19}-q^{18}-q^{17}+q^{16}+q^{15}-q^{13}+q^{11}-q^{10}\\
&&\quad \quad \quad + q^8 +q^7-q^6 -2q^5 +2q^4+q^2-2q+1\bigr)
\cdot \bigl(q^2(q^3-1)^2(q^2-1)^2(q-1)^4 \bigr)^{-1}
\end{eqnarray*}
\begin{eqnarray*}
&& \cW_{(2),(2),(2,1,1,1)}(q)=\tilde{\cW}_{(2),(2),(2,1,1,1)}(q)\\
&=&\bigl(q^{22} -q^{21} -2q^{20} +3q^{19} +q^{18} -3q^{17}
         +3q^{15} -q^{14} -2q^{13} + q^{12}+q^{11} +q^{10} -2q^9-q^8\\
&& +3q^7-3q^5+q^4+3q^3-2q^2-q+1\bigr)
\cdot \bigl(q^{7/2}(q^5-1)(q^3-1)(q^2-1)^3(q-1)^4\bigr)^{-1}\\
&& \cW_{(1),(2,2),(3,2)}(q)=\tilde{\cW}_{(1),(2,2),(3,2)}(q)\\
&=& \bigl(q^{23}-2q^{22}+q^{21}+q^{20}-q^{19}+q^{18}-2q^{17}+q^{16} 
+q^{15} +q^{13}-3q^{12}+q^{10}+ 2q^9 + q^8\\
&&  -2q^7-2q^6+2q^4+2q^3-2q^2-q+1\bigr)
\cdot \bigl(q(q^4-1)(q^3-1)^2(q^2-1)^3(q-1)^4\bigr)^{-1}
\end{eqnarray*}